\documentclass[10pt, a4paper]{compositio}

\usepackage[english]{babel}

\usepackage{amsmath}
\usepackage{amssymb}
\usepackage{amsthm}
\usepackage{mathabx}
\usepackage{tikz}
\usepackage{tikz-cd}
\usepackage{todonotes}
\usepackage{hyperref}
\usepackage{cleveref}
\usepackage{footmisc}
\usepackage{bbm}
\usepackage{afterpage}
\usepackage{mathrsfs}
\usepackage{sseq}
\usepackage{extarrows}
\usepackage{enumitem}
\usepackage{pdflscape}
\usepackage{graphicx}
\usepackage{caption}

\newcommand{\Baker}{\mathrm{Baker}}
\newcommand{\can}{\mathrm{can}}
\newcommand{\cl}{\mathrm{cl}}
\newcommand{\fet}{\mathrm{f\acute{e}t}}
\newcommand{\et}{\mathrm{\acute{e}t}}
\newcommand{\fin}{\mathrm{fin}}
\DeclareMathOperator{\Fin}{Fin}

\newcommand{\fold}{\mathrm{fold}}
\newcommand{\GHM}{\mathrm{GHM}}
\DeclareMathOperator{\lcm}{lcm}
\newcommand{\open}{\mathrm{open}}
\renewcommand{\top}{\mathrm{top}}

\DeclareMathOperator{\Tw}{Tw}

\newcommand{\BTn}{{{\mathrm{BT}^p_n}}}
\newcommand{\BTtwo}{{{\mathrm{BT}^p_2}}}
\DeclareMathOperator{\Spec}{Spec}
\DeclareMathOperator{\Spf}{Spf}

\DeclareMathOperator{\GL}{GL}

\DeclareMathOperator{\MF}{MF}
\DeclareMathOperator{\mf}{mf}

\newcommand{\ELL}{{{E}{\ell\ell}}}
\DeclareMathOperator{\KU}{KU}
\DeclareMathOperator{\KO}{KO}

\newcommand{\Sph}{\mathbf{S}}

\DeclareMathOperator{\TMF}{TMF}
\DeclareMathOperator{\tmf}{tmf}
\DeclareMathOperator{\Tmf}{Tmf}

\DeclareMathOperator{\CAlg}{CAlg}
\DeclareMathOperator{\CRing}{CRing}

\DeclareMathOperator{\Fun}{Fun}
\DeclareMathOperator{\Isog}{Isog}

\DeclareMathOperator{\Mod}{Mod}
\DeclareMathOperator{\Span}{Span}
\DeclareMathOperator{\Shv}{Shv}
\DeclareMathOperator{\Sp}{Sp}
\newcommand{\Spc}{\mathcal{S}}

\DeclareMathOperator{\colim}{colim}
\newcommand{\id}{\mathrm{id}}

\newcommand{\op}{\mathrm{op}}

\newcommand{\A}{\mathcal{A}}
\newcommand{\B}{\mathcal{B}}
\renewcommand{\b}{\mathfrak{B}}

\newcommand{\E}{\mathbf{E}}
\newcommand{\EE}{\mathscr{E}}

\newcommand{\F}{\mathbf{F}}

\newcommand{\ff}{\mathcal{F}}
\newcommand{\G}{\mathbf{G}}

\renewcommand{\H}{\mathcal{H}}
\renewcommand{\j}{\mathrm{j}}
\newcommand{\J}{\mathrm{J}}
\newcommand{\M}{\mathcal{M}}
\newcommand{\N}{\mathbf{N}}

\renewcommand{\O}{\mathcal{O}}
\newcommand{\OO}{\mathbf{O}}
\newcommand{\Q}{\mathbf{Q}}

\newcommand{\T}{\mathrm{T}}
\newcommand{\TT}{\mathcal{T}}
\renewcommand{\u}{\mathrm{u}}
\newcommand{\U}{\mathrm{U}}
\newcommand{\x}{\mathfrak{X}}
\newcommand{\xx}{\mathsf{X}}
\newcommand{\yy}{\mathsf{Y}}
\newcommand{\Z}{\mathbf{Z}}

\newcommand{\al}{\alpha}
\newcommand{\be}{\beta}

\newcommand{\Ga}{\Gamma}

\usepackage{pict2e}

\makeatletter
\newcommand{\bigcomp}{%
  \DOTSB
  \mathop{\vphantom{\sum}\mathpalette\bigcomp@\relax}%
  \slimits@
}
\newcommand{\bigcomp@}[2]{%
  \begingroup\m@th
  \sbox\z@{$#1\sum$}%
  \setlength{\unitlength}{0.9\dimexpr\ht\z@+\dp\z@}%
  \vcenter{\hbox{%
    \begin{picture}(1,1)
    \bigcomp@linethickness{#1}
    \put(0.5,0.5){\circle{1}}
    \end{picture}%
  }}%
  \endgroup
}
\newcommand{\bigcomp@linethickness}[1]{%
  \linethickness{%
      \ifx#1\displaystyle 2\fontdimen8\textfont\else
      \ifx#1\textstyle 1.65\fontdimen8\textfont\else
      \ifx#1\scriptstyle 1.65\fontdimen8\scriptfont\else
      1.65\fontdimen8\scriptscriptfont\fi\fi\fi 3
  }%
}

\makeatother

\usepackage{mathtools}

\theoremstyle{theorem}\numberwithin{equation}{section}
\newtheorem{theorem}[equation]{Theorem}
\crefname{theorem}{{th}.\!\!}{{ths}.\!\!}
\Crefname{theorem}{{Th}.\!\!}{{Ths}.\!\!}
\newtheorem*{theoremnonum}{Theorem}
\newtheorem{theoremalph}{Theorem}

\crefname{theoremalph}{{th}.\!\!}{{ths}.\!\!}
\Crefname{theoremalph}{{Th}.\!\!}{{Ths}.\!\!}

\Crefname{problem}{{Prb}.\!\!}{{Prbs}.\!\!}
\newtheorem{prop}[equation]{Proposition}
\Crefname{prop}{{Pr}.\!\!}{{Prs}.\!\!}
\newtheorem{lemma}[equation]{Lemma}
\Crefname{lemma}{{Lm}.\!\!}{{Lms}.\!\!}
\newtheorem{cor}[equation]{Corollary}
\Crefname{cor}{{Cor}.\!\!}{{Cors}.\!\!}
\newtheorem{conjecture}[equation]{Conjecture}
\Crefname{conjecture}{{Conj}.\!\!}{{Conjs}.\!\!}
\newtheorem*{conjecturemaeda}{Maeda's conjecture}

\theoremstyle{definition}\numberwithin{equation}{section}
\newtheorem{mydef}[equation]{Definition}
\Crefname{mydef}{{Df}.\!\!}{{Dfs}.\!\!}

\Crefname{recall}{{Rcl}.\!\!}{{Rcls}.\!\!}
\newtheorem{construction}[equation]{Construction}
\Crefname{construction}{{Con}.\!\!}{{Cons}.\!\!}

\Crefname{ass}{{As}.\!\!}{{As}.\!\!}

\Crefname{notation}{{Nt}.\!\!}{{Nts}.\!\!}

\Crefname{situation}{{St}.\!\!}{{Sts}.\!\!}

\theoremstyle{remark}\numberwithin{equation}{section}
\newtheorem{example}[equation]{Example}
\Crefname{example}{{Ex}.\!\!}{{Exs}.\!\!}

\Crefname{nonexample}{{NonEx}.\!\!}{{NonEx}.\!\!}

\Crefname{claim}{{Clm}.\!\!}{{Clms}.\!\!}
\newtheorem{remark}[equation]{Remark}
\Crefname{remark}{{Rmk}.\!\!}{{Rmks}.\!\!}

\Crefname{idea}{{Id}.\!\!}{{Ids}.\!\!}

\Crefname{warn}{{Warn}.\!\!}{{Warns}.\!\!}
\newtheorem{var}[equation]{Variant}
\Crefname{var}{{Var}.\!\!}{{Vars}.\!\!}

\Crefname{figure}{{Fig.}\!\!}{{Figs.}\!\!}
\Crefname{footnote}{{Fn.}\!\!}{{Fn.}\!\!}

\Crefname{part}{{\textsection}\!\!}{{\textsection}\!\!}
\Crefname{chapter}{{\textsection}\!\!}{{\textsection}\!\!}
\Crefname{section}{{\textsection}\!\!}{{\textsection}\!\!}
\Crefname{subsection}{{\textsection}\!\!}{{\textsection}\!\!}
\Crefname{appendix}{{\textsection}\!\!}{{\textsection}\!\!}

\begin{document}
\title{Hecke operators on topological modular forms}
\author{Jack Morgan Davies}
\email{davies@math.uni-bonn.de}
\address{Mathematisches Institut\\Universit\"{a}t Bonn\\Endenicher Alle 60\\53115 Bonn, Germany}

\begin{abstract}
The cohomology theory $\TMF$ of topological modular forms is a derived algebro-geometric interpretation of the ring of complex modular forms from number theory. In this article, we refine the classical Adams operations, Hecke operators, and Atkin--Lehner involutions from endomorphisms of modular forms to stable operators on $\TMF$. Our algebro-geometric formulation of these operators leads to simple proofs of their many remarkable properties and computations. From these properties, we use techniques from homotopy theory to make simple number-theoretic deductions, including a rederivation of some classical congruences of Ramanujan and providing new infinite families of Hecke operators which satisfy Maeda's conjecture.
\end{abstract}


\maketitle

\vspace{0.8cm}

\begin{center}
\includegraphics[width=0.85\textwidth]{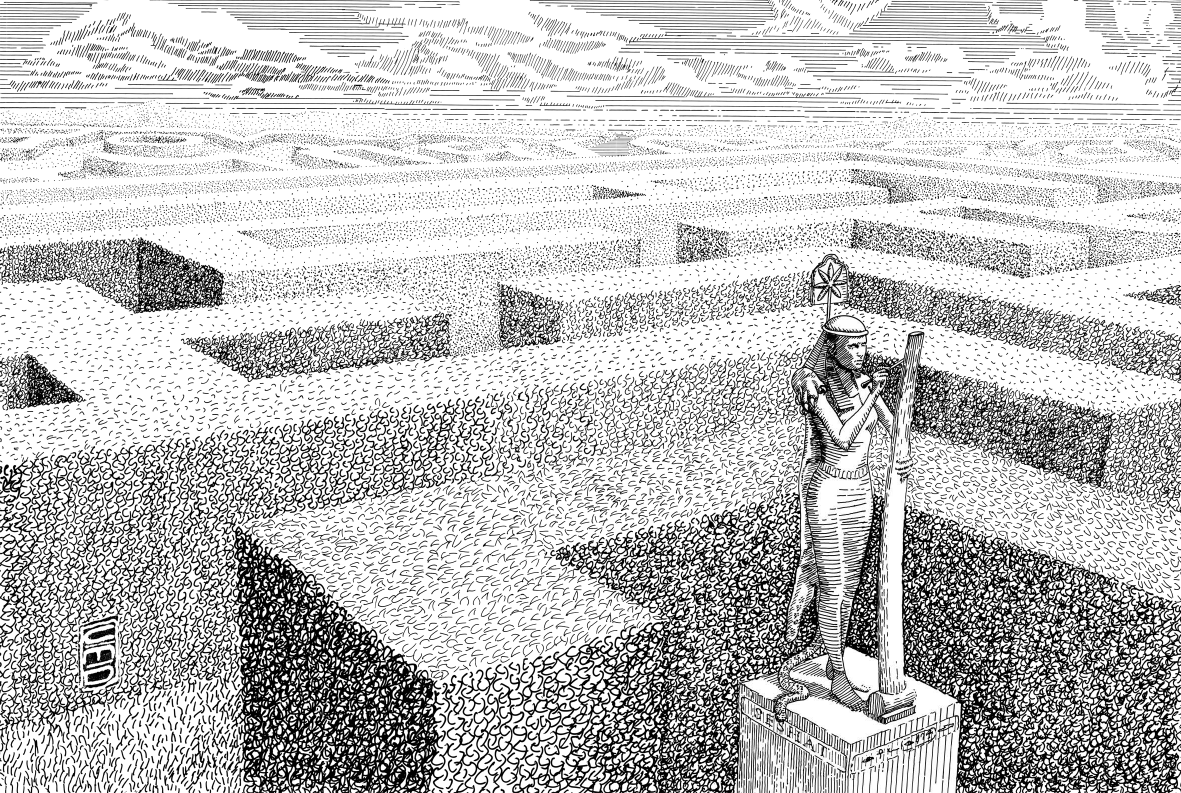}\\
Seshat amidst a hedge maze---Carl Davies 2022
\end{center}

\newpage

\setcounter{tocdepth}{2}
\tableofcontents

%\newpage

\addcontentsline{toc}{section}{Introduction}
\section*{Introduction}

In this article, the connection between number theory and homotopy theory through modular forms is strengthened and refined. We do this by constructing three explicit families of stable operations on the cohomology theory $\TMF$ of \emph{topological modular forms} lifting their classical analogues: \emph{stable Adams operations} $\psi^k$, \emph{stable Atkin--Lehner involutions} $w_Q$, and \emph{stable Hecke operators} $\T_n$. Aside from the conceptual foundation (\Cref{extensiontheoremintro}) giving rise to these operations, which itself may be of independent interest in homotopy theory, and our explicit constructions, this article aims to advertise the elegance, simplicity, and utility of these operations. To this end, we prove that these collections of stable operations:

\begin{itemize}
\item agree with the classical operations from the theory of modular forms when evaluated on the $\TMF$-cohomology of spheres (\Cref{comparisonintro}).
\item mutually commute with one another up to natural (coherent) homotopy and that their compositions can also be succinctly expressed (\Cref{interrationsintro,naturalsplittings,heckeoperatorscompositionintro}).
\item cannot be further refined to stable operations localised \textbf{at} the given prime (\Cref{nonexistenceintro}).
\item have implications in number theory, from rederiving known congruences of Ramanujan to provide further evidence for Maeda's conjecture (\Cref{congruencetheoremintro,primethreemaedaintro}).
\end{itemize}

Further applications of these operations, as well as the methods used below to construct and utilise them, also seem inevitable.

\subsection*{Motivation}

Cohomology theories have long been the main tools at the disposal of algebraic topologists, and any extra structure of these cohomology theories adds to their potency. Such refinements include multiplicative structures, equivariance with respect to a group, or \emph{cohomology operations}. Steenrod powers on mod $p$ cohomology and Adams operations on topological $K$-theory are examples of such operations, and one does not have to look far from the basic properties of these operations to spot some of their many applications. There is an important cohomology theory that has (up until now) not received any cohomology operations: namely, the cohomology theory $\TMF$ of \emph{(periodic) topological modular forms}. \\

The concept of a \emph{modular form} has been fundamental to advancements in number theory for the past century. Originally defined as holomorphic complex-valued functions on the upper half-plane satisfying a modularity property and a growth condition, modern algebraic geometry views modular forms as global sections of a line bundle $\omega^{\otimes d}$ on the \emph{moduli stack of smooth %\footnote{The Goerss--Hopkins--Miller theorem found in \cite[Th.1.2]{bourbakigoerss} or \cite[\textsection 12]{tmfbook} more generally constructs an \'{e}tale sheaf of $\E_\infty$-rings on the moduli stack of \emph{generalised} elliptic curves, also known as the Deligne--Mumford compactification of the moduli stack of smooth elliptic curves. We will only be interested in the latter throughout this article as this lends itself best to Lurie's reinterpretation.}
 elliptic curves} $\M$. There exists a natural homotopical reincarnation of this sheaf $\omega^{\otimes d}$ due to Goerss--Hopkins--Miller and Lurie. To discuss this result we must work in the $\infty$-category of $\E_\infty$-rings---the homotopy-theoretic world whose objects represent cohomology theories with highly structured multiplication.

\begin{theoremnonum}[{[Goerss--Hopkins--Miller, Lurie (see \cite[Th.1.2]{bourbakigoerss} or \cite[Th.7.0.1]{ec2name})]}]\label{ghm}
There is a sheaf $\O^\top_\GHM$ of $\E_\infty$-rings on the moduli stack $\M$ of elliptic curves and an isomorphism of sheaves of abelian groups $\pi_{2d}\O^\top_\GHM\simeq \omega^{\otimes d}$.
\end{theoremnonum}

Just as classical modular forms can be defined as the global sections $\Gamma(\M,\omega^{\otimes d})$, we define $\TMF$ as the global sections $\Gamma(\M,\O^\top_\GHM)$ in the world of $\E_\infty$-rings. Much of the initial fervour behind this construction of $\TMF$ lay in the connection to string manifolds and the Witten genus (\cite{landweberconferencebook,hopkinsfirsttmficm,hopkinssecondtmficm}) and interest in $\TMF$ has only grown in recent years. This is in part due to an elegant construction of Lurie using spectral algebraic geometry \cite{lurieecsurveyname,ec2name}, the emergence of equivariant elliptic cohomology of Lurie and Gepner--Meier (\cite{ec3name,davidandlennart}), and the applications to computations in the stable homotopy groups of spheres due to Goerss--Henn--Mahowald--Rezk and Isaksen--Wang--Xu \cite{ktwolocalsphere,iwxpisph}, just to name a few.\\

As a rich cohomology theory, $\TMF$ begs for cohomology operations. This article provides multiple families of such operations inspired by analogous operations on spaces of modular forms. Let us take Hecke operators as our leading example. One cannot construct such operations on $\TMF$ by manipulating power series, as one might define Hecke operators on classical complex modular forms. Instead, we should look to an algebro-geometric perspective.\\

Fix a prime $\ell$, which we implicitly invert for this example, and write $\M_0(\ell)$ for the \emph{moduli stack of elliptic curves $E$ with a chosen subgroup $H\leq E$ of order $\ell$}. There are two maps of stacks $p,q\colon \M_0(\ell)\to \M$ we would like to consider: the first $p$ sends a pair $(E,H)$ to $E$, and the second sends a pair $(E,H)$ to quotient elliptic curve $E/H$. Following \cite[\textsection 4.5]{conrad}, one can define the \emph{$\ell$th classical Hecke operator} $\T_\ell^\cl$ on $\MF_k$ as the composition
\begin{equation}\label{classicaldefinition}\T^\cl_\ell\colon H^0(\M,\omega^{\otimes d})\xrightarrow{q^\ast} H^0(\M_0(\ell), q^\ast \omega^{\otimes d})\xrightarrow{\xi^{\otimes d}} H^0(\M_0(\ell),p^\ast\omega^{\otimes d})\xrightarrow{\frac{1}{\ell}p_!} H^0(\M,\omega^{\otimes d})\end{equation}
where $\xi\colon q^\ast\omega\to p^\ast\omega$ is the map induced by the quotient $\EE_0(\ell)\to \EE_0(\ell)/\H$ of the universal elliptic curve over $\M_0(\ell)$ by the universal subgroup $\H$ of order $\ell$, and $p_!$ is the transfer map on global sections. One can attempt to retread this path with $\TMF$ in place of $H^0(\M,\omega^{\otimes\ast})$. The morphism $q^\ast$ lifts to $\TMF$ and the realm of higher algebra using the functoriality of $\O^\top_\GHM$. With considerable effort (showing that $\TMF\to \TMF_0(\ell)$ views the target as dualisable over the source), one can also lift $p_!$ to $\TMF$ with the caveat that no functoriality of this lift is guaranteed. More problematic is the lack of a homotopical analogue of the map $\xi^{\otimes k}$. Indeed, the quotient map $\EE_0(\ell)\to \EE_0(\ell)/\H$ does not exist in the small \'{e}tale site of $\M$ upon which the sheaf $\O^\top_\GHM$ is defined---only isomorphisms of elliptic curves are permitted. This is where the results of the article begin.

\subsection*{Main results}

Using techniques from spectral algebraic geometry and the modern homotopy theory surrounding $\infty$-categories of spans, we lift $\O^\top_\GHM$ to a sheaf $\O^\top$ which is functorial with respect to isogenies of elliptic curves of invertible degree equipped with transfer maps which are functorial up to coherent homotopy. To encode these coherences, let us write:

\begin{itemize}
\item $\M_\et$ for the small \'{e}tale site over the moduli stack $\M$ of elliptic curves.
\item $\Isog$ for the site whose objects are those of $\M_\et$ and whose morphisms are isogenies of elliptic curves of invertible degree (\Cref{sitewelike}).
\item $\Span_\fin(\Isog)$ for the $\infty$-category of spans in $\Isog$ where the underlying map of Deligne--Mumford stacks on the right leg is finite and the map of elliptic curves is an isomorphism (\Cref{spectralmackextension}).
\end{itemize}

\begin{theoremalph}\label{extensiontheoremintro}
There are functors $\O^\top$ and $\OO^\top$ such that the diagram of $\infty$-categories
\[\begin{tikzcd}
{\M_\et^\op}\ar[dr, "{\O^\top_\GHM}", swap]\ar[r]	&	{\Isog^\op}\ar[d, "{\O^\top}"]\ar[r]	&	{\Span_\fin(\Isog)}\ar[d, "{\OO^\top}"]	\\
										&	{\CAlg}\ar[r]					&	{\Sp}
\end{tikzcd}\]
 commutes.
\end{theoremalph}

Using this extension of $\O^\top_\GHM$ to $\Isog$, one obtains a morphism of $\E_\infty$-rings
\[\Xi\colon \O^\top_\GHM(q\colon \M_0(\ell)\to \M)\to \O^\top_\GHM(p\colon \M_0(\ell)\to \M)\]
refining $\xi^{\otimes\ast}$, simply by applying $\O^\top$ to the universal isogeny $\EE_0(\ell)\to \EE_0(\ell)/\H$ of order $\ell$. The stable Hecke operators $\T_\ell$ on $\TMF[\frac{1}{\ell}]$ can now be defined as the composite $p_!\circ \Xi\circ q^\ast$ following (\ref{classicaldefinition}), where $p_!$ is a homotopy coherent transfer coming from the extension of $\O^\top_\GHM$ to $\Span_\fin(\Isog)$. Using this extra functoriality provided by \Cref{extensiontheoremintro}, we define \emph{stable Adams operations} $\psi^k$, \emph{stable Atkin--Lehner involutions} $w_Q$, and \emph{stable Hecke operators} $\T_n$:

\[\psi^k\colon \TMF[\frac{1}{k}]\to \TMF[\frac{1}{k}]\qquad\qquad w_Q\colon \TMF_0(N)\to \TMF_0(N)\]
\[\T_n\colon \TMF[\frac{1}{n}]\to \TMF[\frac{1}{n}]\]

More generally, we define operations on various topological modular forms with level structures akin to \emph{double coset operators}. Given a triple of congruence subgroups $\Ga_1,\Ga_2\leq \Ga$, an isomorphism $\phi\colon \Ga_1\simeq \Ga_2$, and an isogeny $\Phi\colon \EE_{\Ga_2}\to \phi^\ast\EE_{\Ga_1}$ of invertible degree, we obtain a stable operation
\[\TMF(\Ga)\xrightarrow{p^\ast_{\Ga_2}}\TMF(\Ga_2)\xrightarrow{(\phi,\Phi)^\ast}\TMF(\Ga_1)\xrightarrow{p_!^{\Ga_1}} \TMF(\Ga)\]
which specialises to the other examples above.\\

Aside from the abstract homotopy theory behind \Cref{extensiontheoremintro} and defining these families of stable operations, the main goal of this article is to advertise these operations by highlighting their formal properties and computational aspects. First, the formalism of \Cref{extensiontheoremintro} and its proof immediately imply that our stable operations agree with their classical counterparts.

\begin{theoremalph}\label{comparisonintro}
Let $d,k,N,n$ be integers, with $N,n$ both positive, and $Q$ a divisor of $N$ such that $\gcd(Q,N/Q)=1$. Then diagrams of abelian groups
\[\begin{tikzcd}
{\pi_{2d}\TMF[\frac{1}{k}]}\ar[r, "{\psi^k}"]\ar[d, "{e}"]	&	{\pi_{2d}\TMF[\frac{1}{k}]}\ar[d, "{e}"]	\\
{\MF_d^{\Z[\frac{1}{k}]}}\ar[r, "{\cdot k^d}"]			&	{\MF_d^{\Z[\frac{1}{k}]}}
\end{tikzcd}\qquad
\begin{tikzcd}
{\pi_{2d}\TMF_0(N)}\ar[r, "{w_Q}"]\ar[d, "{e}"]	&	{\pi_{2d}\TMF_0(N)}\ar[d, "{e}"]	\\
{\MF_0(N)_d}\ar[r, "{w_Q^\cl}"]				&	{\MF_0(N)_d}
\end{tikzcd}\]
\[\begin{tikzcd}
{\pi_{2d}\TMF[\frac{1}{n}]}\ar[r, "{\T_n}"]\ar[d, "{e}"]	&	{\pi_{2d}\TMF[\frac{1}{n}]}\ar[d, "{e}"]	\\
{\MF_d^{\Z[\frac{1}{n}]}}\ar[r, "{n\T_n^\cl}"]		&	{\MF_d^{\Z[\frac{1}{n}]}}
\end{tikzcd}\]
all commute, where $e\colon \pi_{2d}\TMF\to H^0(\M,\omega^{\otimes d})=\MF_d$ is the natural edge map.
\end{theoremalph}

The factor of $n$ popping up in the commutative diagram above involving Hecke operators is only a superficial difference between $\T_n$ and $\T_n^\cl$, as $n$ is always invertible in this case. One should be aware of it during explicit computations though.\\

We have not yet been able to explicitly compare our stable Hecke operators with those defined by Baker \cite{bakerhecketwo} on $\TMF[\frac{1}{6}]$, although the above theorem says these operators agree on coefficient rings up to the unit $n$.\\

Next, we provide the following series of homotopies comparing the various compositions of these operators.

\begin{theoremalph}\label{interrationsintro}
Let $k,\ell$ be integers, $m,n,N$ be positive integers, and $Q,R$ be positive divisors of $N$ such that $\gcd(Q,N/Q)=1$ and $\gcd(R,N/R)=1$. Then there exist natural homotopies 
\[\psi^1\simeq \psi^{-1}\simeq \id \qquad \psi^k\circ \psi^\ell\simeq \psi^{k\ell}\simeq \psi^\ell\circ \psi^k \qquad w_Q\circ \psi^k\simeq \psi^k\circ w_Q\]
\[w_1\simeq \id\qquad w_Q\circ w_Q\simeq \psi^Q\qquad w_Q\circ w_R\simeq w_R\circ w_Q\]
of morphisms of $\E_\infty$-rings. If $\gcd(N,n)=1$ and $\gcd(m,n)=1$, there also exist natural homotopies
\[\T_1\simeq \id\qquad \T_n\circ \psi^k\simeq \psi^k\circ \T_n\qquad \T_n\circ w_Q\simeq w_Q\circ \T_n\qquad \T_n\circ \T_m\simeq \T_m\circ \T_n\]
 of morphisms of spectra.
\end{theoremalph}

Constructing homotopies between morphisms of $\E_\infty$-rings or spectra, especially spectra such as $\TMF$, is highly nontrivial. However, our definitions of these operations using \Cref{extensiontheoremintro} lead to relatively simple proofs of the above facts. The adjective \emph{natural} appearing above is also not used lightly: the particular choice of these homotopies are 1-homotopies in a coherent collection of higher homotopies; this is discussed in more detail in \Cref{coherenthomotopiesremark}. As an explicit example of this, it follows from \Cref{interrationsintro} that the $C_2$-equivariant versions of $\TMF_1(N)$ (as used in \cite{lennartandhill}) obtain stable Adams operations and Hecke operators.\\

In fact, the homotopies above can be controlled so precisely that they look as close to equalities as one can hope in higher category theory. We illustrate this idea with prime decomposition formulae for stable Adams operations and Hecke operators.

\begin{theoremalph}\label{naturalsplittings}
Let $k,n$ be integers with $n\geq 1$. Then we have the natural homotopies
\[\psi^k\simeq \bigcomp_{\substack{p|k \\ p\,\mathrm{prime}}}\psi^{p^{\nu_p(k)}}\qquad\qquad \T_n\simeq \bigcomp_{\substack{p|n \\ p\,\mathrm{prime}}}\T_{p^{\nu_p(n)}}\]
the first being a homotopy between endomorphisms the $\E_\infty$-ring $\TMF[\frac{1}{k}]$ and the second of the spectrum $\TMF[\frac{1}{n}]$, where $\nu_p(-)$ indicates $p$-adic valuation.
\end{theoremalph}

The symbol $\bigcomp$ indicates that all possible orders of compositions are homotopic up to higher coherent homotopy. For us, this natural homotopy between $\T_n\circ \T_m$ and $\T_m\circ \T_n$ for coprime $m$ and $n$ is evidence that a (speculative) spectral Hecke algebra should have an $\E_\infty$-structure---this is work in progress.\\

We also have an expression for $\T_n\circ \T_m$ for $\gcd(m,n)\neq 1$ so long as we further invert $\gcd(6,\phi(mn))$, where $\phi(-)$ is \emph{Euler's totient function}. The following formula is inspired by Baker's definition of stable Hecke operators on elliptic cohomology \cite{bakerhecketwo} and the classical formula \cite[\textsection5.3]{diamondshur}.

\begin{theoremalph}\label{heckeoperatorscompositionintro}
Let $m,n$ be two positive integers and set $\phi=\gcd(6,\phi(mn))$. Then there is a homotopy between morphisms of spectra
\[\T_n\circ \T_m\simeq \sum_{d|m,n}d\psi^d \T_{\frac{mn}{d^2}}\colon \TMF[\frac{1}{mn\phi}]\to \TMF[\frac{1}{mn\phi}]\]
where the sum above ranges over those positive integers $d$ which divide both $m$ and $n$.
\end{theoremalph}

The last property we prove about these operators is that they are as integral as possible. In more detail, we also show the necessity of inverting $k$ (resp.\ $n$) in our definition of the stable Adams operation $\psi^k$ (resp.\ the Hecke operator $\T_n$).

\begin{theoremalph}\label{nonexistenceintro}
Let $p$ be a prime and $e\geq 1$. There does not exist a stable operation $\psi^{p^e}$ on $\TMF_p$ which agrees with our stable operation $\psi^{p^e}$ on $\TMF[\frac{1}{p}]$ upon taking homotopy groups over $\Q_p$. There does not exist a stable operation $\T_{p^e}$ on $\TMF_p$ which agrees with our stable operation $\T_{p^e}$ on $\TMF[\frac{1}{p}]$ upon taking homotopy groups over $\Q_p$ as long as either $p=2,3$ or we further demand $\T_{p^e}$ to be $\Z_p^\times$-equivariant with respect to the $p$-adic Adams operations.
\end{theoremalph}

Finally, we explore some utility of these stable Hecke operators with two applications in number theory. The first is a rederivation of some classical congruences due to Ramanujan. Recall the \emph{divisor sum function} $\sigma(n)=\sum_{d|n}d$ and the \emph{Ramanujan $\tau$-function} $\tau(n)$, defined as the coefficients of the $q$-expansion of the \emph{discriminant modular form} $\Delta$ of weight $12$.

\begin{theoremalph}[{[Congruences of Ramanujan]}]\label{congruencetheoremintro}
If $n$ is odd, then $n\tau(n)\equiv \sigma(n)$ modulo $8$. If $n$ is not divisible by $3$, then $n\tau(n)\equiv \sigma(n)$ modulo $3$.
\end{theoremalph}

These congruences are usually expressed as $\tau(n)\equiv_8 \sigma(n)$ for odd $n$, and $\tau(n)\equiv_3 \sigma(n)$ for $n$ not divisible by $3$. These statements are equivalent; see \cite[p.\ 142]{mythesis}. Curiously, the congruence $n\tau(n)\equiv_{16}\sigma(n)$ holds at least for odd $n\leq 10^6$, but classically we have the non-congruence $\tau(5)=4830\equiv_{16} 14 \not\equiv_{16} 6 =\sigma(5)$.\\

Of course, these particular congruences are also simple consequences of the theory of modular forms, however, the fact we can prove these statements using the torsion in $\pi_\ast\TMF$, rather than manipulating $q$-expansions, is notable. Moreover, these results easily generalise using the periodicity of $\pi_\ast\TMF$ which leads us to our second application: \emph{Maeda's conjecture}. Write $S_d$ for the space of \emph{cusp forms} inside the group of rational weight $d$ holomorphic modular forms, and $\T_{n,d}(X)$ for the characteristic polynomial of the operation $\T_n^\cl$ acting upon the rational vector space $S_d$.

\begin{conjecturemaeda}
Let $d,n$ be positive integers, where $d$ is even and greater than $24$. Then $\T_{n,d}(X)$ is irreducible over $\Q$ and has Galois group the full symmetric group on $D$ letters, where $D=\dim_\Q S_d$.
\end{conjecturemaeda}

Empirical evidence exists for this conjecture (\cite{expevimaedaconjecture}), and the congruences appearing out of our calculations of $\T_n$ on $\pi_\ast\TMF$ add to this empirical evidence.

\begin{theoremalph}\label{primethreemaedaintro}
Let $d,n\geq 2$ be two coprime integers with $n$ not divisible by $3$ satisfying the following conditions:
\begin{enumerate}
\item $d\leq 1,000$ and for all $1\leq i\leq d-1$, the coefficient of $q^d$ in the $q$-expansion of $\Delta^i$ is divisible by $3$.
\item For each prime factor $p$ of $n$ with exponent $e$, if $p\equiv_3 1$ then $e\equiv_6 0,1,3,4$, and if $p\equiv_3 2$ then $e$ is even. 
\end{enumerate}
Then $\T_{dn,12d}(X)$ is irreducible over $\Q$ and has Galois group $\Sigma_{d}$, so Maeda's conjecture holds for $\T_{dn}^\cl$ acting on $S_{12d}$.
\end{theoremalph}

For example, this produces an infinite collection of Hecke operators $\T_{2p}^\cl$ acting on $S_{24}$ which satisfy Maeda's conjecture, as long as $p$ is congruent to $1$ modulo $3$. This is the first infinite collection of classical Hecke operators $\T_n^\cl$ satisfying Maeda's conjecture in a fixed degree where $n$ is composite. We also have a result at the prime $2$, although it is more complicated to state; see \Cref{primetwomaedaintext}. In general, this article highlights the following guiding principle:

\begin{center}\emph{
The more sophisticated the homotopy theory surrounding $\TMF$ becomes,\\
the stronger the number-theoretic results we obtain.
}\end{center}

For example, a more in-depth study into stable Adams operations, Hecke operators, and Atkin--Lehner involutions on $\TMF$ (and its cousins $\TMF(\Ga)$, $\Tmf$, and $Q(N)$) could yield further evidence for Maeda's conjecture at primes $p\geq 5$ (\Cref{higherprimesmaedaremark}), and hopefully lead to other applications as well.

\section*{Related work}
Although the families of stable operations in this article are the first of their kind on $\TMF$ itself, there has been a long history of cohomology operations on \emph{elliptic cohomology}. A direct precursor to this article is Baker's work on stable Hecke operators $\T_n^\Baker$ on the elliptic cohomology theory $\ELL[\frac{1}{n}]=\TMF[\frac{1}{6n}]$; see \cite{bakerhecketwo,bakeropsandcoops,bakeretwoterm}. Many of our results here are inspired and guided by his work. Recent work of Candelori--Salch \cite{candelorisalchone} uses Baker's stable Hecke operators on the elliptic cohomology of various spaces and draws number-theoretic conclusions about these groups---as stated in \cite[Rmk.3.9]{candelorisalchone}, their study of stable Hecke operators is somewhat orthogonal to our work here. \emph{Unstable} Hecke operators have also appeared in the work of Ando \cite{andopo} and Rezk \cite{isogpoht}, and both authors highlight the importance of using isogenies to construct power operations. Unstable power operations and Hecke operators also appear in work by Ganter \cite{gantertatektheory} on a form of elliptic cohomology known as \emph{Tate $K$-theory}. In their work on topological automorphic forms, Behrens--Lawson \cite[\textsection11]{taf} construct stable Hecke operators associated with double coset elements similar to our \Cref{doublecosetoperationdef}. Considering that our approach and that of \cite{taf} both depend on \emph{Lurie's theorem} (\Cref{luriestheoremintext}), the stable operations constructed in this article can be viewed as an integral and homotopy coherent refinement of those of Behrens--Lawson at height two. We have also further explored Adams operations on $\TMF$ and $\Tmf$ in \cite{luriestheorem,adamsontmf,realspectra}.

\section*{Outline}

In \Cref{extensionsection}, we prove \Cref{extensiontheoremintro}. In \Cref{luriestheoremsection}, we discuss \emph{Lurie's theorem} (\Cref{luriestheoremintext}) from spectral algebraic geometry and its relevance to $\TMF$ (\ref{diagramfromluriestheorem}). In \Cref{extensionofisgosection}, we prove the half of \Cref{extensiontheoremintro} extending $\O^\top_\GHM$ to be functorial with respect to isogenies of invertible degree (\Cref{isogextension}). In \Cref{mackeyfunctorsection}, we show $\O^\top_\GHM$ can be further extended to a \emph{spectral Mackey functor}, meaning $\O^\top_\GHM$ can be equipped with homotopy coherent transfers (\Cref{spectralmackextension}), the second half of \Cref{extensiontheoremintro}. Throughout \Cref{extensionofisgosection,mackeyfunctorsection}, we emphasise the computability of these extensions on homotopy groups; see \Cref{pullbacketwocompiarson,classicaltransfers}.\\

In \Cref{zooofoperations}, we construct our three families of stable operations on $\TMF$ and prove many of their basic properties, including \Cref{interrationsintro,heckeoperatorscompositionintro,naturalsplittings,comparisonintro}. In \Cref{definitionssection}, these families of stable operations are constructed (\Cref{adamsdefinition,stableheckedefinition,stableaktinlehnerdefinition}) and we show their compatibility with the classical operations (\Cref{comparisonintro})---in both cases using the foundations of \Cref{extensionsection}. In \Cref{interrelationssection}, the interrelations (\Cref{interrationsintro}) between these operations are proven. To study the composition of stable Hecke operators, we spend \Cref{decompiingstackssection} carefully studying the stacks which define them. In \Cref{heckecompsection}, we prove \Cref{heckeoperatorscompositionintro}, and in \Cref{pprimaryfactorisationsection} we prove \Cref{naturalsplittings}---both refining the methods of \Cref{interrelationssection}. Finally, in \Cref{nonexistencesection} we prove (\Cref{nonexistenceintro}) that in many cases it is necessary to invert an integer when constructing stable operations on $\TMF$ (\Cref{nonexistenceofadams,nonexistenceofhecke,heckenonexstiencecomplicated}).\\

In \Cref{applicationssection}, we explore some basic consequences of the above families of stable operations on $\TMF$, including proofs of \Cref{congruencetheoremintro,primethreemaedaintro}. In \Cref{congrunecessection}, we use a basic computation of the stable Hecke operators on the homotopy groups of $\TMF$ to deduce a collection of congruences of classical modular forms (\Cref{congruencetheoremintext}). In \Cref{maedasconjecturesubsection}, these congruences are combined with theory of modular forms to deduce some extra cases of Maeda's conjecture (\Cref{primetthreemaedaintext,primetwomaedaintext}).\\

In \Cref{combinatorialappendix}, we provide a small appendix of combinatorial results used elsewhere.

\section*{Notation}

All of our mathematics is taking place in the language of $\infty$-categories, and much of our notation follows that of \cite{httname,haname,ec2name}; see \cite{luriestheorem} for a more detailed summary.\\

Our algebraic geometry will take place in the $\infty$-category of presheaves of discrete commutative rings $\Fun(\CRing,\Spc)$, which contains the classical $2$-category of functors of presheaves of rings with values in groupoids as a full $\infty$-subcategory. In particular, our stacks will be viewed primarily as functors $X\colon \CRing\to \Spc$ and morphisms of stacks will be natural transformations. When discussing morphisms of stacks, the discrete ring $R$ will often be suppressed, and we will also occasionally replace $R$ with a more general scheme $S$, for simplicity of notation.\\

Write $\M$ for the \emph{moduli stack of elliptic curves}, which is Deligne--Mumford and smooth of finite presentation over $\Spec \Z$; see \cite{delignemumford}. Write $\widehat{\M}$ for the product of $\M$ with $\Spf \Z_p$ for an ambient prime $p$. We will also write $\MF^\Z_d$ for the groups of weight $d$ \emph{meromorphic modular forms} over $\Z$, defined as the sheaf cohomology group $H^0(\M,\omega^{\otimes d})$; see \cite[\textsection12]{diamondim}. If $R$ is flat over $\Z$, then we will write $\MF^R_d=\MF^\Z_d\otimes R$---in this article, $R$ will usually be a localisation or completion of $\Z$.\\

We also assume the reader has been exposed to the definition of $\TMF$ found in \cite{handbooktmf}, \cite{tmfbook}, or \cite{ec2name}. In particular, we will often use the fact that the natural map of $\E_\infty$-rings
\[\TMF[\frac{1}{n}]=\O^\top(\M)[\frac{1}{n}]\to \O^\top(\M_{\Z[\frac{1}{n}]})\]
is an equivalence; this follows by setting $R=\Z[\frac{1}{n}]$ in \cite[Lm.8.1]{meieroroztmfo7}. Moreover, the \emph{descent spectral sequence} (DSS) for sections of $\O^\top$ will be used multiple times, as well as the identification of the DSS for $\TMF$ with its \emph{Adams--Novikov spectral sequence}; see \cite[Cor.5.3]{tmfhomology} for $\tmf$-case. These spectral sequences will appear in more detail in \Cref{pullbacketwocompiarson}, the proof of \Cref{comparisonintro}, \Cref{nonexistencesection}, and \Cref{applicationssection}.

\section*{Acknowledgements}
This article was born out of my PhD thesis at Utrecht University, and I have to thank my supervisor Lennart Meier for suggesting this topic and for his persistent help---from Bonn to Utrecht (there and back again). Also thank you to the members of my reading committee, Mark Behrens, David Gepner, Paul Goerss, Gijs Heuts, and Gerd Laures for their feedback, corrections, and encouragement. Thank you to Miguel Barerro, Venkata Sai Narayana Bavisetty, Dominic Culver, Mike Hill, Tommy Lundemo, Theresa Rauch, Stefan Schwede, Yuqing Shi, and Liz Tatum for your enlightening conversations and suggestions. Finally, I would like to earnestly thank the anonymous referee for their careful reading and suggestions which greatly improved the quality of this article.

%%%%%%%%%%%%%%%%%%%%%%%%%%%%%%%%%%%%%%%%%%%%%%%%%%%%%%%%%%
%%%%%%%%%%%%%%%%%%%%%%%%%%%%%%%%%%%%%%%%%%%%%%%%%%%%%%%%%%
%%%%%%%%%%%%%%%%%%%%%%%%%%%%%%%%%%%%%%%%%%%%%%%%%%%%%%%%%%

\section{Extending the functorality of $\O^\top$}\label{extensionsection}

The cohomology theory $\TMF$ is defined as the global sections of a certain sheaf $\O^\top_\GHM$ of $\E_\infty$-rings on the small \'{e}tale site of the moduli stack $\M$ of elliptic curves. To define our desired operations and prove \Cref{extensiontheoremintro}, we will use a powerful theorem due to Jacob Lurie known as \emph{Lurie's theorem}, which constructs sheaves of $\E_\infty$-rings over moduli stacks of $p$-divisible groups. As the $p$-divisible group associated with an elliptic curve $E$ has many more automorphisms than there are automorphisms of $E$, the sections of a sheaf on a moduli stack of $p$-divisible groups have many more symmetries than $\O^\top_\GHM$.

%%%%%%%%%%%%%%%%%%%%%%%%%%%%%%%%%%%%%%%%%%%%%%%%%%%%%%%%%%

\subsection{Lurie's theorem and topological modular forms}\label{luriestheoremsection}

Fix a prime $p$ and a positive integer $n$. Recall \cite{tatepdiv} that a \emph{$p$-divisible group} $\G$ over a scheme (or Deligne--Mumford stack) $S$ is a collection of finite flat $S$-group schemes $\G_m$ for every positive integer $m$, and closed immersions $\G_m\to \G_{m+1}$ which witness the equality $\G_m=\G_{m+1}[p^m]$. If each $\G_m$ has degree $p^{mn}$ over $S$, we say $\G$ has \emph{height} $n$. A prototypical example of a $p$-divisible group is the collection of $p$-power torsion of an elliptic curve $E$, usually denoted by $E[p^\infty]$, which has height $2$. Lurie's theorem states the existence of a sheaf of $\E_\infty$-rings on a particular site of $p$-divisible groups and that this sheaf has properties similar to $\O^\top_\GHM$. The original statement can be found in \cite[Th.8.1.4]{taf}, but we will be following the notation and language of \cite{luriestheorem} where a proof first appeared; see \cite[\textsection1.1]{luriestheorem} for a more general summary.\\

Write $\M_\BTn$ for the moduli stack of $p$-divisible groups of fixed height $n$, $\widehat{\M}_\BTn$ for the product of $\M_\BTn$ with $\Spf \Z_p$, and $C$ for the site of formal Deligne--Mumford stacks $\x$ over $\widehat{\M}_\BTn$ such that:
\begin{itemize}
\item the structure map $\x\to \widehat{\M}_\BTn$ is formally \'{e}tale, and
\item the formal Deligne--Mumford stack $\x$ is of finite presentation over $\Spf \Z_p$.
\end{itemize}
We give $C$ the \'{e}tale topology through the forgetful functor to formal Deligne--Mumford stacks. The following is a specialised version of \cite[Th.1.6]{luriestheorem} using \cite[Pr.1.8]{luriestheorem}; see \emph{ibid} for the definitions of any unknown words below.

\begin{theorem}\label{luriestheoremintext}
There is an \'{e}tale hypersheaf of $\E_\infty$-rings $\O^\top_\BTn$ on $C$ such that when evaluated on affine objects $\G\colon \Spf R\to \widehat{\M}_\BTn$ in $C$, the $\E_\infty$-ring $A=\O^\top_\BTn(R)$ has the following properties:
\begin{enumerate}
\item $A$ is weakly $2$-periodic, meaning the homotopy group $\pi_2 A$ is a projective $\pi_0 A$-module of rank one and for every integer $d$ the canonical map of $\pi_0 A$-modules
\[\pi_2 A\underset{\pi_0 A}{\otimes} \pi_n A\to \pi_{n+2} A\]
is an isomorphism.
\item There is a chosen natural (in $\Spf R$) isomorphism $\pi_0 A\simeq R$.
\item The groups $\pi_k A$ vanish for odd $k$ and otherwise there is a natural isomorphism of $R$-modules $\pi_{2}A\simeq \omega_{\G}(R)$ where $\omega_\G$ is the \emph{dualising line} of the formal group $\G^\circ$ associated with $\G$, defined by the formula $I/I^2$ where $I$ is the kernel of the augmentation $\O_{\G^\circ}(R)\to R$ determined by the identity element.
\item There is a chosen natural (in $\Spf R$) isomorphism of formal groups between $\G^\circ$ and the \emph{classical Quillen formal group} of $A$.
\end{enumerate}
\end{theorem}

By \cite[Th.5.17]{luriestheorem} (which is essentially \cite[Th.7.0.1]{ec2name}), the $p$-completion of the Goerss--Hopkins--Miller sheaf $\O^\top_\GHM$ (of \cite{bourbakigoerss} or \cite[\textsection 7]{ec2name}) factors through $C$ (where $n=2$). In other words, the diagram of $\infty$-categories
\begin{equation}\label{diagramfromluriestheorem}\begin{tikzcd}
{\M_\et^\op}\ar[r, "{[p^\infty]}"]\ar[d, "{\O^\top_\GHM}", swap]	&	{C^\op}\ar[d, "{\O^\top_\BTtwo}"]	\\
{\CAlg}\ar[r, "{(-)^\wedge_p}"]							&	{\CAlg}
\end{tikzcd}\end{equation}
commutes up to (non-canonical) homotopy, where the upper horizontal functor indicates base-change over $\Spf \Z_p$ followed by the associated $p$-divisible group functor. Intuitively, we need to $p$-complete $\O^\top_\GHM$, or use the original $p$-complete construction of this sheaf, to obtain this comparison with $\O^\top_\BTtwo$, as the later is based on $p$-divisible groups for a fixed prime. Technically, this is because (\ref{diagramfromluriestheorem}) relies on the \emph{Serre--Tate} theorem, which states that the map of stacks $\M\to \M_\BTtwo$ is formally \'{e}tale \textbf{after} base-change over $\Spf \Z_p$; see \cite[Ex.2.7]{luriestheorem} for a counter-example if we do not apply this base-change.\\

The goal of the next section is to show (\ref{diagramfromluriestheorem}) can be further refined in the integral setting and can be interpreted \textbf{without} reference to $p$-divisible groups.

%%%%%%%%%%%%%%%%%%%%%%%%%%%%%%%%%%%%%%%%%%%%%%%%%%%%%%%%%%

\subsection{Extension to isogenies of invertible degree}\label{extensionofisgosection}

The functor sending an elliptic curve to its associated $p$-divisible group can be used to detect properties in morphisms between elliptic curves.

\begin{prop}\label{plocalpdivisiblegroupcheck}
Let $\phi\colon E\to E'$ be a morphism of elliptic curves over a formal Deligne--Mumford stack $\x$ over $\Spf \Z_p$. Then $\phi$ induces an isomorphism between the associated $p$-divisible groups if and only if $\phi$ is an isogeny of degree prime to $p$.
\end{prop}

Recall that an isogeny of elliptic curves is a nonzero homomorphism of elliptic curves---by \cite[Th.2.4.2]{km}, a nonzero homomorphism of elliptic curves finite locally free. The degree of an isogeny is this locally constant degree.

\begin{proof}
It suffices to work locally on $\x$, so we are reduced to the case where $\x=\Spf R$ for some adic ring $R$ with adic topology generated by an ideal $I$ containing $p$. Writing $\Spf R=\colim \Spec R/I^n$, we are further reduced to the case of an affine scheme $X$. There is a short exact sequence of abelian \'{e}tale sheaves on $X$
\[0\to K_p\to E[p^\infty]\xrightarrow{\phi[p^\infty]} E'[p^\infty]\to 0\]
where $K_p$ is the maximal $p$-power component of the finite flat group scheme $K=\ker(E\to E')$. It is now clear that $\phi[p^\infty]$ is an equivalence if and only if the degree of $\phi$ is prime to $p$.
\end{proof}

There is also an integral version of this statement.

\begin{cor}\label{integralpdivisiblegroupcheck}
Let $\phi\colon E\to E'$ be a morphism of elliptic curves over a Deligne--Mumford stack $\xx$. Then $\phi$ induces an isomorphism on the associated $p$-divisible groups after the base change to $\Spf \Z_p$ for every prime $p$ if and only if $\phi$ is an isogeny of invertible degree.
\end{cor}

\begin{proof}
Let us assume $\xx$ is connected. If $p$ is invertible on $\xx$, then $\xx\times \Spf\Z_p$ vanishes and there is nothing to prove. Otherwise, if $p$ is not invertible on $\xx$, we base change over $\Spf \Z_p$ and apply \Cref{plocalpdivisiblegroupcheck}. If $\xx$ is not connected, we apply the above argument to each connected component.
\end{proof}

Inspired by the previous two statements and \Cref{luriestheoremsection}, we define the following categories.

\begin{mydef}\label{sitewelike}
Let $\Isog$ be the following $2$-category:
\begin{itemize}
\item Objects are those in the small \'{e}tale site of $\M$, so pairs $(\xx,E)$ of a Deligne--Mumford stack $\xx$ and an elliptic curve $E$ over $\xx$ such that the defining morphism $\xx\to \M$ is \'{e}tale.
\item 1-morphisms are pairs $(f,\phi)\colon (\xx,E)\to (\xx',E')$ where $f\colon \xx\to \xx'$ is a morphism of Deligne--Mumford stacks (\textbf{not} necessarily over $\M$) and $\phi\colon E\to f^\ast E'$ is an isogeny of elliptic curves over $\xx$ of invertible degree.
\item 2-morphisms $\al\colon (f,\phi)\to (g,\psi)$ are isomorphisms $\al\colon f^\ast E'\to g^\ast E'$ of elliptic curves over $\xx$ such that $\al\circ \phi=\psi$.
\end{itemize}
For an ambient prime $p$, define the $2$-category $\widehat{\Isog}$ by modifying the definition above, by replacing $\M$ with $\widehat{\M}:=\M\times \Spf \Z_p$, Deligne--Mumford stacks with formal Deligne--Mumford stacks, and isogenies of invertible degree with isogenies of degree prime to $p$. Give $\Isog$ (resp.\ $\widehat{\Isog}$) the \'{e}tale topology through the forgetful functor to Deligne--Mumford stacks (resp.\ formal Deligne--Mumford stacks)---this is easily checked to be a morphism of sites; see \cite[Rmk.6.1.5]{mythesis}.
\end{mydef}

Notice that the wide $2$-subcategory of $\Isog$ whose 1-morphisms are those pairs $(f,\phi)$ where $\phi$ is an isomorphism is precisely the usual small \'{e}tale site of $\M$. The same also goes for $\widehat{\Isog}$ and the small \'{e}tale site of $\widehat{\M}$.

\begin{cor}\label{pcompletetheorema}
Fix a prime $p$. There is an \'{e}tale hypersheaf of $\E_\infty$-rings $\O^\top_p$ on $\widehat{\Isog}$ such that its restriction to the small \'{e}tale site of $\widehat{\M}$ is (non-canonically) homotopic to the $p$-complete version of the sheaf $\O^\top_\GHM$.
\end{cor}

To prove the above statement we will appeal to the uniqueness of $\O^\top_\GHM$, which in turn requires the concept of a \emph{natural elliptic cohomology theory}. We will follow the definition of \cite[Df.1]{uniqueotop} specialised to the case of a smooth elliptic curve; there are many variants in the literature.

\begin{mydef}
Let $E$ be an elliptic curve over a ring $R$ defining a morphism $E\colon \Spec R\to \M$. We say a homotopy commutative ring spectrum $A$ is an \emph{elliptic cohomology theory} for $E$ if we have the following data:
\begin{enumerate}
\item $A$ is weakly $2$-periodic.
\item The groups $\pi_d A$ vanish for all odd integers $d$.
\item There is a chosen isomorphism of rings $\pi_0 A\simeq R$.
\item There is a chosen isomorphism of formal groups between the formal group $\widehat{E}$ of $E$ and the classical Quillen formal group of $A$.
\end{enumerate}
\end{mydef}

Notice how closely these conditions match those of \Cref{luriestheoremintext}. The fact that the sections of $\O^\top$ define elliptic cohomology theories is a characterising feature of this sheaf; the following is mentioned in \cite[Th.1.2]{bourbakigoerss} and proven in \cite[Th.A]{uniqueotop}.

\begin{theorem}\label{uniquenessofotopreferencetheorem}
The sheaf $\O^\top_\GHM$ (resp.\ its $p$-completion) are uniquely defined as sheaves of $\E_\infty$-rings on the small \'{e}tale site of $\M$ (resp.\ $\widehat{\M}$) by the property that it defines natural elliptic cohomology theories on the subsite of affine schemes.
\end{theorem}

\begin{proof}[Proof of \Cref{pcompletetheorema}]
Using \Cref{luriestheoremintext}, it is easy to construct $\O^\top_p$. Indeed, writing $C$ for the site of \Cref{luriestheoremintext} (with $n=2$), we can define a $2$-functor $[p^\infty]\colon\widehat{\Isog}\to C$ sending a pair $(\x,E)$ to $(\x,E[p^\infty])$---this is well-defined as isogenies of elliptic curves of degree prime to $p$ induces isomorphisms on associated $p$-divisible groups by \Cref{plocalpdivisiblegroupcheck}. Define the functor $\O^\top_p$ as the composite
\[\O^\top_p\colon\widehat{\Isog}^\op\xrightarrow{[p^\infty]}C^\op\xrightarrow{\O^\top_\BTtwo} \CAlg.\]
Notice that $[p^\infty]$ is a morphism of \'{e}tale sites as both sites obtain this Grothendieck topology by the forgetful functor to formal Deligne--Mumford stacks and $[p^\infty]$ does not change the underlying formal Deligne--Mumford stack. From this, and the given fact that $\O^\top_\BTtwo$ is an \'{e}tale hypersheaf, we see that $\O^\top_p$ is also an \'{e}tale hypersheaf. We will now appeal to \Cref{uniquenessofotopreferencetheorem} to prove that the restriction of $\O^\top_p$ to the small \'{e}tale site of $\widehat{\M}$ yields the $p$-complete version of $\O^\top_\GHM$. After identifying the formal completion $\widehat{E}$ of an elliptic curve $E$ at its identity section with the identity component of the $p$-divisible group $E[p^\infty]$ (\cite[Pr.7.4.1]{ec2name}) one sees the conditions 1-4 of an elliptic cohomology theory correspond exactly to conditions 1-4 given in \Cref{luriestheoremintext} and satisfied by $\O^\top_\BTtwo$. This implies that $\O^\top_p$ restricts to a sheaf of elliptic cohomology theories and the uniqueness of the $p$-complete version of $\O^\top_\GHM$ yields the desired (non-canonical) homotopy.
\end{proof}

Patching together these sheaves $\O^\top_p$ with some rational data produces leads to an integral statement.

\begin{theorem}\label{isogextension}
There is an \'{e}tale sheaf of $\E_\infty$-rings $\O^\top$ on $\Isog$ such that restriction to the small \'{e}tale site of $\M$ is (non-canonically) homotopic to the sheaf $\O^\top_\GHM$.
\end{theorem}

\begin{proof}
Let us begin with the rational construction. We will implicitly use the Schwede--Shipley equivalence $\Mod_\Q\simeq D(\Q)$ of symmetric monoidal $\infty$-categories; see \cite{schwedeshipleystablemodulecategoresiarecategirofmodules} for the original statement and \cite[Th.7.1.2.13]{haname} for an $\infty$-categorical version. Define a sheaf $\O^\top_\Q$ on $\Isog$ valued in rational cdgas by assigning to each pair $(\xx,E)$ the formal $\Q$-cdga $\omega_{E/\xx}^{\otimes \ast}(\xx)\otimes \Q$, itself defined by placing each $\omega_{E/\xx}^{\otimes k}(\xx)$ in degree $2k$ for each integer $k$. Here $\omega_{E/\xx}$ is the line bundle on $\xx$ defined by pushing forward $\Omega^1_{E/\xx}$ through the structure map $E\to \xx$, or equivalently as the dualising line for the formal group $\widehat{E}$. To see this is well-defined on $\Isog$, notice that any nonzero homomorphism of elliptic curves induces an equivalence on global sections of dualising lines over $\Q$, as these rational vector spaces are $1$-dimensional. The functor $\O^\top_\Q$ defines an \'{e}tale hypersheaf of $\E_\infty$-rings as $\omega_E$ is such a hypersheaf of abelian groups. Consider the morphisms of sites
\[\pi^p\colon \Isog\to \widehat{\Isog}\]
induced by base change along $\Spf \Z_p\to \Spec \Z$ for each prime $p$. This yields an \'{e}tale hypersheaf $\prod_p \pi_\ast^p \O^\top_p$ on $\Isog$. We would now like to construct a morphism of \'{e}tale hypersheaves of $\E_\infty$-rings
\[\al\colon\O^\top_\Q\to \left(\prod_p \pi_\ast^p \O^\top_p \right)_\Q\]
where the codomain is the rationalisation of a previously defined hypersheaf. By standard arguments (see \cite[Lm.6.1.10]{mythesis} for example), to define this map of hypersheaves, it suffices to restrict one's attention to affine objects in $\Isog$, so those pairs $(\xx,E)$ where $\xx=\Spec R$ is an affine scheme. When evaluated on $(\Spec R,E)$ the domain is by definition the formal rational cdga $\omega_{E/R}^{\otimes \ast}(R)$. Similarly, by \cite[Pr.4.8]{lennartconnective} each $\pi_\ast^p \O^\top_p(\Spec R,E)[\frac{1}{p}]$ is a formal rational cdga, and by \Cref{luriestheoremintext} we see that when evaluated on $(\Spec R,E)$ the codomain is naturally equivalent to the formal rational cdga
\[\left(\prod_p \left(\omega^{\otimes \ast}_{E/R}(R)\right)^\wedge_p\right)_\Q.\]
The map $\al$, when evaluated on $(\Spec R,E)$, can then be defined as the rationalisation of the product over all primes of the natural completion map $\omega_{E/R}^{\otimes \ast}(R)\to (\omega_{E/R}^{\otimes\ast}(R))^\wedge_p$. The desired \'{e}tale hypersheaf $\O^\top$ is then defined by the pullback in the diagram
\begin{equation}\label{arithmeticsfracturesquare}\begin{tikzcd}
{\O^\top}\ar[r]\ar[d]		&	{\prod_p \pi_\ast^p \O^\top_p}\ar[d]	\\
{\O^\top_\Q}\ar[r, "\al"]	&	{\left(\prod_p \pi_\ast^p \O^\top_p\right)_\Q}
\end{tikzcd}\end{equation}
of \'{e}tale hypersheaves of $\E_\infty$-rings on $\Isog$. To check the restriction of this $\O^\top$ is (non-canonically) homotopic to the sheaf $\O^\top_\GHM$, we again appeal to \Cref{uniquenessofotopreferencetheorem}. Indeed, this follows as $\O^\top$ defines a sheaf of elliptic cohomology theories when restricted to the small \'{e}tale site of $\M$ by \cite[Lm.4.5]{hilllawson}, or more accurately, by \cite[Lm.6.1.11]{mythesis}. 
\end{proof}

The fact that $\O^\top$ defines a natural elliptic cohomology theory shows that applying $\O^\top$ to a morphism in $\Isog$ is still computationally tractable. To this end, recall the following construction of the \emph{descent spectral sequence} (DSS) for a section $\O^\top(\xx,E)$: first, take an affine \'{e}tale hypercover $\Spec A^\bullet$ of $\xx$. As $\O^\top$ is an \'{e}tale hypersheaf, the canonical map of $\E_\infty$-rings
\[\O^\top(\xx,E)\xrightarrow{\simeq} \lim \O^\top(\Spec A^\bullet,E_{A^\bullet}) \]
is an equivalence. The Bousfield--Kan spectral sequence associated to this cosimplicial limit yields a spectral sequence with $E_1$-page
\[E_1^{s,t}\simeq \pi_t \O^\top(\Spec A^s,E_{A^s})\]
converging to $\pi_{t-s}\O^\top(\xx,E)$. As there are natural isomorphisms of $\pi_0 A^s$-modules
\[\pi_{2t} \O^\top(\Spec A^s,E_{A^s})\simeq \omega^{\otimes t}_{E_{A^s}/A^s}(A^s)\]
and the dualising lines $\omega$ satisfy descent, we see the $E_2$-page of the DSS takes the form
\[E_2^{s,2t}\simeq H^s(\xx,\omega_{E/\xx}^{\otimes t})\]
and $E_2^{s,t}=0$ for odd $t$. This construction of the DSS is taken up in \cite[\textsection 5.3]{tmfbook}, for example.

\begin{prop}\label{pullbacketwocompiarson}
Let $(f,\phi)\colon (\xx,E)\to (\xx',E')$ be a morphism in $\Isog$. Then the morphism of $\E_\infty$-rings
\[(f,\phi)^\ast\colon \O^\top(\xx',E')\to \O^\top(\xx,E)\]
induces the composition of algebraic pullbacks on cohomology 
\[H^\ast(\xx',\omega_{E'/\xx'})\xrightarrow{f^\ast} H^\ast(\xx,\omega_{f^\ast E/\xx})\xrightarrow{\phi^\ast} H^\ast(\xx,\omega_{E/\xx})\]
between the $E_2$-pages of DSSs. In particular, if both $\xx$ and $\xx'$ are affine, then $\pi_{2\ast}(f,\phi)^\ast$ is the algebraic pullback on global sections under the natural isomorphisms
\[\pi_{2\ast}\O^\top(\xx,E)\simeq \omega^{\otimes \ast}_{E/\xx}(\xx)\qquad \pi_{2\ast}\O^\top(\xx',E')\simeq \omega^{\otimes \ast}_{E'/\xx'}(\xx').\]
\end{prop}

At its core, the agreement above comes from \Cref{luriestheoremintext}.

\begin{proof}
First, let us assume $\xx'=\Spec A$ and $\xx=\Spec B$ are affine. Using the naturality of the isomorphisms found in \Cref{luriestheoremintext}, we see that $(f,\phi)^\ast$, expressed as the composition $(\id,\phi)^\ast\circ(f,\can)^\ast$, induces the morphisms of $A$-modules on $\pi_2$
\[\omega_{E'/A}(A)\to f_\ast f^\ast\omega_{E'/A}(A)\simeq f_\ast \omega_{f^\ast E'/B}(A)\xrightarrow{f_\ast \xi} f_\ast\omega_{E/B}(A)\]
where the first map is the unit, the second map is the canonical base change isomorphism of \cite[Rmk.4.2.4]{ec2name}, and the third is induced by $\psi$. Indeed, the fact that $(f,\can)^\ast$ induces the first two morphisms follows exactly from the fact that $\O^\top$ defines a natural elliptic cohomology theory when restricted to $\M_\et$---this much is also true for $\O^\top_\GHM$. To see that $(\id,\phi)^\ast$, which we will now write as $\phi^\ast$, induces the morphisms suggested above, first use the arithmetic fracture square (\ref{arithmeticsfracturesquare}) applied to this morphism, hence it suffices to prove the rational case and the $p$-complete cases for every prime $p$. The rational case follows by construction. In the $p$-complete case, so $A$ and $B$ are both $p$-complete rings for some prime $p$, then we can use the natural isomorphism between the formal group of an abelian variety over $A$ (or $B$) and the formal group associated with the $p$-divisible group of that variety; see \cite[Pr.7.3.1]{ec2name}. The line bundles $\omega_{f^\ast E'/B}$ and $\omega_{E/B}$ only depend on the formal group of $f^\ast E'$ and $E$, respectively, so we are reduced to show that $\psi^\ast$ induces the algebraic map $\xi$ above. This follows precisely from the naturality of the isomorphisms found in part 4 of \Cref{luriestheoremintext} and our construction of $\O^\top$ using this theorem. This finishes the affine case.\\

For the general case, take an affine \'{e}tale hypercover $\Spec A^\bullet\to \xx'$, which exists as $\xx'$ is qcqs---a consequence of being \'{e}tale over $\M$. Pulling back this hypercover along $\xx\to \xx'$ we obtain an \'{e}tale hypercover $P_\bullet\to \xx$, which we inductively refine to an \emph{affine} \'{e}tale hypercover $\Spec B^\bullet\to \xx$. The composite
\[f_\bullet\colon \Spec B^\bullet \xrightarrow{g_\bullet} P_\bullet\xrightarrow{h_\bullet} \Spec A^\bullet\]
gives a morphism of simplicial stacks over $\xx'$. Consider the induced morphism of simplicial objects in $\Isog$
\[(f_\bullet, \phi_\bullet)\colon (\Spec B^\bullet, E_{B^\bullet})\to (\Spec A^\bullet, E'_{A^\bullet})\]
where $\phi^\bullet\colon E_{B^\bullet}\to (f^\bullet)^\ast E'_{A^\bullet}$ is the composite of the canonical isomorphism $E_{B^\bullet}\simeq g_\bullet^\ast E_{P_\bullet}$ and the base change of $\phi$ along $h^\bullet$. As the DSS for $\xx$ (resp.\, $\xx'$) is constructed as the Bousfield--Kan spectral sequence associated to an affine \'{e}tale hypercover with $E_1$-page given by $\pi_t \O^\top(\Spec B^s)$ (resp.\, $\pi_t\O^\top(\Spec A^s)$). We know that the morphism between these $E_1$-pages induced by $(f^\bullet,\phi^\bullet)^\ast$ agrees with the composition of the algebraic pullbacks on by the affine case above. By descent, we see these algebraic morphisms on the $E_1$-page induce the desired algebraic morphisms on the $E_2$-page, hence we obtain the desired agreement of morphisms between $E_2$-pages as these morphisms already agree on $E_1$-pages.
\end{proof}

The result above does not readily extend to the compactification $\overline{\M}$ of $\M$, as we currently lack a spectral algebro-geometric interpretation of $\O^\top$ in that case. %Some work in this direction is done in \cite[\textsection1]{adamsontmf}, using Lurie's theorem at heights $1$ and $2$, then forming a neighbourhood around the height $1$ part, and then gluing the results together with Goerss--Hopkins obstruction theory. The results are not at clear These (or a generalisation of Lurie's work in \cite{ec2name}) and currently does not produce a statement as clean as \Cref{isogextension}.\\

As a simple application of \Cref{isogextension}, let us also write down the moduli and (almost) integral description of Behrens' $Q(N)$ spectra; see \cite{ktwospheremark,dividedbetafamily}.

\begin{construction}\label{behrensqnconstruction}
Let $N$ be a positive integer which we implicitly invert, $\M_0(N)$ be the moduli stack of elliptic curves with a chosen \emph{cyclic} subgroup $H$ of order $N$, and $p,q\colon \M_0(N)\to \M$ be the structure and quotient maps given by $(E,H)\mapsto E$ and $(E,H)\mapsto E/H$, respectively; these stacks and morphisms are discussed in detail in \Cref{heckecompsection}. Moreover, write $\tau\colon \M_0(N)\to \M_0(N)$ for the involution given by $(E,H)\mapsto (E/H, \widecheck{H})$---this is discussed in the proof of \Cref{propertiesoflevelnstacks}. This yields the diagram of stacks
\[\begin{tikzcd}
{\M_0(N)}\ar[r, "{p}"{description}]\ar[r, "{\id}", swap, shift right = 3]\ar[r, "{\tau}", shift left = 3]	&	{\M\sqcup \M_0(N)}\ar[r, "{\id\sqcup q}", shift left = 2]\ar[r, "{\id\sqcup p}", swap, shift right = 2]	&	{\M}
\end{tikzcd}\]
indexed over the truncation $\Delta^{\op}_{\leq 2}$ of the simplex category and ignoring degeneracies. We can enhance the above diagram to one in $\Isog$ as follows:
\begin{itemize} 
\item Equip each object $\M$ and $\M_0(N)$ with their respective universal elliptic curves $\EE$ and $\EE_0(N)$. 
\item Equip each occurrence of $p$ with the canonical isomorphism $\EE_0(N)\simeq p^\ast\EE$, the left and upper-right occurrences of $\id$ with the identity on the associated elliptic curves, the lower-right occurrence of $\id$ with the $N$-fold multiplication map $[N]\colon \EE\to \EE$, and finally equip both $\tau$ and $q$ with the canonical quotient map $Q\colon \EE_0(N)\to \EE_0(N)/\H$ of the universal elliptic curve by its universal cyclic subgroup $\H$.
\end{itemize}
This produces the diagram in $\Isog$
\[\begin{tikzcd}
{(\M_0(N), \EE_0(N))}\ar[rr, "{(p,\can)}"{description}]\ar[rr, "{(\id,\id)}", swap, shift right = 3]\ar[rr, "{(\tau,Q)}", shift left = 3]	&&	{(\M,\EE)\sqcup (\M_0(N),\EE_0(N))}\ar[d, "{(\id,[N])\sqcup (q,Q)}", shift left = 2]\ar[d, "{(\id,\id)\sqcup (p,\can)}", swap, shift right = 2]		\\
&&	{(\M, \EE)}
\end{tikzcd}\]
again indexed over $\Delta^\op_{\leq 2}$ without degeneracies. As in \cite[\textsection1.1]{ktwospheremark}, one can check these morphisms satisfy the semisimplicial relations up to homotopy in $\Isog$---the $n$-simplices for $n\geq 3$ are all degenerate. Moreover, all of the morphisms above lie in $\Isog$, so we can apply $\O^\top$ to the above diagram by \Cref{isogextension} and obtain a semicosimplicial $\E_\infty$-ring $Q(N)^\bullet$. Behrens defines $Q(N)$ as the limit of this diagram of $\E_\infty$-rings and uses this $\E_\infty$-ring to study the $K(2)$-local stable homotopy category.
\end{construction}

Using this construction of $Q(N)$ and the operations of \Cref{zooofoperations}, the author and Venkata Sai Narayana Bavisetty are currently exploring various stable operations on these spectra.

%%%%%%%%%%%%%%%%%%%%%%%%%%%%%%%%%%%%%%%%%%%%%%%%%%%%%%%%%%

\subsection{Extension to a spectral Mackey functor}\label{mackeyfunctorsection}

The classical definitions of Hecke operators from number theory involve an averaging process. In algebraic geometry, these processes can be interpreted as a kind of transfer map, and classically such transfer maps are functorial and satisfy base change; see \cite[Expos\'{e} IX, \textsection5]{sgafour} or \cite[\href{https://stacks.math.columbia.edu/tag/03SH}{03SH}]{stacks}. These transfers also exist in spectral algebraic geometry, and one can ask if these spectral transfers satisfy base change and are functorial up to coherent homotopy. To encode our homotopy coherent transfers, we will use the $\infty$-categories of spans originally defined by Barwick \cite{spectralmackeyi} and the norm constructions of Bachmann--Hoyois \cite{normsmotivic}. For more on the following definition, see \cite[\textsection5]{spectralmackeyi} or \cite[\textsection C]{normsmotivic}.

\begin{mydef}
Given an $\infty$-category $C$ with pullbacks, equipped with a class of morphisms $M$ closed under composition and pullback, then there is an $\infty$-category $\Span_M(C)$ of \emph{$M$-spans in $C$}. The objects of $\Span_M(C)$ are those of $C$, 1-morphisms from $X$ to $Y$ are spans
\[X\xleftarrow{f} Z\xrightarrow{g} Y\]
where $f$ is any map in $C$ and $g$ belongs to $M$, and composition is given by pullback. Functors $\Span_M(C)\to \Sp$ are called \emph{spectral Mackey functors} \`{a} la \cite{spectralmackeyi}.
\end{mydef}

This leads us to our second extension to the traditional functorality of $\O^\top_\GHM$.

\begin{theorem}\label{spectralmackextension}
Let $\fin$ be the class of morphisms $(f,\phi)\colon(\xx,f^\ast E')\to (\xx',E')$ in $\Isog$ such that $f$ is finite and $\phi$ is naturally isomorphic to the identity. Then there is a functor $\OO^\top\colon \Span_\fin(\Isog)\to \Sp$ such that the diagram of $\infty$-categories
\[\begin{tikzcd}
{\Isog^\op}\ar[d]\ar[r, "{\O^\top}"]		&	{\CAlg}\ar[d]	\\
{\Span_\fin(\Isog)}\ar[r, "{\OO^\top}"]	&	{\Sp}
\end{tikzcd}\]
commutes. This also holds for $\widehat{\Isog}$.
\end{theorem}

Essentially, \cite[\textsection C]{normsmotivic} states that because finite \'{e}tale morphisms are \'{e}tale locally fold maps $X^{\sqcup n}\to X$ and $\O^\top$ is an \'{e}tale sheaf, we obtain a canonical extension of $\O^\top$ to a spectral Mackey functor.

\begin{proof}
Let $C$ be either $\Isog$ or $\widehat{\Isog}$ for a fixed prime $p$. This theorem will be an application of \cite[Cor.C.13]{normsmotivic}. Using the notation \emph{ibid}, this means that $C=C$, $t$ will be the \'{e}tale topology, $m$ the collection $\fin$, $F$ the sheaf $\O^\top$ or $\O^\top_p$, and $D$ the $\infty$-category with finite products $\Sp$. To apply \emph{ibid}, we need to check the hypotheses contained within:
\begin{itemize}
\item The site $C$ admits finite coproducts, that for all objects $X,Y$ in $C$ the fibre product $X\times_{X\sqcup Y}Y$ exists and is the initial object, and finite coproduct decompositions are stable under pullback. In other words, $C$ is \emph{extensive}.
\item There is a containment of $\infty$-categories $\Shv^\et(C)\subset \Shv^\sqcup(C)$, where the latter is the $\infty$-category of sheaves of spaces on $C$ with Grothendieck topology defined by setting covers to be finite disjoint unions.
\item All finite fold maps $\nabla \colon \coprod X\to X$ are contains in $\fin$.
\item Every morphism in $\fin$ is \'{e}tale locally such a fold morphism.
\end{itemize}
The first condition follows as finite coproducts in $\Isog$ are given in the stack and elliptic curve variable separately. The second and third conditions are clear. The fourth is classical (see \cite[\href{https://stacks.math.columbia.edu/tag/04HN}{04HN}]{stacks}) as finite \'{e}tale morphisms are \'{e}tale locally disjoint unions. It directly follows from \cite[Cor.C.13]{normsmotivic} that there exists an essentially unique dashed functor in the commutative diagram of $\infty$-categories
\[\begin{tikzcd}
{C^\op}\ar[d]\ar[r, "{\O^\top}"]	&	{\CAlg}\ar[d]		\\
{\Span_\fin(C)}\ar[r, "{\OO^\top}", dashed]	&	{\CAlg^\times(\Sp)}
\end{tikzcd}\]
where $\CAlg^\times \Sp$ is the $\infty$-category of $\E_\infty$-objects in $\Sp$ with respect to the Cartesian monoidal structure. To finish our proof, note the following natural equivalences of $\infty$-categories:
\[\CAlg^\times(\Sp)\xrightarrow{\simeq} \CAlg^\sqcup(\Sp)\xleftarrow{\simeq} \Sp\]
The first equivalence above follows from the fact that $\Sp$ is stable, so the canonical map from a finite coproduct to a finite product is an equivalence, and the second as every object in a cocomplete $\infty$-category has an essentially unique coCartesian monoidal structure.
\end{proof}

\begin{remark}%thanks mike!
One might first guess to define $\fin$ as the collection of all $(f,\phi)$ such that $f$ is finite without imposing restrictions on $\phi$, however, our proof of \Cref{spectralmackextension} would not carry over in this generality. In essence, one would have to show that any isogeny of elliptic curves $E\to E'$ of invertible degree over a base scheme $S$ is \'{e}tale locally on $S$ a fold map from a disjoint union. This cannot be true though, as over algebraically closed fields there are isogenies that are not fold maps. It is true that \'{e}tale locally on $E'$ this isogeny is a fold map, but our Grothendieck topology on $\Isog$ is defined on $S$ not $E$.
\end{remark}

\begin{remark}
The use of \cite[Cor.C.13]{normsmotivic}, as we have done in the proof above, applies to many situations in homotopy theory and spectral algebraic geometry. For example, it also applies to $\O^\top_\GHM$ considered as a log-\'{e}tale sheaf on the small log-\'{e}tale site of the compactification $\overline{\M}$ as constructed in \cite{hilllawson}. The same argument also works for any of the \'{e}tale hypersheaves $\O^\top_\BTn$ occuring in Lurie's theorem \Cref{luriestheoremintext} and their more general counterparts in \cite[Df.1.5]{luriestheorem}. 
\end{remark}

Together, the work above prove the main theorem of this section.

\begin{proof}[Proof of \Cref{extensiontheoremintro}]
Combine \Cref{isogextension,spectralmackextension}.
\end{proof}

Transfers are famously fickle in homotopy theory, however, those constructed from the above theorem are not so mysterious.

\begin{mydef}\label{algebraictransfers}
Let $\xx$ be a Deligne--Mumford stack and $\A$ be an $\O_\xx$-algebra such that the underlying $\O_{\xx}$-module is finite locally free of rank $d$. There is an \emph{algebraic transfer map} $\A\to \O_{\xx}$ of $\O_{\xx}$-modules defined by the composite of natural morphisms
\[\A\to \mathcal{H}\mathrm{om}_{\O_{\xx}}(\A,\A)\xleftarrow{\simeq}\mathcal{H}\mathrm{om}_{\O_{\xx}}(\A,\O_{\xx})\otimes_{\O_{\xx}}\A\to \O_{\xx}.\]
If $f\colon \yy\to \xx$ is a finite morphism of Deligne--Mumford stacks, we write $f_!\colon f_\ast\O_\yy\to \O_{\xx}$ for the transfer map associated with the map of $\O_{\xx}$-algebras $\O_{\xx}\to f_\ast\O_\yy$. Given another $\O_{\xx}$-module $\B$, we will also write
\[f_!\colon \B\otimes_{\O_{\xx}}f_\ast\O_\yy\to \B\]
for the tensor product of $f_!$ with $\B$. Moreover, we will further write
\[f_!\colon H^\ast(\xx,\B\otimes_{\O_{\xx}}f_\ast\O_\yy)\to H^\ast(\xx,\B)\]
for the induced transfer map on sheaf cohomology.
\end{mydef}

By \cite[Expos\'{e} IX, \textsection5]{sgafour} or \cite[\href{https://stacks.math.columbia.edu/tag/03SH}{03SH}]{stacks}, the algebraic transfer maps $f_!$ are uniquely determined by the fact that \'{e}tale locally on $\xx'$ they are the $n$-fold addition map.\\

Recall that given an object $(\xx,E)$ in $\Isog$, the \emph{descent spectral sequence} for $\O^\top(\xx,E)$ is the Bousfield--Kan spectral sequence associated to a choice of affine \'{e}tale hypercover $\Spec R^\bullet\to \xx$ as $\O^\top(\xx,E)$ is equivalent to the limit of the cosimplicial $\E_\infty$-ring $\O^\top(\Spec R^\bullet, E_{R^\bullet})$.

\begin{prop}\label{classicaltransfers}
Let $(f,\phi)\colon (\xx,E)\to (\xx',E')$ be a morphism in $\fin$. Then the morphism of $\O^\top(\xx',E')$-modules
\[(f,\phi)_!\colon \O^\top(\xx,E)\to \O^\top(\xx',E')\]
induces the algebraic transfer map
\[f_!\colon H^\ast(\xx,\omega_{E/\xx}^{\otimes\ast})\to H^\ast(\xx', \omega_{E'/\xx'}^{\otimes\ast})\]
between the $E_2$-pages of the DSSs for $\O^\top(\xx,E)$ and $\O^\top(\xx',E')$. In particular, if $\xx'$ is affine, hence $\xx$ is also affine, then $(f,\phi)_!$ induces the algebraic transfer map between the homotopy groups of $\O^\top(\xx,E)$ and $\O^\top(\xx',E')$.
\end{prop}

\begin{proof}
Let us drop the elliptic curves from our notation---no confusion should arise as $\phi$ is isomorphic to the identity.\\

By \cite[Cor.C.13]{normsmotivic}, the spectral Mackey functor $\OO^\top$ can be characterised as the right Kan extension of the functor
\[\O^\top_\fold\colon \Span_\fold(\Isog)\to \Sp\]
whose values on fold maps $\xx^{\sqcup n}\to \xx$ is given by the addition map of spectra:
\[\O^\top_\fold(\xx^{\sqcup n})\simeq \O^\top_\fold(\xx)^n\to \O^\top_\fold(\xx)\]
This also induces the addition map on homotopy groups. Let us now assume that $\xx'=\Spec A$ is affine, hence $\xx=\Spec B$ is also affine. There is an \'{e}tale cover $\Spec C\to \Spec A$ such that $C\otimes_A B$ is equivalent to $C^d$ as $C$-algebra, and $f_!\colon \OO^\top(B)\to \OO^\top(A)$ is \'{e}tale locally on $A$ the $d$-fold addition map. As there is a natural isomorphism $\pi_0\OO^\top(A)\simeq A$ and the map $f^\ast\colon \OO^\top(A)\to \OO^\top(B)$ is a finite \'{e}tale morphism of $\E_\infty$-rings, we see that $(f,\phi)_!$ is \'{e}tale locally on $\OO^\top(A)$ the $d$-fold addition map. Taking $\pi_0$, we see the morphism of $A$-modules $\pi_0 f_!\colon B\to A$ is \'{e}tale locally on $A$ the $d$-fold addition map, hence it agrees with the algebraic transfer map as they do so \'{e}tale locally. The same goes for the morphisms of invertible $A$-modules $(f,\phi)_!\colon \pi_k \O^\top(B)\to \pi_k\O^\top(A)$ for every $k\in \Z$. Indeed, there are the following natural (in $\M_\et$) isomorphism of $A$-modules for every integer $k$:
\[\pi_{2k}\O^\top(B)\simeq \omega_{E/B}^{\otimes k}(B)\simeq \omega_{E'/A}(A)^{\otimes k}\otimes_A B\qquad \pi_{2k}\O^\top(A)\simeq \omega_{E'/A}^{\otimes k}(A)\]
It follows that $\pi_{2k}$ of $(f,\phi)_!$ can be identified with the algebraic transfer map from \Cref{algebraictransfers}. Otherwise, both $\pi_{2k+1}\O^\top(B)$ and $\pi_{2k+1}\O^\top(A)$ vanish. In particular, we have proven \Cref{classicaltransfers} in the affine case.\\

For general $\xx'$, take an affine \'{e}tale hypercover $\Spec A^s\to \xx'$, and as $\xx\to \xx'$ is finite and hence also affine, we see that $\Spec B^s\simeq \Spec A^s\times_{\xx'}\xx$ is also affine. The $E_1$-page of the DSS of $\O^\top(\xx')$ (resp.\ $\O^\top(\xx)$) is given by the groups $\pi_t \O^\top(B^s)$ (resp. $\pi_t\O^\top(A^s)$). We have uniquely identified the map induced by the spectral transfers between these $E_1$-pages with the algebraic transfers in the affine case above. It follows that the induced map between $E_2$-pages also agrees with the algebraic transfer, as the functoriality of sheaf cohomology can be summarised by taking a \v{C}ech nerve, pulling back all the transfers along this nerve, and then using this resolution to calculate cohomology.
\end{proof}

The following corollary is a key step to proving \Cref{heckeoperatorscompositionintro}.

\begin{cor}\label{compositionsareconstants}
Suppose we are in the situation of \Cref{classicaltransfers}, where we further assume $\phi$ has constant degree $d$. If the column $E_\infty^{s,s}$ in the DSS for $\O^\top(\xx',E')$ is concentrated in degree $0$, then there is a (non-canonical) homotopy
\[(f,\phi)_!\circ (f,\phi)^\ast\simeq d\colon \OO^\top(\xx',E')\to \OO^\top(\xx',E')\]
of maps of $\OO^\top(\xx',E')$-modules, where $d$ indicates the $d$-fold addition map.
\end{cor}

The homotopy above is not canonical in any sense---we only claim these two maps agree in the homotopy category of $\OO^\top(\xx',E')$-modules. As the proof below will show, if the $E_\infty^{s,s+1}$ column vanishes for $s\geq0$, this homotopy will be unique up to some (unspecified) $2$-homotopy. One could continue this line of argument, but unless $E_\infty^{s,s+k}$ vanishes for all $s\geq 0$ and $t\geq 1$, this argument will never lead to a canonical choice of homotopy. Moreover, in this case the homotopy groups of $\pi_\ast\O^\top(\xx',E')$ is concentrated in nonpositive degrees, which is not a feature any of the explicit spectra in this article have. This point highlights the necessity of our spectral algebro-geometric approach to the functoriality of $\O^\top$, which provides a chosen lift of $\O^\top_\GHM$ to $\Isog$ and $\Span_\fin(\Isog)$. The traditional obstruction theoretic approaches require showing such a space of lifts is contractible, which is not true.

\begin{remark}\label{whentoapplyadditionmaplemma}
The condition on the descent spectral sequence is satisfied in many examples of interest for us:
\begin{enumerate}
\item If $\xx'$ is affine the descent spectral sequence is concentrated in degree zero and collapses on the $E_2$-page, so the condition holds. 
\item If $\xx'=\M_{R}$ where $R$ is any localisation of the integers, then by the calculation of \cite{bauer} (which applies to $\TMF$ as discussed in \cite{konter} or Chapter 13 of \cite{tmfbook}), we see the $E_\infty$-page of the associated descent spectral sequence is concentrated in filtration zero in the desired column.
\item If $\xx'=\M_\Ga$ is a moduli stack of elliptic curves with $\Ga$-level structure and we further assume that $\Ga$ is \emph{tame} as defined in \cite[Df.2.2]{tmfwls}, then by part 4 of \cite[Pr.2.5]{tmfwls}, the descent spectral sequence for $\TMF(\Ga)$ is concentrated in filtration zero. This applies to $\TMF(n)$ and $\TMF_1(n)$ for $n\geq 2$, as well as $\TMF_0(n)$ if we invert $\gcd(6,\phi(n))$, where $\phi(n)$ is Euler's totient function.
\end{enumerate}
\end{remark}

We will use the latter example to prove \Cref{heckeoperatorscompositionintro}. To prove a version of \Cref{heckeoperatorscompositionintro} without having to invert $\gcd(6,\phi(n))$, one might want the following conjecture to hold.

\begin{conjecture}\label{edgemaphasnotorsion}
For every $n\geq 1$ the natural edge map $\pi_0 \TMF_0(n)\to \MF_0(n)_0$ is an isomorphism.
\end{conjecture}

We have no reason to expect this conjecture to be true except for the empirical evidence supplied by others; it has been shown that \Cref{edgemaphasnotorsion} holds for $n=1$ \cite{bauer}, $n=2$ \cite{ktwospheremark}, $n=3$ \cite{levelonethree}, $n=5$ \cite{markkyle} at the prime $2$, and $n=7$ \cite{meieroroztmfo7} away from the prime $2$.

\begin{proof}[Proof of \Cref{compositionsareconstants}]
The group of homotopy classes of $\OO^\top(\xx',E')$-module endomorphisms of $\OO^\top(\xx',E')$ is naturally isomorphic to $\pi_0 \OO^\top(\xx',E')$, and our hypotheses show that the edge map $H^0(\xx',\O_{\xx'})\to \pi_0\OO^\top(\xx',E')$, from the descent spectral sequence for $\O^\top(\xx',E')$, is an isomorphism. It follows from \Cref{classicaltransfers} that the unit followed by a transfer induces multiplication by $d$ on sheaf cohomology, as this is true for the algebraic transfer.
\end{proof}

Let us outline how the existence of $\OO^\top$ endows the sections $\O^\top(\xx,E)$ with the structure of genuinely equivariant spectrum for the profinite group $\pi_1^\et \xx$, as defined by Barwick in \cite[\textsection B]{spectralmackeyi}.

\begin{remark}
Fix an object $(\xx,E)$ of $\Isog$. The extra structure of $\Isog$ will not be used in this example, so let us consider $\xx$ as an object of $\M_\et$, the small \'{e}tale site of $\M$, and drop $E$ from our notation. Let us write $\xx_\fet$ for the subsite of the small \'{e}tale site of $\xx$ spanned by Deligne--Mumford stacks $\yy$ whose structure morphism $\yy\to \xx$ is finite. If $\xx$ is a connected scheme and $x$ is a chosen geometric point of $\xx$, then following \cite[Ex.D.24]{spectralmackeyi} and Grothendieck, we see that there is a canonical equivalence of categories between $\xx_\fet$ and the category of finite $\pi_1^\et(\xx,x)$-sets with open stabiliser, denoted as $\Fin_{\pi_1^\et(\xx,x)}^\open$. Using this equivalence and the inclusion $\xx_\fet\to \M_\et$ we obtain a functor
\[\Span(\Fin_{\pi_1^\et(\xx,x)}^\open)\simeq \Span(\xx_\fet)\to \Span_\fin(\M_\et)\to \Span_\fin(\Isog)\xrightarrow{\OO^\top} \Sp\]
whose value on the one point set is $\O^\top(\xx)$. This construction is natural in morphisms $\xx\to \xx'$ in $\M_\et$ meaning the map $\O^\top(\xx')\to\O^\top(\xx)$ is $\pi_1^\et(\xx,x)$-equivariant. One can try to push this observation further, but unless the reader has something specific in mind, a certain global homotopy type \`{a} la Schwede \cite{s} for example, the framework that encodes the most symmetries of sections of $\O^\top$ is simply the functor $\OO^\top$ itself.
\end{remark}

%%%%%%%%%%%%%%%%%%%%%%%%%%%%%%%%%%%%%%%%%%%%%%%%%%%%%%%%%%
%%%%%%%%%%%%%%%%%%%%%%%%%%%%%%%%%%%%%%%%%%%%%%%%%%%%%%%%%%
%%%%%%%%%%%%%%%%%%%%%%%%%%%%%%%%%%%%%%%%%%%%%%%%%%%%%%%%%%

\section{Into the zoo of stable operations}\label{zooofoperations}

The goal of this section is to demonstrate how \Cref{extensiontheoremintro} leads to simple definitions of stable operations on $\TMF$, and that these operations possess remarkable coherence properties amongst one another.

%%%%%%%%%%%%%%%%%%%%%%%%%%%%%%%%%%%%%%%%%%%%%%%%%%%%%%%%%%

\subsection{Definitions and comparisons}\label{definitionssection}

The simplest (non-identity) isogeny between elliptic curves is probably the $k$-fold multiplication isogeny $[k]\colon E\to E$, where $k$ is any integer. This exists for all elliptic curves, hence it exists on the universal elliptic curve $\EE$ over the moduli stack of elliptic curves $\M$. By \cite[Th.2.3.1]{km}, this isogeny has degree $k^2$, so working over $\M_{\Z[\frac{1}{k}]}$ we obtain an isogeny of invertible degree $[n]\colon \EE\to \EE$.

\begin{mydef}\label{adamsdefinition}
Let $k$ be an integer. The \emph{stable Adams operation}
\[\psi^k\colon \TMF[\frac{1}{k}]\to \TMF[\frac{1}{k}]\]
is defined by applying $\O^\top$ to the morphism
\[(\M_{\Z[\frac{1}{k}]},\EE)\xrightarrow{(\id,[k])}(\M_{\Z[\frac{1}{k}]},\EE)\]
inside $\Isog$. One can also define \emph{$p$-adic stable Adams operations} $\psi^\lambda\colon \TMF_p\to \TMF_p$ for each prime $p$ and each $p$-adic unit $\lambda\in \Z_p^\times$---this was done in \cite[\textsection5.5]{luriestheorem} by applying $\O^\top_\BTtwo$ (of \Cref{luriestheoremintext}) to the $\lambda$-fold multiplication map on the $p$-divisible group $\EE[p^\infty]$ of the universal elliptic curve $\EE$.
\end{mydef}

The natuality of the above definition in the category $\Isog$ should be stressed: for any pair $(\xx,E)$ in $\Isog$ such that the defining map $f\colon \xx\to\M$ factors through $\M_{\Z[\frac{1}{k}]}$, the diagram in $\Isog$ commutes
\[\begin{tikzcd}
{(\xx,E)}\ar[d, "{(f,\can)}"]\ar[r, "{(\id,[k])}"]		&	{(\xx,E)}\ar[d, "{(f,\can)}"]	\\
{(\M_{\Z[\frac{1}{n}]},\EE)}\ar[r, "{(\id,[k])}"]	&	{(\M_{\Z[\frac{1}{n}]},\EE).}
\end{tikzcd}\]
Indeed, this boils down to the fact that $f^\ast[k]\colon f^\ast \EE\to f^\ast\EE$ is precisely the $k$-fold multiplication map on $f^\ast\EE\simeq E$. Applying $\O^\top$ to the above diagram provides a natural homotopy witnessing that the diagram of $\E_\infty$-rings
\begin{equation}\label{naturalityofadamsoperations}\begin{tikzcd}
{\TMF[\frac{1}{k}]}\ar[r, "{\psi^k}"]\ar[d, "{f^\ast}"]	&	{\TMF[\frac{1}{k}]}\ar[d, "{f^\ast}"]	\\
{\O^\top(\xx,E)}\ar[r, "{\psi^k}"]				&	{\O^\top(\xx,E)}
\end{tikzcd}\end{equation}
commutes. \\

Universal isogenies of degree $n$ lead to the definition of the $n$th stable Hecke operator.\\

Let $n$ be a positive integer. Write $\M_n$ for the \emph{moduli stack of elliptic curves with chosen finite subgroup of order $n$} over $\Spec \Z[\frac{1}{n}]$. In other words, the functor $\M_n\colon \CRing\to \Spc$ defined by sending a ring $R$ wherein $n$ is invertible to the 1-truncated space of pairs $(E,H)$ of elliptic curves $E$ over $R$ with a finite subgroup $H$ of order $n$, and for all rings $R$ wherein $n$ is not invertible to the empty space. We will see in \Cref{propertiesoflevelnstacks} that the stack $\M_n$ is a Deligne--Mumford stack and that the two maps $p,q\colon \M_n\to \M_{\Z[\frac{1}{n}]}$ are \'{e}tale, defined on $R$-valued points by $p(E,H)=E$ and $q(E,H)=E/H$. Write $(\EE_n,\H_n)$ for the universal pair over $\M_n$, so there is a canonical equivalence $\can\colon \EE_n\simeq p^\ast\EE$ and $\H_n$ is the universal subgroup of $\EE_n$ of order $n$. Let us also write $Q\colon \EE_n\to \EE_n/\H_n\simeq q^\ast\EE$ for the canonical quotient.

\begin{mydef}\label{stableheckedefinition}
Let $n$ be a positive integer. The \emph{stable Hecke operator}
\[\T_n\colon \TMF[\frac{1}{n}]\to \TMF[\frac{1}{n}]\]
is defined by applying $\OO^\top$ to the morphism
\[(\M_{\Z[\frac{1}{n}]},\EE) \xleftarrow{(q,Q)} (\M_n, \EE_n) \xrightarrow{(p,\can)} (\M_{\Z[\frac{1}{n}]},\EE)\]
inside $\Span_\fin(\Isog)$.
\end{mydef}

Just as the stable Adams operations are ``natural in $\Isog$'' (see (\ref{naturalityofadamsoperations})), the stable Hecke operators have similar naturality properties, although there is a subtlety with regards to \emph{descent data}; see \Cref{generalisedheckeoperators}.

\begin{remark}
For the moment, let use write $\M_{\Isog_n}$ for the \emph{moduli stack of elliptic curves and isogenies of degree $n$} over $\Spec \Z[\frac{1}{n}]$, meaning that an $R$-valued point of this stack is a triple of $\phi\colon E_s\to E_t$ of two elliptic curves over $\Spec R$ together with an isogeny $\phi$ of degree $n$. This stack can be identified with $\M_n$ by sending a triple $(E_s,E_t,\phi)$ to the pair $(E_s,\ker \phi)$, and conversely sending a pair $(E,H)$ to the triple $(E,E/H,q)$, where $q$ is the evident quotient. The span of \Cref{stableheckedefinition} can now be written in terms of the \emph{universal isogeny} $\Phi\colon \EE_s\to \EE_t$ of degree $n$ over $\M_{\Isog_n}$
\[(\M_{\Z[\frac{1}{n}]},\EE) \xleftarrow{(t,\Phi)} (\M_{\Isog_n}, \EE_s) \xrightarrow{(s,\can)} (\M_{\Z[\frac{1}{n}]},\EE)\]
where $s,t\colon \M_{\Isog_n}\to \M_{\Z[\frac{1}{n}]}$ are the source and target maps, sending a triple $(E_s,E_t,\phi)$ to either $E_s$ or $E_t$.
\end{remark}

Now to stable Atkin--Lehner involutions. Fix an integer $N\geq 2$ and a positive divisor $Q$ of $N$ such that $\gcd(Q,N/Q)=1$. Define the \emph{Atkin--Lehner involution} $w_Q$ on the \emph{moduli stack $\M_0(N)$ of elliptic curves with chosen cyclic subgroup of order $N$} as follows: given a pair $(E,H)$ inside $\M_0(N)(S)$ for a fixed base scheme $S$, then $H$ splits uniquely into the product of subgroups $H_Q$ and $H_{N/Q}$ of order $Q$ and $N/Q$, respectively. The quotient isogeny determined by the subgroup $H=H_Q\times H_{N/Q}$ can be factorised in two different ways:
\[\begin{tikzcd}
	&	{E/H_Q}\ar[rd]	&	\\
{E}\ar[ru, "{\phi_Q}"]\ar[rd, "{\phi_{N/Q}}", swap]\ar[rr]	&&	{E/H.}	\\
	&	{E/H_{N/Q}}\ar[ru]	&
\end{tikzcd}\]
We then define $w_Q(E,H)=(E/H_Q,K)$ where $K$ is the kernel of the isogeny $\phi_{N/Q}\circ \widecheck{\phi}_Q$ and $\widecheck{(-)}$ denotes the dual isogeny. It is clear that $K$ is cyclic of degree $N$ as $\phi_m$ has cyclic kernel of degree $m$ for $m=Q$ and $N/Q$, and $\gcd(Q,N/Q)=1$. If $Q=N$, then we call $w_N$ the \emph{Fricke involution} of $\M_0(N)$. By \Cref{propertiesoflevelnstacks} (to come), we see that $\M_0(N)$ naturally lies in $\Isog$ and one easily checks $w_Q\circ w_Q$ is naturally isomorphic to the identity.

\begin{mydef}\label{stableaktinlehnerdefinition}
Let $N\geq 2$ be an integer and $Q$ be a positive divisor of $N$ such that $\gcd(Q,N/Q)=1$. The \emph{stable Atkin--Lehner involution}
\[w_Q\colon \TMF_0(N)\to \TMF_0(N)\]
is defined by applying $\O^\top$ to the morphism
\[(w_Q, \phi_{Q})\colon (\M_0(N), \EE_0(N))\to (\M_0(N),\EE_0(N))\]
inside $\Isog$, where $\phi_Q\colon \EE_0(N)\to \EE_0(N)/\H_Q=w_Q^\ast \EE_0(N)$ is the natural quotient by the $Q$-primary part of the universal cyclic subgroup of $\EE_0(N)$ of order $N$. If $N=Q$, we call $w_N$ the \emph{stable Fricke involution} of $\TMF_0(N)$.
\end{mydef}

The three definitions above give rise to a general pattern. Let $N$ be a positive integer and $\Ga$ be a \emph{congruence subgroup} of level $N$ as in \cite[Df.1.2.1]{diamondshur}, ie, a subgroup of $\GL_2(\Z)$ containing the kernel of the reduction map $\GL_2(\Z)\to \GL_2(\Z/N\Z)$. There are moduli stacks of elliptic curves $\M_\Ga$ with \emph{$\Ga$-level structure} which are finite \'{e}tale over $\M_{\Z[\frac{1}{N}]}$; see \cite[\textsection IV.3]{dr}. We will write $p_\Ga\colon \M_\Ga\to \M_{\Z[\frac{1}{N}]}$ for the structure maps and $\EE_\Ga=p_\Ga^\ast\EE$ for the universal elliptic curve over $\M_\Ga$. Let us write $\TMF(\Ga)$ for the $\E_\infty$-$\TMF[\frac{1}{N}]$-algebra $\O^\top(\M_\Ga)$.

\begin{mydef}\label{doublecosetoperationdef}
Let $\Ga$ be a congruence subgroup of level $N\geq 1$, and let $\Ga_1,\Ga_2\leq \Ga$ be two finite index subgroups together with an isomorphism $\phi\colon\Ga_1\to \Ga_2$ and an isogeny $\Phi\colon \EE_{\Ga_1}\to \phi^\ast\EE_{\Ga_2}$ of invertible degree. The \emph{stable operator of the triple $(\Ga,\phi, \Phi)$}
\[(\Ga,\phi,\Phi)\colon \TMF(\Ga)\to \TMF(\Ga)\]
is defined by applying $\OO^\top$ to the span
\[(\M_{\Ga},\EE_{\Ga})\xleftarrow{(p_{\Ga_2},\can)}(\M_{\Ga_2},\EE_{\Ga_2})\xleftarrow{(\phi,\Phi)}(\M_{\Ga_1},\EE_{\Ga_1})\xrightarrow{(p_{\Ga_1},\can)} (\M_{\Ga},\EE_{\Ga})\]
inside $\Isog$.
\end{mydef}

Stable Adams operations and stable Atkin--Lehner involutions are of the above form. The same is true for the $\ell$th stable Hecke operator for a prime $\ell$, if one takes $\Ga=\GL_2(\Z)$,
\[\Ga_1=\left\{\left.\begin{pmatrix}	a & b \\ c & d	\end{pmatrix}\,\right|\, c\equiv_\ell 0\right\}\qquad \Ga_2=\left\{\left.\begin{pmatrix}	a & b \\ c & d	\end{pmatrix}\,\right|\, b\equiv_\ell 0\right\}\]
$\phi$ to be conjugation by $\begin{pmatrix}	\ell & 0 \\ 0 & 1	\end{pmatrix}$, and $\Phi$ to be the quotient of the universal elliptic curve over $\M_{\Ga_1}=\M_0(\ell)$ by the universal cyclic $\ell$-subgroup, which has invertible degree when working over $\Z[\frac{1}{\ell}]$.\\

For maximal generality, let us define operations directly associated with spans inside $\Isog$.

\begin{mydef}\label{maximumaction}
Consider a morphism $(\xx,E)\xleftarrow{(f,\phi)}(\yy,F)\xrightarrow{(g,\psi)}(\xx',E')$ inside $\Span_\fin(\Isog)$. The stable operation
\[\T^{(f,\phi)}_{(g,\psi)}\colon \O^\top(\xx,E)\to \O^\top(\xx',E')\]
is defined by applying $\O^\top$ to the span above.
\end{mydef}

The notation above is supposed to suggest variance.\\

By construction, we see that \Cref{adamsdefinition,stableheckedefinition,stableaktinlehnerdefinition,doublecosetoperationdef} all fall under the umbrella of \Cref{maximumaction}. This latter approach can also be used to generalise \Cref{adamsdefinition,stableheckedefinition,stableaktinlehnerdefinition} to other sections of $\O^\top$ given that they come equipped with certain \emph{descent data}; let us write down this generalisation in the case of stable Hecke operators.

\begin{var}\label{generalisedheckeoperators}
Let $f\colon \xx\to \M_{\Z[\frac{1}{n}]}$ be an object inside $\Isog$ and define $\xx_n$ using the pullback of stacks
\[\begin{tikzcd}
{\xx_n}\ar[r, "{p'}"]\ar[d, "{f'}"]	&	{\xx}\ar[d, "{f}"]	\\
{\M_n}\ar[r, "{p}"]			&	{\M_{\Z[\frac{1}{n}]}}
\end{tikzcd}\]
where $p\colon \M_n\to \M_{\Z[\frac{1}{n}]}$ is the structure map of \Cref{stableheckedefinition}. Suppose that the pair $(\xx,E)$ comes equipped with \emph{descent data with respect to subgroups of order $n$}, meaning there is a morphism of stacks $q'\colon \xx_n\to \xx$ and a homotopy $\al$ witnessing that the diagram of stacks
\begin{equation}\label{homotopywitnessingdescentdata}\begin{tikzcd}
{\xx_n}\ar[r, "{q'}"]\ar[d, "{f'}"]	&	{\xx}\ar[d, "{f}"]	\\
{\M_n}\ar[r, "{q}"]			&	{\M_{\Z[\frac{1}{n}]}}
\end{tikzcd}\end{equation}
commutes. We will discuss some examples and non-examples of such data shortly. In this case, we define a stable Hecke operator for $(\xx,E)$ by applying $\OO^\top$ to the composite
\[(\xx,E)\xleftarrow{(q',\can)} (X_n,q'^\ast E\simeq f'^\ast q^\ast \EE)\xleftarrow{(\id,f'^\ast \Phi)} (\xx,f'^\ast p^\ast \EE \simeq p'^\ast E)\xrightarrow{(p',\can)} (\xx,E)\]
inside $\Span_\fin(\Isog)$, where $\Phi$ is the universal isogeny of \Cref{stableheckedefinition}. Notice that both $q'$ and $\al$ are implicitly used in the above definition. For readability, we will not continue a general study of these operations and continue our study of the universal example, however, all of the results of \Cref{interrelationssection,decompiingstackssection,heckecompsection,pprimaryfactorisationsection} generalise once one has a grasp on the functoriality of above descent data; see \cite[\textsection3.1]{sigmaishinfty}.
\end{var}

Substituting $\xx=\M_n$ into this variant above gives some intution regarding the existence or nonexistence of cetrain descent data.

\begin{example}\label{descentdataexample}
The pair $(\M_n,\EE_n)$ has a canonical choice of descent data for subgroups of order $m$ so long as $m$ and $n$ are coprime. Indeed, if $\gcd(m,n)=1$ then the stack $(\M_n)_m$, defined by substituting $\M_n=\xx$ into \Cref{generalisedheckeoperators}, has $R$-valued points given by pairs $(E,H)$ where $H\leq E$ is a finite flat subgroup of degree $mn$. The morphism $q'\colon (\M_n)_m\to \M_n$ defined by sending a pair $(E,H)$ to $(E/H_m,H/H_m)$, where $H_m$ is the unique subgroup of $H$ of order $m$, provides half of the descent data. Unwinding the definitions, the other half of the descent data, the $\al$ witnessing that (\ref{homotopywitnessingdescentdata}) commutes, is defined by the isomorphism of elliptic curves $q'^\ast \EE_n\simeq \EE_n/(\H_n)_m$ confirming that $q'$ ``is defined by taking a quotient by the subgroup $H_m$''. More generally, any of the stacks of the form $\M_A$ from \Cref{stacksforheckedecomposition} come equipped with descent data for subgroups of order prime to the order of $A$.\\

If $\gcd(m,n)\neq 1$, the above definitions are not possible. Let us take $m=n$ for simplicity. In this case, the stack $(\M_n)_n$ has $R$-valued points triples $(E,K,H)$, where $K,H\leq E$ are both subgroups of order $n$. The na\"{i}ve guess at a map $q'\colon (\M_n)_n\to \M_n$ now fails, as sending a triple $(E,K,H)$ to $(E/K,H/K)$ does not make sense, and $K$ might not be contained in $H$, if even if $K\leq H$ there is no guarantee that $H/K$ would have order $n$. This is why we cannot define $\T_n$ operators on $\TMF_0(n)$ in \Cref{stableheckeandothers}, for instance.
\end{example}

Inspired by \cite{derivedheckealgebra}, we also mention degree-shifting versions of stable Hecke operators.

\begin{var}\label{derivedhcekeoperators}
Fix a positive integer $n$ and a class $\al$ in $\pi_d \TMF_n$ for some integer $d$, where $\TMF_n=\O^\top(\M_n)$. Define the \emph{derived stable Hecke operator} as the composite
\[\T^\al_n\colon \TMF[\frac{1}{n}][d]\xrightarrow{(q,Q)^\ast} \TMF_n[d]\xrightarrow{\cdot\al}\TMF_n\xrightarrow{(p,\can)_!}\TMF[\frac{1}{n}] \]
where $X[d]=\Sigma^d X$ denotes the $d$-fold suspension endofunctor of $\Sp$. This definition has the advantage of being a cohomology operation of (potentially) nonzero degree, and we might speculate that it should be related to a spectral version of a \emph{derived Hecke algebra}. We leave a general study of such operations for another time.
\end{var}

With all of these definitions out of the way, let us show how the set-up of \Cref{extensionsection} easily leads us to the comparison result from the introduction.

\begin{proof}[Proof of \Cref{comparisonintro}]
The diagrams in question commute if we replace the lower horizontal maps with that induced on DSSs by the upper horizontal map---this is just the naturality of the edge morphisms---hence we are reduced to studying these induced morphisms on DSSs. Using \Cref{pullbacketwocompiarson,classicaltransfers}, these induced maps on DSSs are identified with the associated algebraic morphisms.\\

For the stable Adams operations, we are left to calculate what the morphism $(\id,[k])^\ast$ induces on the sheaf $\omega_{E/\M}$, and this is classical: writing $\omega_{E/\M}$ as the dual of the Lie algebra of $\widehat{E}$, then the morphism induced by $[k]$ on this Lie algebra is multiplication-by-$k$; for more details, see \cite[Pr.5.28]{luriestheorem}. For Hecke operators, we note that our Hecke composition formula of \Cref{heckeoperatorscompositionintro} (proven independently to \Cref{comparisonintro}) and the classical Hecke composition formula reduce us to the case of primes, so $n=\ell$ is a prime. In this case, we note that the classical construction of $\ell\T_\ell^\cl$ takes the form
\[H^0(\M,\omega_{\EE/\M}^{\otimes\ast})\xrightarrow{q_\ell^\ast} H^0(\M_0(\ell),\omega_{q_\ell^\ast\EE/\M_0(\ell)}^{\otimes\ast})\xrightarrow{\xi^{\otimes\ast}}H^0(\M_0, \omega^{\otimes\ast}_{\EE/\M_0(\ell)})\xrightarrow{p^\ell_!}H^0(\M,\omega_{\EE/\M}^{\otimes\ast})\]
as in the introduction; see \cite[(4.5.1)]{conrad} in the compactified case. These three maps above are identified with the maps induced on $E_2$-pages by $(q,\can)^\ast$, $(\id,Q)^\ast$, and $(p_\ell,\id)_!$. A similar argument works for our stable Atkin--Lehner involutions using \Cref{pullbacketwocompiarson,classicaltransfers}; the relevant algebraic constructions can be found \cite{atkinlehnermodular}.
\end{proof}

Let us write down a corollary of our proof above for later reference.

\begin{cor}\label{etwoaggreement}
In the situation of \Cref{comparisonintro}, the operations $\psi^k$, $w_Q$, and $\T_n$ induce their algebraic counterparts on the $E_2$-page of the associated DSSs.
\end{cor}

By ``algebraic counterparts'' we cannot forget the factor of $n$ between the stable $\T_n$ and $\T_n^\cl$ as in \Cref{comparisonintro}.

\begin{remark}
\Cref{comparisonintro} does not constitute a complete calculation of our stable operations on the homotopy groups of $\TMF$. Indeed, this theorem does not say that $\T_n(c_4)$ is $\epsilon$-torsion or not, for odd $n$. Here $c_4$ is the unique class in $\pi_8 \TMF$ detecting the integral Eisenstein series $E_4$ under the comparison map $\pi_{2d}\TMF\to \MF_d$ which is $\kappa$-torsion, and $\epsilon$ and $\kappa$ are the unique nonzero classes in the image of the unit map $\Sph\to \TMF$ in degrees $8$ and $14$, respectively. Some of these deeper calculations are discussed in \cite[\textsection2]{adamsontmf} in the case of $\psi^k$, for example.
\end{remark}

%As we are comparing our operations to the classical ones, let us go one step further and compare our stable operations with those %stable operations of Andrew Baker on Landweber exact elliptic cohomology.

%\begin{prop}
%Let $n$ be a positive integer. Write $\T_n^\Baker\colon \TMF[\frac{1}{6n}]\to \TMF[\frac{1}{6n}]$ for Baker's $n$th stable Hecke operator of \cite{bakerhecketwo}. Then there is the following noncanonical homotopy of morphisms of spectra:
%\[n\T_n^\Baker\simeq \T_n\colon \TMF[\frac{1}{6n}]\to \TMF[\frac{1}{6n}]\]
%\end{prop}

%\begin{proof}
%\end{proof}

\begin{remark}
The stable Hecke operators of Baker \cite{bakerhecketwo} $\T^\Baker_n\colon \TMF[\frac{1}{6n}]\to \TMF[\frac{1}{6n}]$ on elliptic cohomology with $6$ inverted have many of the same features as our operations; here we use the modern notation of $\TMF[\frac{1}{6n}]$ for what is often written as $\ELL[\frac{1}{6n}]$. For example, Baker's operators agree with the classical Hecke operators on homotopy groups and there is a decomposition formula akin to \Cref{heckeoperatorscompositionintro}. It seems likely that there is a homotopy $\T_n\simeq n\T_n^\Baker$ as endomorphisms of the spectrum $\TMF[\frac{1}{6n}]$, although we do not claim to have a proof of this. %Indeed, as this spectrum is Landweber exact, it suffices to check these two endomorphisms induce the same morphism on homology theories. Using \Cref{heckeoperatorscompositionintro}, we reduce ourselves to the case where $n=\ell$ is a prime and then one should explicitly compare our geometric construction with Baker's on the Tate curve over $\Z[\frac{1}{6\ell}]((q))$.
\end{remark}

%%%%%%%%%%%%%%%%%%%%%%%%%%%%%%%%%%%%%%%%%%%%%%%%%%%%%%%%%%

\subsection{Interrelations between stable operators}\label{interrelationssection}

Just as the various $n$-fold multiplication maps on an elliptic curve interact with one another, there is a collection of immediate relations between our various operations defined in the previous section.\\

All of the proofs in this subsection will follow the same outline: to find the desired homotopy, we use the definitions of our operations as the image of morphisms in $\Span_\fin(\Isog)$ under the functor $\OO^\top$ to reduce ourselves to finding homotopies in $\Span_\fin(\Isog)$, and quite often in the wide subcategory $\Isog$. We are then reduced to a purely algebraic question. As a springboard for increasingly complicated versions of this outline, let us start by detailing the simplest case of the composition of two stable Adams operations.

\begin{prop}\label{stableadamscompintext}
Let $k,\ell$ be two integers. Then there exist natural homotopies
\[\psi^1\simeq \psi^{-1}\simeq \id\colon \TMF\to \TMF\]
\[\psi^k\circ \psi^\ell\simeq \psi^{k\ell}\simeq \psi^\ell\circ \psi^k \colon \TMF[\frac{1}{\ell k}]\to \TMF[\frac{1}{\ell k}]\]
of morphisms of $\E_\infty$-rings.
\end{prop}

The second homotopy above $\psi^{-1}\simeq \id$ is similar to the fact that $\psi^{-1}$ acts trivially on $\KO$. This property is not shared by all sections of $\O^\top$ though. For example, $\psi^{-1}$ acts nontrivially on $\TMF_1(N)$ for $N\geq 3$; the $N=3$ case is explored at great length in \cite{lennartandhill} and more examples (and non-examples) are given in \cite[\textsection3.1]{realspectra} under the guise of \emph{real spectra}.

\begin{proof}
As $\psi^1=\O^\top(\id,[1])$ and $(\id,[1])=(\id,\id)$ in $\Isog$, then the functoriality of $\O^\top$ shows $\psi^1$ is naturally homotopic to $\id$. For $\psi^{-1}$, consider the homotopy
\[[-1]\colon (\id,\id)\to (\id,[-1])\]
inside $\Isog$. It looks a little silly in formulae, but it defines a homotopy in $\Isog$ and hence we obtain a natural homotopy from $\id$ to $\psi^{-1}$ in $\CAlg$. Next, we use the functoriality of $\O^\top$ in isogenies of invertible degree to obtain the desired chain of natural homotopies:
\[\psi^k\circ \psi^\ell=\O^\top(\id,[k])\circ \O^\top(\id, [\ell])\simeq \O^\top(\id, [\ell]\circ [k])\simeq \O^\top(\id,[k \ell])=\psi^{k\ell}\qedhere\]
\end{proof}

Let us give a formulation for how coherent the above homotopies are.

\begin{remark}\label{coherenthomotopiesremark}
Say we are given three integers $j,k,\ell$, and write $H_{k,\ell}$ for the homotopy of \Cref{stableadamscompintext} between $\psi^k\circ \psi^\ell$ and $\psi^{k\ell}$. Then one can homotope between $\psi^j\circ\psi^k\circ\psi^\ell$ and $\psi^{jk\ell}$ as morphisms of $\TMF[\frac{1}{jk\ell}]$ in two different ways:
\begin{equation}\label{wantatwocellplease}\begin{tikzcd}
{\psi^j\circ\psi^k\circ\psi^\ell}\ar[r, "{H_{j,k}\circ\psi^\ell}"]\ar[d, "{\psi^j\circ H_{k,\ell}}"]	&	{\psi^{jk}\circ \psi^\ell}\ar[d, "{H_{jk,\ell}}"]	\\
{\psi^j\circ\psi^{k\ell}}\ar[r, "{H_{j,k\ell}}"]								&	{\psi^{jk\ell}.}
\end{tikzcd}\end{equation}
We claim that there exists a natural 2-cell witnessing the above commutative diagram as one in the space of maps of $\E_\infty$-rings on $\TMF[\frac{1}{jk\ell}]$. Indeed, the morphisms and homotopies above are all obtained by applying $\O^\top$ to the diagram
\[\begin{tikzcd}
{[j]\circ[k]\circ[\ell]}\ar[r]\ar[d]	&	{[jk]\circ[\ell]}\ar[d]	\\
{[j]\circ[k\ell]}\ar[r]			&	{[jk\ell]}
\end{tikzcd}\]
in the set of morphisms of $\EE$ over $\M_{\Z[\frac{1}{jk\ell}]}$. As this set of morphisms has contractible higher cells there is a canonical choice of 2-cells for the above diagram, and applying $\O^\top$ yields a canonical 2-cell for (\ref{wantatwocellplease}). We give some alternative ways to view this coherence in \Cref{pprimaryfactorisationsection}. Rather than write down more of these coherences, we just state that they follow from the definitions of our stable operations and the fact that $\O^\top$ and $\OO^\top$ are functors of $\infty$-categories.
\end{remark}

Let us move on to the Atkin--Lehner involutions, which is a misleading name for these morphisms; see the second homotopy in the following.

\begin{prop}\label{stablealcompintext}
Let $k$ be an integer, $N$ be a positive integer, and $Q,R$ be two divisors of $N$ such that $\gcd(Q,N/Q)=\gcd(R,N/R)=1$. Then there exist natural homotopies
\begin{align*}
					w_1\simeq \id \colon										&  \TMF\to \TMF\\
      					w_Q\circ w_Q\simeq \psi^Q\colon						 	& \TMF_0(N)\to \TMF_0(N)\\
					w_Q\circ w_R\simeq w_R\circ w_Q \colon						& \TMF_0(N)\to \TMF_0(N)\\
					w_Q\circ \psi^k\simeq \psi^k\circ w_Q\colon 					& \TMF_0(N)[\frac{1}{k}]\to \TMF_0(N)[\frac{1}{k}]
\end{align*}
of morphisms of $\E_\infty$-rings.
\end{prop}

\begin{proof}
As in the proof of \Cref{stableadamscompintext}, we will simply check some homotopies in $\Isog$, and leave it to the reader to apply $\O^\top$ and obtain the homotopies above. The first homotopy is clear as if $Q=1$ we quotient by a group of order $1$, which is uninteresting. For the second homotopy, use the notation of \Cref{stableaktinlehnerdefinition} and consider the composition
\[(\M_0(N),\EE_0(N))\xrightarrow{(w_Q,\phi_Q)}(\M_0(N),\EE_0(N))\xrightarrow{(w_Q,\phi_Q)}(\M_0(N),\EE_0(N))\]
inside $\Isog$. This is naturally equivalent to $(\id,[Q])$ using the definition of composition in $\Isog$. Indeed, the composition $w_Q\circ w_Q$ is naturally equivalent to the identity by construction and the fact that the dual of a dual isogeny is naturally equivalent to the original isogeny. For the maps of elliptic curves over $\M_0(N)$, first note that in the composition
\[\EE_0(N)\xrightarrow{\phi_Q} \EE_0(N)/\H_Q\simeq w^\ast_Q \EE_0(N)\xrightarrow{w_Q^\ast \phi_Q} \EE_0(N)\]
we identify $w_Q^\ast\phi$ with $\widecheck{\phi}_Q$. By \cite[Th.2.6.1]{km}, the above composite is then identified with $[Q]$. For the third homotopy, if $Q=R$ we are done by the previous case, so we maybe assume $Q\neq R$ which from our assumptions imply that $\gcd(Q,R)=1$. In this case, the universal cyclic subgroup $\H_0(N)$ of $\EE_0(N)$ of order $N$ splits uniquely into the product $\H_Q\times \H_R\times \H_M$ where $N=QRM$. Using this, we quickly see that $w_Q\circ w_R$ and $w_R\circ w_Q$ are naturally equivalent as morphisms of stacks, so we are left with the diagram
\[\begin{tikzcd}
{\EE_0(N)}\ar[r]\ar[d]		&	{\EE_0(N)/\H_Q}\ar[d]	\\
{\EE_0(N)/\H_R}\ar[r]	&	{\EE_0(N)/(\H_Q\times H_R)}
\end{tikzcd}\]
of elliptic curves. The maps above are all the suggestive quotients and it easily follows that the above diagram commutes, giving us our third homotopy. For the fourth and final homotopy, consider the diagram
\[\begin{tikzcd}
{(\M_0(N),\EE_0(N))}\ar[r, "{(\id,[k])}"]\ar[d, "{(w_Q,\phi_Q)}"]	&	{(\M_0(N),\EE_0(N))}\ar[d, "{(w_Q,\phi_Q)}"]	\\
{(\M_0(N),\EE_0(N))}\ar[r, "{(\id,[k])}"]\					&	{(\M_0(N),\EE_0(N))}
\end{tikzcd}\]
in $\Isog$. The maps of stacks tautologically commute. The morphisms of elliptic curves also commute, as quotient maps are homomorphisms of elliptic curves, and hence commute with multiplication maps.
\end{proof}

The following statement is proven using the same techniques found in the proofs of \Cref{stableadamscompintext,stablealcompintext}. There are extra subtlies due to the presence of spans and the existence of $\T_n$ on $\TMF_0(N)$. The latter is constructed using \Cref{generalisedheckeoperators} and \Cref{descentdataexample}.

\begin{prop}\label{stableheckeandothers}
Let $k$ be an integer, $n, N$ be positive integers with $\gcd(n,N)=1$, and $Q$ be a divisor of $N$ such that $\gcd(Q,N/Q)=1$. Then there exist natural homotopies
\begin{align}
					\T_1\simeq \id\colon 										&	\TMF\to \TMF\\
					\T_n\circ \psi^k\simeq \psi^k\circ \T_n\colon 					&	\TMF[\frac{1}{kn}]\to \TMF[\frac{1}{kn}]\\
					\T_n\circ w_Q\simeq w_Q\circ \T_n\colon 						&	\TMF_0(N)[\frac{1}{n}]\to \TMF_0(N)[\frac{1}{n}]
\end{align}
of morphisms of spectra.
\end{prop}

For the following proof, recall that $\Isog^\op$ embeds into $\Span_\fin(\Isog)$ by sending a morphism $A\to B$ to the span $B\gets A\xrightarrow{=}A$. This means we can define stable Adams operations and Atkin--Lehner involutions on $\TMF$ via $\OO^\top$, the only difference is that once we have forgotten these maps of spectra were once morphisms of $\E_\infty$-rings. 

\begin{proof}
The first homotopy is a consequence of the tautological equivalence $\M_1\simeq \M$. For the second homotopy, consider the following two compositions in $\Span_\fin(\Isog)$
\begin{equation}\label{firstoftwospans}
\M\xleftarrow{(q,Q)} \M_n\xrightarrow{(p,\can)} \M\xleftarrow{(\id,[m])} \M \xrightarrow{(\id,\id)} \M
\end{equation}
\begin{equation}\label{secondoftwospans}
\M\xleftarrow{(\id,[m])} \M \xrightarrow{(\id,\id)} \M\xleftarrow{(q,Q)} \M_n\xrightarrow{(p,\can)} \M
\end{equation}
where we have suppressed the elliptic curves in the notation of our objects and implicitly inverted $kn$. The composition in $\Span_\fin(\Isog)$ is performed by taking fibre products, and the objects in the fibre product of the middle cospans of both (\ref{firstoftwospans}) and (\ref{secondoftwospans}) are naturally equivalent to $(\M_n,\EE_n)$. Focusing on the morphisms of stacks for a moment, both (\ref{firstoftwospans}) and (\ref{secondoftwospans}) compose to naturally equivalent spans as in both cases we are pulling back along the identity. The right legs of the composites of both (\ref{firstoftwospans}) and (\ref{secondoftwospans}) on elliptic curves are the canonical map identification $\EE_n\simeq p^\ast \EE$. To compare the two left legs on elliptic curves, we must show that the quotient map $\EE_n\to \EE_n/\H_n$ of $\EE_n$ by its universal subgroup commutes with multiplication by $[k]$---however, this is always true, as morphisms of elliptic curves are homomorphisms of abelian group schemes and commute with multiplication maps. The third homotopy exists for similar reasons---the isogenies in this case commute as their degrees are coprime, just as part 3 of \Cref{stablealcompintext}.
\end{proof}

The last family of natural homotopies of \Cref{interrationsintro} are between $\T_n\circ \T_m$ and $\T_m\circ \T_n$.

\begin{prop}\label{heckecommutingstatement}
Let $m,n$ be coprime positive integers. Then there is the natural homotopy
\[\T_n\circ \T_m\simeq \T_m\circ \T_n\colon \TMF[\frac{1}{mn}]\to \TMF[\frac{1}{mn}]\]
of morphisms of spectra.
\end{prop}

\begin{proof}
This is a consequence of \Cref{naturalsplittings} which is independently proven in \Cref{pprimaryfactorisationsection}.
\end{proof}

We now have a proof of \Cref{interrationsintro}.

\begin{proof}[Proof of \Cref{interrationsintro}]
Combine \Cref{stableadamscompintext,stablealcompintext,stableheckeandothers,heckecommutingstatement}.
\end{proof}

%%%%%%%%%%%%%%%%%%%%%%%%%%%%%%%%%%%%%%%%%%%%%%%%%%%%%%%%%%

\subsection{Decomposing the stacks $\M_n$}\label{decompiingstackssection}

To better analyse our stable Hecke operators and prepare ourselves for the proof of \Cref{heckeoperatorscompositionintro}, we want to study the stacks $\M_n$ more closely. 

\begin{mydef}\label{stacksforheckedecomposition}
Let $A,B$ be finite abelian groups and $d,e,m,n\geq 1$ be positive integers. 
\begin{enumerate}
\item Let $\M_A$ denote the \emph{moduli stack of elliptic curves with a chosen subgroup $H$ of type $A$} over $\Z[\frac{1}{|A|}]$. In particular, if $A=C_n$ then $\M_A=\M_0(n)$ (see the paragraph preceding \Cref{stableaktinlehnerdefinition}). We abbreviate $\M_{C_m\times C_n}$ to $\M_{(m,n)}$.
\item Let $\M_n$ denote the \emph{moduli stack of elliptic curves with a chosen finite subgroup $H$ of order $n$} over $\Z[\frac{1}{n}]$. 
\item Let $\M_{A\leq B}$ denote the \emph{moduli stack of elliptic curves with a chosen pair of nested subgroups $K\leq H$ of type $A$ and $B$, respectively}, over $\Z[\frac{1}{|A|\cdot |B|}]$. For $A=C_d\times C_m$ and $B=C_e\times C_n$ we will make the abbreviation
\[\M_{C_d\times C_m\leq C_e\times C_n}=\M_{(d,m)\leq (e,n)}.\]
\item Let $\M_{m\leq n}$ denote the \emph{moduli stack of elliptic curves with a chosen pair of nested finite subgroups $K\leq H$ of order $m$ and $n$, respectively}, over $\Z[\frac{1}{mn}]$.
\end{enumerate}
These stacks come with canonical structure maps $p$ to $\M$, and those of the form $\M_A$ and $\M_n$ also come with maps $q$ to $\M$ which take a quotient by the given subgroup.
\end{mydef}

In the definition above, we leave open the possibility that these moduli stacks are empty, for example, if $A$ is not isomorphic to a finite subgroup of $S^1\times S^1$ in part 1, or if $m\nmid n$ in part 4.\\

Let us explore some properties of these stacks. The first follows from their definition.

\begin{prop}\label{obvioussplitting}
Let $m,n\geq 1$ be positive integers. Then there are the natural decompositions of stacks
\[\M_n\simeq \coprod_{|A|=n}\M_A\qquad \M_{m\leq n}\simeq \coprod_{\substack{|A|=m \\ |B|=n}}\M_{A\leq B}\]
the first over $\Z[\frac{1}{n}]$ and second over $\Z[\frac{1}{mn}]$.
\end{prop}

Our second statement shows these stacks live in $\Isog$.

\begin{prop}\label{propertiesoflevelnstacks}
All of the stacks of \Cref{stacksforheckedecomposition} are separated Deligne--Mumford stacks of finite presentation over $\Spec \Z$. All of the maps $p$ and $q$ to $\M$ are finite \'{e}tale.
\end{prop}

\begin{proof}
Fix two finite abelian groups $A$ and $B$. By \cite[Pr.1.6.4]{km}, we see that for a fixed scheme $S$ and an $S$-point $E$ inside $\M_{\Z[\frac{1}{|A|}]}(S)$, the moduli scheme of $A$-structures on $E$ is finite \'{e}tale (although potentially empty) over $S$. Therefore the structure map $p\colon\M_A\to \M_{\Z[\frac{1}{|A|}]}$ is finite \'{e}tale and $\M_A$ is a Deligne--Mumford stack. The same argument used to prove \cite[Pr.1.6.4]{km} shows that the moduli schemes of $(A\leq B)$-structures on elliptic curves are finite \'{e}tale over $S$, hence the structure map $p\colon \M_{A\leq B}\to \M_{\Z[\frac{1}{|A|\cdot|B|}]}$ is finite \'{e}tale and $\M_{A\leq B}$ is also Deligne--Mumford. This covers the stacks and structure maps in parts 1 and 3 of \Cref{stacksforheckedecomposition} and parts 2 and 4 follow by the decompositions of \Cref{obvioussplitting}. To see the quotient morphisms are finite \'{e}tale, note there exists an involution of stacks
\[\tau\colon \M_A\to \M_A\qquad (E,H)\mapsto (E/H,\widecheck{H})\]
and that the quotient map $q\colon\M_A\to \M_{\Z[\frac{1}{|A|}]}$ is a composition of $\tau$ with the finite \'{e}tale projection $p$. In more detail, we have written $\widecheck{H}$ for the kernel of the isogeny $E/H\to E$ dual to the quotient isogeny, which can be identified with the Cartier dual of $H$ by \cite[(2.8.2.1)]{km}. Notice that $\widecheck{H}$ has the same type as $H$ as are working where $|A|$ is invertible. We also see that $\tau^2$ is naturally equivalent to the identity on $\M_A$, as the composition of isogenies $E\to E/A\to E$ is equal to the $|A|$-fold multiplication map (\cite[Th.2.6.1]{km}) which induces a natural equivalence $E/E[|A|]\simeq E$, and the double dual of $H$ is naturally equivalent to $H$.
\end{proof}

Our third statement shows how one can shave off common divisors from the stacks $\M_{(m,n)}$.

\begin{prop}\label{shavingoffdivisors}
Let $k,m,n$ all be integers. The morphism of stacks over $\Z[\frac{1}{kmn}]$
\[\M_{(km,kn)}\to \M_{(m,n)}\qquad (E,H)\mapsto(E,[k]H)\]
is an equivalence whose inverse is given by the following map:
\[ \M_{(m,n)}\to \M_{(km,kn)}\qquad (E,K)\mapsto (E,[k]^\ast K)\]
\end{prop}

To be clear: $[k]H$ is the image of $H$ in $E$ under the $k$-fold multiplication map $[k]$, and $[k]^\ast K$ is the pullback of a subgroup $K\leq E$ along $[k]$.

\begin{proof}
To justify that the second map is well-defined, note that when given a pair $(E,K)$ inside $\M_{(m,n)}(S)$ for a scheme $S$, we have the commutative diagram of schemes
\[\begin{tikzcd}
{E[k]}\ar[d]\ar[r]	&	{[k]^\ast K}\ar[d, "{[k]}"]\ar[r]	&	{E}\ar[d, "{[k]}"]	\\
{S}\ar[r]		&	{K}\ar[r]				&	{E.}
\end{tikzcd}\]
The right and outer rectangles are Cartesian by definition, hence the left square is also Cartesian. This leads to the short exact sequence of finite \'{e}tale group schemes over $S$
\[0\to E[k]\to [k]^\ast K\xrightarrow{[k]} K\to 0\]
occurring inside $E[kmn]$. Suppose for a moment that $S$ is connected. Choosing a geometric point $s\colon \Spec \kappa\to S$, consider the equivalence of categories between finite \'{e}tale commutative group schemes over $S$ and finite abelian groups with an action of $\pi_1^\et(S,s)$; see \cite[\href{https://stacks.math.columbia.edu/tag/03VD}{03VD}]{stacks} or \cite{sgaone}. From this equivalence, we see that $[k]^\ast K$ must have type $C_{km}\times C_{kn}$, as that is the unique subgroup of $C_{kmn}\times C_{kmn}$ whose image under the $k$-fold multiplication map is $C_m\times C_n$. For a general scheme $S$ we may apply this argument over each connected component, and we see that $[k]^\ast K$ has the correct type.\\

The above argument shows that $[k][k]^\ast K$ is isomorphic to $K$, given a point $(E,K)$ in $\M_{(m,n)}(S)$ for any scheme $S$. Conversely, take a point $(E,H)$ in $\M_{(km,kn)}(S)$ for a general scheme $S$. The finite \'{e}tale subscheme $[k]^\ast[k] H$ of $E$ is defined by the Cartesian diagram 
\[\begin{tikzcd}
{[k]^\ast[k] H}\ar[r]\ar[d]	&	{[k]H}\ar[d]	\\
{E}\ar[r, "{[k]}"]			&	{E}
\end{tikzcd}\]
of schemes over $S$. There is also a natural commutative diagram of schemes
\begin{equation}\label{littlediagramtokeepmumhappy}\begin{tikzcd}
{H}\ar[r, "{[k]}"]\ar[d]	&	{[k]H}\ar[d]\ar[r]	&	{S}\ar[d]	\\
{E}\ar[r, "{[k]}"]	&	{E}\ar[r]		&	{E/[k]H}
\end{tikzcd}\end{equation}
so it suffices to see the left square above is Cartesian. The right square is Cartesian by inspection, so it suffices to show the whole rectangle above is Cartesian. This leads us to the natural commutative diagram of schemes over $S$
\begin{equation}\label{biggiediagram}\begin{tikzcd}
{H}\ar[r, "{[k]}"]\ar[d]		&	{[k]H}\ar[d]		&			\\
{E[kmn]}\ar[r, "{[k]}"]\ar[d]	&	{E[mn]}\ar[r]\ar[d]	&	{S}\ar[d]	\\
{E}\ar[r, "{[k]}"]			&	{E}\ar[r, "{[mn]}"]	&	{E}
\end{tikzcd}\end{equation}
which we enlarge a little for later use. The lower-right square and lower rectangle are Cartesian by definition, hence the lower-left square is also Cartesian. The upper-left square is also Cartesian, which can be checked using the equivalence of categories between finite \'{e}tale commutative group schemes over $S$ and finite abelian groups with a $\pi_1^\et(S,s)$-action, as discussed earlier in this proof. In total, we see the desired left rectangle is Cartesian, hence $[k]^\ast[k]H$ is naturally equivalent to $H$.
\end{proof}

Our fourth statement shows the stacks $\M_A$ can be identified with certain $\M_0(n)$.

\begin{prop}\label{moregammazeronstructures}
Let $A$ be a finite abelian group and write $|A|=n$. If $A$ is not isomorphic to a subgroup of $C_n\times C_n$, then $\M_A=\varnothing$, and otherwise there exists a unique positive integer $d$ such that $A$ is isomorphic to $C_d\times C_{\frac{n}{d}}$ and $d^2|n$. In this case, the equivalence of \Cref{shavingoffdivisors} takes the form
\[\M_A\simeq \M_{(d,\frac{n}{d})}\xrightarrow{\simeq} \M_{(1,\frac{n}{d^2})}=\M_0\left(\frac{n}{d^2}\right).\]
\end{prop}

\begin{proof}
The first statement is clear. For the second statement, it is clear there exists an $d$ such that $A\simeq C_d\times C_\frac{n}{d}$. The minimal such $d$ has the property that $\gcd(d,\frac{n}{d})=d$ which is equivalent to the condition that $d^2|n$. The last statement is a special case of \Cref{shavingoffdivisors}.
\end{proof}

Our fifth statement utilises the decompositions and identifications made above to find relationships between various alternative stable Hecke operators (see \Cref{generalisedheckeoperators}), which we will now make explicit. Let us use the following notation:

\begin{itemize}
\item $\M_\Ga$ for any of the stacks of \Cref{stacksforheckedecomposition} defined by a single subgroup and write $n$ for the size of $\Ga$.
\item $p_\Ga,q_\Ga\colon \M_\Ga\to \M_{\Spec \Z[\frac{1}{n}]}$ for the structure and quotient maps associated with $\M_\Ga$, respectively.
\item $\EE_\Ga$ for the universal elliptic curve over $\M_\Ga$, which is canonically isomorphic to $p_\Ga^\ast \EE$.
\item $Q_\Ga\colon \EE_\Ga\to \EE_\Ga/\H_\Ga$ for the quotient by the universal subgroup.
\end{itemize}

\begin{mydef}
For $\Ga$ and $n$ as above, define the \emph{stable Hecke operator}
\[\T_\Ga\colon \TMF[\frac{1}{n}]\to \TMF[\frac{1}{n}]\]
by applying $\OO^\top$ to the morphism
\[(\M_{\Z[\frac{1}{n}]},\EE)\xleftarrow{(q_\Ga,Q_\Ga)} (\M_\Ga,\EE_\Ga)\xrightarrow{(p_\Ga,\can)} (\M_{\Z[\frac{1}{n}]},\EE)\]
inside $\Span_\fin(\Isog)$.
\end{mydef}

\begin{prop}\label{relationsbetweenalternativehecekeoperators}
Let $k,m,n$ be positive integers. Then there exists the natural homotopies
\[\sum_{d^2|n}\T_{(d,\frac{n}{d})}\simeq \T_n\qquad \T_{(km,kn)}\simeq \psi^k\circ \T_{(m,n)}\]
between endomorphisms of the spectra $\TMF[\frac{1}{n}]$ and $\TMF[\frac{1}{kmn}]$, respectively.
\end{prop}

\begin{proof}
The first statement follows directly from \Cref{moregammazeronstructures}. The second statement more-or-less follows from \Cref{shavingoffdivisors}. Consider the diagram
\begin{equation}\label{importmentannouncement}\begin{tikzcd}
	&& {\M_{(km,kn)}} \\
	& {\M_{(m,n)}} && {\M_{(m,n)}} \\
	\M &&{\M_{(m,n)}}&& \M.
	\arrow["{(q_2,Q_2)}"', bend right = 30, from=1-3, to=3-1]
	\arrow["{(f,[k])}"{description}, from=1-3, to=2-2]
	\arrow["{(\id, [k])}"{description}, from=3-3, to=2-2]
	\arrow["{(\id, \id)}"{description}, from=3-3, to=2-4]
	\arrow["{(f, \can)}"{description}, from=1-3, to=3-3]
	\arrow["{(p_2,\can)}", bend left = 30, from=1-3, to=3-5]
	\arrow["{(f, \can)}"{description}, from=1-3, to=2-4]
	\arrow["{(q_1,Q_1)}", from=2-2, to=3-1]
	\arrow["{(p_1,\can)}"', from=2-4, to=3-5]
\end{tikzcd}\end{equation}
in $\Isog$. Above we have suppressed the universal elliptic curves from our notation and abbreviated $\M_{\Z[\frac{1}{kmn}]}$ as $\M$, the maps $(q_i,Q_i)$ and $(p_i,\can)$ are the obvious maps, and $f$ is the equivalence of \Cref{shavingoffdivisors}. The two middle triangles tautologically commute. It is also tautological that the right-most region commutes, so we are left with the left-most region. This commutes by inspection:
\begin{itemize}
\item First, take a point $(E,H)$ inside $\M_{(km,kn)}(S)$ for some fixed scheme $S$. Applying $q_1\circ f$ to this pair yields $E/[k]H$. Consider the solid diagram of finite \'{e}tale commutative group schemes over $S$
\[\begin{tikzcd}
{K}\ar[r]\ar[d, dashed]	&	{H}\ar[r, "{[k]}"]\ar[d]	&	{[k]H}\ar[d]	\\
{E[k]}\ar[r]			&	{E[kmn]}\ar[r, "{[k]}"]	&	{E[mn]}
\end{tikzcd}\]
In the proof of \Cref{shavingoffdivisors}, we showed all the squares in (\ref{biggiediagram}) are Cartesian which also shows the right square above is Cartesian. This fact, and the fact the above rows above are short exact sequences, shows that the dashed morphism exists and is an isomorphism. In particular, $[k]H$ fits into the exact sequence of finite \'{e}tale commutative groups schemes over $S$
\[0\to E[k]\to H\xrightarrow{[k]} [k]H\to 0\]
which provides an isomorphism $[k]H\simeq H/E[k]$ compatible with the identification of $E$ with $E/E[k]$. Using this, we obtain the first natural isomorphism
\[E/[k]H\xleftarrow{\simeq,[k]} \frac{E/E[k]}{H/E[k]}\simeq E/H\]
and the second natural isomorphism above is one of the usual numbered isomorphism theorems from algebra. This yields a canonical isomorphism $q_1\circ f\simeq q_2$.
\item Next, we consider the maps of universal elliptic curves. This is comparing the quotient $Q_2\colon\EE\to \EE/\H$ of the universal elliptic curve over $\M_{(km,kn)}$ by the universal subgroup $\H$ with the composite $\EE\xrightarrow{[k]}\EE\to \EE/[k]\H$.
As in the first point above, we naturally identify this composite with $Q_2$, which shows the morphisms of universal elliptic curves are also equivalent.
\end{itemize}
Unwinding the definitions, the two paragraphs above imply the left-most region of (\ref{importmentannouncement}) commutes up to natural equivalence. The fact that the vertical arrow $(f,\can)$ in (\ref{importmentannouncement}) is an equivalence by \Cref{shavingoffdivisors} shows the span defining $\T_{(km,kn)}$ is naturally equivalent to the lower span in (\ref{importmentannouncement}). The same argument as in the proof of \Cref{stableheckeandothers} for the fact that $\T_n$ and $\psi^k$ commute shows this second span is naturally equivalent to the composition of spans defining $\psi^k\circ T_{(m,n)}$. This finishes the proof.
\end{proof}

%%%%%%%%%%%%%%%%%%%%%%%%%%%%%%%%%%%%%%%%%%%%%%%%%%%%%%%%%%

\subsection{Hecke composition formula}\label{heckecompsection}

The most complex composition relation that we explore here is a version of the Hecke composition formula. We know from \Cref{naturalsplittings} that the stable Hecke operators $\T_n$ and $\T_m$ all commute up to natural homotopy when $m$ and $n$ are coprime (\Cref{naturalsplittings}), but there is a more fundamental expression for the composition $\T_n\circ \T_m$, provided we work over $\TMF[\frac{1}{mn\phi}]$; see \Cref{heckeoperatorscompositionintro}. Our proof will still follow the same outline as the proofs appearing in \Cref{interrelationssection}, however, the stacks involved are more intricate and involve a more careful analysis. The methodology we will use below can also be used to prove the classical Hecke composition formula (\Cref{classicalfromthisargument}) and appears to have been largely overlooked in the literature, possibly as it is unnecessary in the classical case.\\

It is clear how to start the proof of \Cref{heckeoperatorscompositionintro}: we write down the two spans defining $\T_m$ and $\T_n$ and then study their composite in $\Span_\fin (\Isog)$. This composite involves a certain pullback of stacks, which we will now decompose.

\begin{prop}\label{splittingoffibreproduct}
Let $m,n$ be positive integers. Then there exists a Cartesian diagram
\[\begin{tikzcd}
{\coprod_{d,e} \M_{(d,\frac{m}{d})\leq(e,\frac{mn}{e})}}\ar[r]\ar[d]	&	{\M_m}\ar[d, "{q}"]	\\
{\M_n}\ar[r, "{p}"]										&	{\M_{\Z[\frac{1}{mn}]}}
\end{tikzcd}\]
of stacks over $\Spec \Z[\frac{1}{mn}]$. The above coproduct ranges over those positive integers $d,e$ such that $d^2|m$, $e^2|mn$, $d|e$, and either $m|de$ or $e|dn$.
\end{prop}

\begin{proof}
Define the stack $F_{m,n}$ over $\Z[\frac{1}{mn}]$ as the pullback in the following diagram:
\[\begin{tikzcd}
{F_{m,n}}\ar[r, "{q'}"]\ar[d, "{p'}"]	&	{{\M}_n}\ar[d, "{p}"]	\\
{{\M}_m}\ar[r, "{q}"]	&	{{\M}_{\Z[\frac{1}{mn}]}.}
\end{tikzcd}\]
This has a modular interpretation: for a fixed scheme $S$, we can identify the groupoid $F_{m,n}(S)$ with that of pentuples
\[(E_m,H_m,E_n,H_n,\al)\]
where $E_m,E_n$ are elliptic curves over $S$ with finite closed subgroups $H_m\leq E_m$ and $H_n\leq E_n$ of order $m$ and $n$, respectively, and $\al\colon E_n\simeq E_m/H_m$ is an isomorphism of elliptic curves over $S$. Given such a pentuple in $F_{m,n}(S)$, one can consider the commutative diagram
\[\begin{tikzcd}
{H_m}\ar[r]\ar[d]	&	{Q_m^\ast H_n}\ar[r]\ar[d]	&	{E_m}\ar[d, "{Q_m}"]	\\
{S}\ar[r]		&	{H_n}\ar[r]				&	{E_n\overset{\alpha}{\simeq} E_m/H_m}
\end{tikzcd}\]
 of schemes over $S$. The right square above is Cartesian by construction, the whole rectangle is Cartesian by inspection, so we see the left square is also Cartesian. This left square then witnesses the short exact sequence
\[0\to H_m\to \pi^\ast H_n\to H_n\to 0\]
 of finite \'{e}tale groups schemes over $S$. From the above short exact sequence, it is clear that $\pi^\ast H_n$ has order $mn$. These observations justify the well-definedness of the map of stacks
\[F_{m,n}\to \M_{m\leq mn},\qquad (E_m,H_m, E_n,H_n, \al)\mapsto (E_m, H_m\leq Q_m^\ast H_n)\]
 over $\Spec \Z[\frac{1}{mn}]$. We claim this map is an equivalence of stacks, which is easy to check using the following inverse:
\[\M_{m\leq mn}\to F_{m,n}\qquad (E,K\leq H)\mapsto (E, K, E/K, H/K, \id).\]
Indeed, all of the homotopies used to show these two functors are inverse to each other are canonical. By \Cref{obvioussplitting}, the stack $\M_{m\leq mn}$ decomposes as
\[\M_{m\leq mn}\simeq \coprod_{\substack{|A|=m \\ |B|=mn}}\M_{A\leq B}\]
where $A\leq C_m\times C_m$ and $B\leq C_{mn}\times C_{mn}$. Elementary group theory (see the proof of \Cref{moregammazeronstructures}, for example) states that for each $A$ (resp.\ $B$) there exists a unique positive integer $d$ (resp.\ $e$) such that $d^2|m$ and $A\simeq C_d\times C_\frac{m}{d}$ (resp.\ $e^2|mn$ and $B\simeq C_e\times C_{\frac{mn}{e}}$). Notice that $d|e$ and either $e|nd$ or $m|de$ by part 1 of \Cref{beautyofsubgroupnumbers}, which gives us the desired indexing set for our coproduct.
\end{proof}

With the above decomposition of stacks, we are now ready to prove our Hecke composition formula.

\begin{proof}[Proof of \Cref{heckeoperatorscompositionintro}]
First, consider the diagram of stacks over $\Z[\frac{1}{mn}]$
\begin{equation}\label{pullbacksofspansequiations}\begin{tikzcd}
	& {\coprod \M'_{e}} &&{\coprod \M'_{e}}\\
	&& {\coprod \M'_{d,e}} \\
	& {\M_n} && {\M_m} \\
	{\M} && {\M} && {\M}
	\arrow["{q_n}", from=3-2, to=4-1]
	\arrow["{p_n}"', from=3-2, to=4-3]
	\arrow["{q_m'}", from=2-3, to=3-2]
	\arrow["{p_n'}"', from=2-3, to=3-4]
	\arrow["{q_m}", from=3-4, to=4-3]
	\arrow["{p_m}"', from=3-4, to=4-5]
	\arrow["p", from=2-3, to=1-2, swap]	
	\arrow["p", from=2-3, to=1-4]
	\arrow["{q_{(e,\frac{mn}{e})}}"', bend left = -30, from=1-2, to=4-1]
	\arrow["{p_{(e,\frac{mn}{e})}}", bend right = -30, from=1-4, to=4-5]
\end{tikzcd}\end{equation}
where we have written $\M_{(d,\frac{m}{d})\leq (e,\frac{mn}{e})}$ as $\M'_{d,e}$, $\M_{(e,\frac{mn}{e})}$ as $\M'_e$, and $p$ for the structure map $\M'_{d,e}\to \M'_e$ sending an $S$ point $(E,K\leq H)$ to $(E,H)$. Moreover, the three coproducts are indexed as in \Cref{splittingoffibreproduct} and $mn$ is implicitly inverted everywhere. To see this diagram commutes (up to equivalence), we only need to check the left and right regions commute---the centre diamond commutes by \Cref{splittingoffibreproduct}. For the right region, we quickly see that given a triple $(E,K\leq H)$ in $\M'_{d,e}(S)$ for a fixed base scheme $S$, then both composites send this triple to $E$ inside $\M(S)$. For the left region, the upper composition takes a triple $(E,K\leq H)$ first to $(E,H)$ and then to $E/H$ inside $\M(S)$. Alternatively, the lower composite first takes $(E,K\leq H)$ to the pair $(E/K,H/K)$, this follows from the expression of $q'_m$ given in the proof of \Cref{splittingoffibreproduct}, and then to the quotient $(E/K)/(H/K)$. This is naturally isomorphic to $E/H$, as desired. We can enhance (\ref{pullbacksofspansequiations}) to a diagram in $\Isog$: simply pair each stack above with the associated universal elliptic curve, and pair each map labelled with a ``$p$'' with a canonical equivalence of elliptic curves, and each map labelled with a ``$q$'' with the suggestive quotient (as in \Cref{stableheckedefinition} for $\T_n$ for example). We can now apply $\OO^\top$ to this diagram, meaning that all of the maps going to the right are realised by transfer maps. This yields the following diagram of spectra
\begin{equation}\label{gleichaberfuertmf}\begin{tikzcd}
	& {\bigoplus TMF'_e} &  & {\bigoplus TMF'_e} \\
	&& {\bigoplus TMF'_{d,e} } \\
	& {TMF_n} && {TMF_m} \\
	TMF && TMF && TMF
	\arrow["{q_n^\ast}"', from=4-1, to=3-2]
	\arrow["{p_!^n}"', from=3-2, to=4-3]
	\arrow["{q_m^\ast}"', from=4-3, to=3-4]
	\arrow["{p^m_!}"', from=3-4, to=4-5]
	\arrow["{q_{(e,\frac{mn}{e})}^\ast}", bend left = 20, from=4-1, to=1-2]
	\arrow["{p^{(e,\frac{mn}{e})}_!}", bend left = 20, from=1-4, to=4-5]
	\arrow["{p^\ast}"', from=1-2, to=2-3]
	\arrow["{q'^\ast_m}"', from=3-2, to=2-3]
	\arrow["{p_!}"', from=2-3, to=1-4]
	\arrow["{(p'_n)_!}"', from=2-3, to=3-4]
	\arrow["{\T_n}", from=4-1, to=4-3, swap]
	\arrow["{\T_m}", from=4-3, to=4-5, swap]
	\arrow["{\bigoplus c_{m,n}(e,d)}", from=1-2, to=1-4]
\end{tikzcd}\end{equation}
where we have used analogous notation to (\ref{pullbacksofspansequiations}) and suppressed inverting $mn$. As $\OO^\top$ is a functor on span categories, we see the centre diagram, as well as the left and right regions of the above diagram naturally commute. The lower two triangles commute by definition. We are left with the composition $p_!\circ p^\ast$ which is \textbf{not} the identity---we claim it is non-canonically homotopic to the map $\TMF'_e\to \TMF'_e$ given by the degree of $p$ on each summand. Indeed, by \Cref{moregammazeronstructures} we obtain the following equivalence of $\E_\infty$-rings:
\[\TMF'_e=\TMF_{(e,\frac{mn}{e})}\simeq\TMF_0\left(\frac{mn}{e^2}\right)\]
If we now invert $\gcd(6,\phi(mn))$, then $\gcd(6,\phi(\frac{mn}{e^2}))$ is also invertible. Indeed, this is a consequence of the standard fact that if $a|b$, then $\phi(a)|\phi(b)$, which itself follows by reduction to the case $a=\ell^e$ and $b=\ell^{e+f}$ for a prime $\ell$, which is checked explicitly:
\[\phi(b)=\phi(\ell^{e+f})=\ell^{e+f}-\ell^{e+f-1}=\ell^f(\ell^e-\ell^{e-1})=\ell^f\phi(\ell^e)=\ell^f\phi(a)\]
We now find ourselves in the situation of part 2 of \Cref{whentoapplyadditionmaplemma}, so \Cref{compositionsareconstants} applies to give a homotopy witnessing the commutativity of the top triangle above.\\

Let us write $c_{m,n}(d,e)$ for this local degree, so the degree of $p$ between $\M'_{d,e}\to \M'_e$. The argument above leads to the first non-canonical homotopy between endomorphisms of the spectrum $\TMF[\frac{1}{mn}]$
\[\T_m\circ \T_n\approx \sum c_{m,n}(d,e)\T_{(e,\frac{mn}{e})}\simeq \sum_{\substack{b|m,n \\ a^2 | \frac{mn}{b^2}}} b \T_{(ab,\frac{mn}{ab})}\]
where the first sum is indexed as the sums in \Cref{gleichaberfuertmf} are, the second natural homotopy comes from the combinatorial bookkeeping in part 4 of \Cref{beautyofsubgroupnumbers} by setting $X^e=\T_{(e,\frac{mn}{e})}$. We have used the notation $\approx$ to denote a non-canonical homotopy. Using the natural homotopies provided in \Cref{relationsbetweenalternativehecekeoperators} and the homotopies from above, we obtain the desired homotopy
\[\T_m\circ \T_n\approx \sum c_{m,n}(d,e)\T_{(e,\frac{mn}{e})}\simeq\sum_{\substack{b|m,n \\ a^2 | \frac{mn}{b^2}}} b \T_{(ab,\frac{mn}{ab})}\simeq \sum_{\substack{b|m,n \\ a^2 | \frac{mn}{b^2}}} b \psi^b \T_{(a,\frac{mn}{ab^2})}\simeq \sum_{b|m,n} b \psi^b \T_{\frac{mn}{b^2}}\]
between endomorphisms of spectra. Once again, the homotopy indicated by the symbol $\approx$ is the only non-canonical homotopy.
\end{proof}

\begin{remark}\label{classicalfromthisargument}
If one runs the above proof, and instead of taking applying $\O^\top$ to (\ref{pullbacksofspansequiations}) one considers sheaf cohomology and the classical pullbacks, transfers, and twisting maps, then the same series of deductions above proof the classical Hecke composition formula on weight $k$ meromorphic modular forms:
\[\T_m^\cl\circ \T_n^\cl=\sum_{d|m,n}d^{k-1}\T_{\frac{mn}{d^2}}^\cl\] 
The statement for holomorphic modular forms also follows as they inject into the ring of meromorphic modular forms. Alternatively, one can use \Cref{comparisonintro} to obtain the classical formula from \Cref{heckeoperatorscompositionintro}, however, as we will prove \Cref{comparisonintro} using \Cref{heckeoperatorscompositionintro}, this is circular as stated.
\end{remark}

%%%%%%%%%%%%%%%%%%%%%%%%%%%%%%%%%%%%%%%%%%%%%%%%%%%%%%%%%%

\subsection{Factorisations into $p$-primary components}\label{pprimaryfactorisationsection}

In this subsection, we will use a natural extension of the ideas occurring in \Cref{interrelationssection,decompiingstackssection,heckecompsection} to prove \Cref{naturalsplittings}, which naturally factorises our stable Adams operations and Hecke operators into $p$-primary pieces. In fact, we have even more control over the higher homotopies between the possible orders of composing these operations than a higher category theorist might expect---this is highlighted by the following definition.

\begin{mydef}
Let $S$ be a set, $1\leq r\leq \infty$, and write $\E_r$ for the little $r$-cubes $\infty$-operad given in \cite[Df.5.1.0.4]{haname}---only the $r=1,\infty$-cases will explicitly interest us in this article. Define the $\E_r$-monoid $\ff^{\E_r}_S$ as the free $\E_r$-monoid generated by $S$ and $\ff_S$ as the free (discrete) commutative monoid generated by $S$. Suppose now that $S$ is a finite set, $C$ is an $\infty$-category, and $f_s$ is a collection of endomorphisms of some object $X$ inside $C$. If the functor $B\ff^{\E_1}_S\to C$, defined by sending the unique object of $B\ff^{\E_1}_S$ to $X$ and the generators $s\in \ff_S$ to $f_s$, has a lift through the natural map $B\ff^{\E_1}_S\to B\ff^{\E_r}_S$ then we will write
\[\bigcomp^{\E_r}_{s\in S} f_s\colon X\to X\]
for the value of this lift $B\ff^{\E_r}_S\to C$ at the map $s_1\circ\cdots \circ s_n$, where $S=\{s_1,\ldots, s_n\}$ is some enumeration of $S$. If the functor $B\ff^{\E_1}_S\to C$ admits a lift through the natural map $B\ff^{\E_1}_S\to B\ff^{\E_{\infty}}_S\to B\ff_S$, where $\ff^{\E_{\infty}}_S\to \ff_S$ is the truncation map, then we write
\[\bigcomp_{s\in S} f_s\colon X\to X\]
for the value of this lift $B\ff_S\to C$ at the map $s_1\circ\cdots \circ s_n$. If $S=\{\ell|n\}$ is the set of primes dividing $n$ for some $n$, write $\ff^{\E_r}_n$ for $\ff^{\E_r}_S$ and $\ff_n$ for $\ff_S$.
\end{mydef}

The point of the notation $\bigcomp^{\E_r}_{s\in S} f_s$ is to suggest that given two enumerations of $S$, there is a natural homotopy from the composition defined by one enumeration to the other, which is coherent in an ``$\E_r$-way''.\\ %One can also look at the map of $\E_1$-monoids induced by a functor $B\ff^{\E_r}_S\to C$
%\[\ff^{\E_r}_S=\Map_{B\ff^{\E_r}_S}(\ast,\ast)\to \Map_C(X,X)\]
%and view $\bigcomp^{\E_r}_{s\in S}f_s$ as a collection of naturally equivalent points in $\Map_C(X,X)$. In other words, the map of $\E_1$-monoids above factors through the $\E_r$-centre of $\Map_C(X,X)$.\\

As the natural map $\ff^{\E_{r}}_S\to \ff_S$ is a morphisms of $\E_r$-monoids  for all $1\leq r\leq \infty$, we see that if $\bigcomp$ is well-defined then so is $\bigcomp^{\E_r}$. We believe the notion $\bigcomp^{\E_r}$ is the natural definition in higher category theory, however, the decompositions of our operations will use the stronger version $\bigcomp$---perhaps this further highlights their exceptional formal properties. Let us start with the easier case of the stable Adams operations.

\begin{prop}\label{adamsoperationspslitting}
Given an integer $k$, there is the following natural homotopy of morphisms of $\E_\infty$-rings:
\[\psi^k\simeq \bigcomp_{\substack{\ell|k \\ \ell\,\mathrm{prime}}}\psi^{\ell^{\nu_\ell(k)}}\colon \TMF[\frac{1}{k}]\to \TMF[\frac{1}{k}]\]
\end{prop}

Throughout the following proof (and the next), we will write integers $k$ is their prime components $k=\prod_{\ell|k}\ell^{e}$, where $e=\nu_{\ell}(k)$.

\begin{proof}
First, use the natural homotopy of $\E_\infty$-rings $\psi^k\simeq \psi^{|k|}$ of \Cref{stableadamscompintext} to reduce ourselves to the case where $k$ is strictly positive. Define a functor $F\colon B\ff_k\to \Isog^\op$ by sending $\ast$ to $(\M_{\Z[\frac{1}{k}]},\EE)$, $1$-morphisms $a=\prod \ell^e$ written in their prime decomposition to $(\id,[\prod \ell^{e \nu_\ell(a)}])$, and where the $2$-functorality comes from the fact that for $a=\prod \ell^e$ and $b=\prod \ell^f$ in $\ff_k$, we have the natural equivalences of morphisms
\[(\id,[\prod \ell^{f\nu_\ell(b)}])\circ (\id,[\prod \ell^{e\nu_\ell(a)}])\simeq (\id,[\prod \ell^{ef\nu_\ell(ab)}])\]
in $\Isog$. We define the desired functor of $\infty$-categories as the composite
\[B\ff_k\xrightarrow{F} \Isog^\op \xrightarrow{\O^\top}\CAlg.\]
This allows us to define $\bigcomp_{\ell|k} \psi^{\ell^{\nu_\ell(k)}}$, and to compare this to $\psi^k$ we refer to \Cref{stableadamscompintext}.
\end{proof}

The situation for our stable Hecke operators is analogous, except, as per usual, we need to be more careful with our span $\infty$-categories.

\begin{prop}\label{heckeoperatorssplitting}
Given a positive integer $n$, there is the following natural homotopy of morphisms of spectra:
\[\T_n\simeq \bigcomp_{\substack{\ell|n \\ \ell\,\mathrm{prime}}}\T_{\ell^{\nu_\ell(n)}}\colon \TMF[\frac{1}{n}]\to \TMF[\frac{1}{n}]\]
\end{prop}

Our proof is rather long, however, it is mostly formal using \Cref{extensiontheoremintro} and category theory.

\begin{proof}
Let us implicitly work with $n$ inverted. To define a functor $\TT_n\colon B\ff_n\to \Sp$ such that $\TT_n(\ast)=\TMF$ and $\TT_n(\{\ell\})=\T_{\ell^{\nu_\ell(n)}}$, we will use the functor $\OO^\top$ and reduce ourselves to constructing a functor $B\ff_n\to \Span_\fin(\Isog)$ which sends $\ast$ to $(\M,\EE)$ and $\{\ell\}$ to the span
\[(\M,\EE) \xleftarrow{(q,Q)} (\M_{\ell^{\nu_\ell(n)}}, \EE_{\ell^{\nu_\ell(n)}}) \xrightarrow{(p,\can)} (\M,\EE)\]
 of \Cref{stableheckedefinition}. To simplify this situation further, let us use that the functor $\Span$ from the $\infty$-category of adequate triples to the $\infty$-category of $\infty$-categories admits a left adjoint: the twisted arrow category functor $\Tw(-)$ (with ingressive and egressive morphisms defined as those which induce equivalences on the target and source, respectively); see \cite[Th.A]{spanadjunction}. Let us unwrap this adjunction as a sanity check. Given a functor $F\colon C\to \Span D$ of $\infty$-categories, then its adjunct $G\colon \Tw(C)\to D$ is given on objects by sending $f\colon c\to c'$ to the object $P_f$ in the span in the image of $f$ under $F$:
\[F_f=\left(F_c\gets P_f\to F_{c'}\right).\]
On 1-morphisms, $G$ is given by considering a commuting square
\[\begin{tikzcd}
{c}\ar[r, "{\al}"]\ar[d, "{f}"]	&	{d}\ar[d, "{g}"]	\\
{c'}					&	{d'}\ar[l, "{\be}", swap]
\end{tikzcd}\]
in $C$, first applying $F$ to this square
\[\begin{tikzcd}
	&&& {P_f}\ar[d, equal]\ar[ddddlll, bend right = 40]\ar[ddddrrr, bend left = 40] \\
	&&& {P_{\alpha,g,\beta}} \\
	&& {P_{\alpha,g}} && {P_{g,\beta}} \\
	& {P_\alpha} && {P_g} && {P_\beta} \\
	{F_c} && {F_d} && {F_{d'}} && {F_{c'}}
	\arrow[from=2-4, to=3-3]
	\arrow[from=2-4, to=3-5]
	\arrow[from=3-5, to=4-4]
	\arrow[from=3-3, to=4-4]
	\arrow[from=3-3, to=4-2]
	\arrow[from=4-2, to=5-1]
	\arrow[from=4-2, to=5-3]
	\arrow[from=4-4, to=5-3]
	\arrow[from=4-4, to=5-5]
	\arrow[from=3-5, to=4-6]
	\arrow[from=4-6, to=5-5]
	\arrow[from=4-6, to=5-7]
\end{tikzcd}\]
and then defining $G$ of $(\al,\be)\colon f\to g$ to be either of the two naturally equivalent composites $P_f\to P_g$ above. Using this translation, we would like to construct a functor $\TT_n'\colon \Tw(B\ff_n)\to \Isog$ of adequate triples as follows: 

\begin{itemize}
\item The generators $\ell\in \ff_n$, viewed as objects in $\Tw(B\ff_n)$, are sent to
\begin{equation}\label{functorongenobjects}(\M_{\ell^{\nu_\ell(n)}},\EE_{\ell^{\nu_\ell(n)}})=:(\M(\ell),\EE(\ell)).\end{equation}
\item The two 1-morphisms of $\Tw(B\ff_n)$
\[\begin{tikzcd}
{\ast}\ar[r, equal]\ar[d, "{\ell}"]	&	{\ast}\ar[d, equal]	\\
{\ast}					&	{\ast}\ar[l, "{\ell}", swap]
\end{tikzcd}\qquad
\begin{tikzcd}
{\ast}\ar[r, "{\ell}"]\ar[d, "{\ell}"]	&	{\ast}\ar[d, equal]	\\
{\ast}						&	{\ast}\ar[l, equal]
\end{tikzcd}\]
in $\Tw(B\ff_n)$ are sent to the morphisms
\[(\M(\ell),\EE(\ell))\xrightarrow{(q,Q)}(\M,\EE)\qquad (\M(\ell),\EE(\ell))\xrightarrow{(p,\can)}(\M,\EE)\]
in $\Isog$, respectively.
\end{itemize}

Note that $\M(\ell)$ should \textbf{not} be confused with $\Ga(\ell)$-level structure---it is simply notation for this proof. As $\Tw(B\ff_n)$ is a 1-category and $\Isog$ is a 2-category, constructing our desired functor $\TT_n'$ is a 2-categorical task; see \cite[\href{https://kerodon.net/tag/009P}{009P}]{kerodon} for more on the $2$-categorical nerve.\\

Given two stacks $\xx,\yy$ which each admit two maps $p,q$ into $\M$, we will write $\xx\widetilde{\times} \yy$ for the pullback of the following span:
\[\xx\xrightarrow{p_\xx} \M\xleftarrow{q_\yy} \yy\]
Notice this definition is not symmetric and that $\xx\widetilde{\times}\yy$ admits maps $p,q$, defined by $p_\yy\circ \pi_2$ and $q_\xx\circ \pi_1$, respectively. Note that the $q$-map of $X\widetilde{\times} Y$ uses the $q$-map of $X$, and \emph{a priori} is independent of the $q$-map of $Y$. Let us define $\TT'_n$ on objects $\ell$ generating $\ff_n$ using (\ref{functorongenobjects}), and for an element of the form $\ell^e$ in $\ff_n$ for some $e\geq 2$, we define $\TT'_n(\ell^r)$ inductively as
\[\M(\ell^r)=\M(\ell^{e-1})\widetilde{\times}\M(\ell)\]
where $\M(\ell)$ has maps $p,q$ from \Cref{stacksforheckedecomposition}. The pasting lemma implies that we could have also glued $\M(\ell)$ on the other side. For example, when $e=3$ we have the diagram of stacks
\begin{equation}\label{definitionofmellpowers}\begin{tikzcd}
	&&& {\M(\ell^3)} \\
	&& {\M(\ell^2)} && {\M(\ell^2)} \\
	& {\M(\ell)} && {\M(\ell)} && {\M(\ell)} \\
	\M && \M && \M && \M
	\arrow[from=3-2, to=4-1]
	\arrow[from=3-2, to=4-3]
	\arrow[from=2-3, to=3-2]
	\arrow[from=3-4, to=4-3]
	\arrow[from=2-3, to=3-4]
	\arrow[from=1-4, to=2-3]
	\arrow[from=1-4, to=2-5]
	\arrow[from=2-5, to=3-4]
	\arrow[from=3-4, to=4-5]
	\arrow[from=2-5, to=3-6]
	\arrow[from=3-6, to=4-5]
	\arrow[from=3-6, to=4-7]
\end{tikzcd}\end{equation}
where all squares are Cartesian, all arrows with positive slope are (pullbacks of) the quotient maps $q$, and those with negative slope are (pullbacks of) the structure maps $p$. For a composite element $a$ inside $\ff_n$, write the prime decomposition of $a=\ell_1^{e_1}\cdots \ell_r^{e_r}$ where $\ell_i<\ell_{i+1}$, $r\geq 2$, and $e_i\geq 1$, and define $\TT'_n(a)$ inductively as
\[\M(a)=\M(\ell_1^{e_1}\cdots \ell_{r-1}^{e_{r-1}})\widetilde{\times}\M(\ell_r^{e_r})\]
with $p,q$-maps inductively inherited from the $\widetilde{\times}$ construction; some discussion of this definition appears in \Cref{littlediscucsisonondefiniion}. We will use the notation $\EE(a)$ for $p^\ast\EE$ for the universal elliptic curve over $\M(a)$. Importantly, notice that $\M(\ell_1)\widetilde{\times}\M(\ell_2)$ and $\M(\ell_2)\widetilde{\times}\M(\ell_1)$ are both isomorphic as stacks over $\M$. Indeed, by \Cref{splittingoffibreproduct} and its proof, these stacks are naturally isomorphic to
\[\M(\ell_1)\widetilde{\times}\M(\ell_2)\simeq \M(\ell_2)\widetilde{\times}\M(\ell_1)\simeq \M_{\ell_1^{\nu_{\ell_1}(n)}\ell_2^{\nu_{\ell_2}(n)}}\] 
as finite \'{e}tale commutative group schemes naturally decompose into coprime components. This inductively yields a collection of specified isomorphisms
\[\sigma_\ast\colon \M(a)\simeq \M(\ell_{\sigma(1)}^{e_{\sigma(1)}})\widetilde{\times}\cdots\widetilde{\times}\M(\ell_{\sigma(r)}^{e_{\sigma(1)}})\]
for any $\sigma\in \Sigma_r$, the symmetric group on $r$-letters. Let us move on to morphisms: given
\begin{equation}\label{morphismsidagram}\begin{tikzcd}
{\ast}\ar[r, "{c}"]\ar[d, "{a}"]	&	{\ast}\ar[d, "{b}"]	\\
{\ast}					&	{\ast}\ar[l, "{d}"]
\end{tikzcd}\end{equation}
a morphism in $\Tw(B\ff_n)$, which implies the relation $a=dbc$ inside $\ff_n$, we want to produce a morphism $\TT^c_d\colon \M(a)\to \M(b)$. If $c=\ell^e$ and $d=1$, so $b=a/\ell^e$, then we define $\TT^{\ell^e}_1$ on stacks as the composite
\[\M(a)\xrightarrow{\sigma,\simeq} \M(\ell^{\nu_\ell(a)})\widetilde{\times} \M(a/\ell^{\nu_\ell(a)})\xrightarrow{p(\ell^e)\widetilde{\times}\id} \M(\ell^{\nu_\ell(a)-e})\widetilde{\times} \M(\ell^{\nu_\ell(a)-e})\xrightarrow{\sigma^{-1},\simeq} \M(a/\ell^e)\]
where $\ell=\ell_i$ in the prime decomposition of $a$, $\sigma$ is the transposition switching $i$ and $1$, and $p(\ell^e)$ is the indicated base change of the $p$-map $\M(\ell^e)\to \M$, arising from the diagram defining $\M(\ell^{\nu_\ell(a)})$. The associated morphism on universal elliptic curves is the canonical isomorphism.\\

If $c=1$ and $d=\ell^e$, then $\TT^1_{\ell^e}$ is defined as the composite
\[\M(a)\xrightarrow{\tau,\simeq} \M(a/\ell^{\nu_\ell(a)})\widetilde{\times} \M(\ell^{\nu_\ell(a)})\xrightarrow{q(\ell^e)} \M(a/\ell^{\nu_\ell(a)})\widetilde{\times} \M(\ell^{\nu_\ell(a)-e})\xrightarrow{\tau^{-1},\simeq} \M(a/\ell^e)\]
where $\ell=\ell_i$ in the prime decomposition of $a$, $\tau$ is the transposition switching $i$ and $r$, and $q(\ell^e)$ is the indicated base change of the $q$-map $\M(\ell^e)\to \M$, as given in the definition of $\M(\ell^{\nu_\ell(a)})$. The associated morphism on universal elliptic curves is the indicated base change the canonical quotient map $Q\colon \EE(\ell^e)\to q^\ast \EE$.\\

For general $c,d$, we write a prime decomposition for $c=\ell_1^{e_1}\cdots \ell_r^{e_r}$ and $d=\ell_1^{f_1}\cdots \ell_r^{f_r}$ with $\ell_i<\ell_{i+1}$ and define $\TT^c_d$ as the following composite:
\begin{equation}\label{definitionofmap}\TT_{\ell_r^{f_r}}^1\circ\cdots\circ \TT_{\ell_1^{f_1}}^1\circ \TT^{\ell_r^{e_r}}_1\circ\cdots\circ \TT^{\ell_1^{e_1}}_1\end{equation}
In other words, first do the $\ell$-primary structure maps $\TT^{\ell^e}_1$, from small primes to large, followed by the $\ell$-primary quotient maps $\TT^1_{\ell^e}$, from small primes to large.\\

To see this assignment is functorial up to natural isomorphism, consider the composition of morphisms
\[\begin{tikzcd}
	\ast & \ast & \ast \\
	\ast & \ast & \ast
	\arrow["a", from=1-1, to=2-1]
	\arrow["b", from=1-2, to=2-2]
	\arrow["c", from=1-3, to=2-3]
	\arrow["d", from=1-1, to=1-2]
	\arrow["e", from=1-2, to=1-3]
	\arrow["f", from=2-3, to=2-2]
	\arrow["g", from=2-2, to=2-1]
\end{tikzcd}\]
in $\Tw(B\ff_n)$, yielding the relations $a=fbd$ and $b=gce$. Our goal is to construct natural $2$-morphisms comparing $\TT^e_f\circ \TT^d_g$ with $\TT^{ed}_{gf}$ as morphisms from $\M(a)$ to $\M(c)$. Writing these maps out using (\ref{definitionofmap}), we see that it suffices to construct the following compatible $2$-isomorphisms:
\[\TT^{\ell_1^e}_1\circ \TT_1^{\ell_2^{e'}}\simeq \TT^{\ell_2^{e'}}_1\circ\TT^{\ell_1^e}_1 \qquad \TT_{\ell_1^f}^1\circ \TT^1_{\ell_2^{f'}}\simeq \TT_{\ell_2^{f'}}^1\circ\TT_{\ell_1^f}^1\qquad \TT^{\ell_1^e}_1\circ \TT^1_{\ell_2^f}\simeq \TT_{\ell_2^f}^1\circ\TT^{\ell_1^e}_1\]
By construction, each of the maps above further decomposes into the case where to the case where the exponents $e=e'=f=f'=1$. Moreover, as these morphisms are defined via base change of a collection of generating morphisms, it suffices to construct the following collection of $2$-morphisms in $\Isog$:
\begin{equation}\label{firstcaseofcompositions}\TT^{\ell_1}_1\circ \TT_1^{\ell_2}\simeq \TT^{\ell_2}_1\circ\TT^{\ell_1}_1\colon (\M(\ell_1\ell_2),\EE(\ell_1\ell_2))\to (\M,\EE)\end{equation}
\begin{equation}\label{secondcaseofcompositions}\TT_{\ell_1}^1\circ \TT^1_{\ell_2}\simeq \TT_{\ell_2}^1\circ\TT_{\ell_1}^1\colon (\M(\ell_1\ell_2),\EE(\ell_1\ell_2))\to (\M,\EE)\end{equation}
\begin{equation}\label{thirdcaseofcompositions}\TT^{\ell_1}_1\circ \TT^1_{\ell_2}\simeq \TT_{\ell_2}^1\circ\TT^{\ell_1}_1\colon (\M(\ell_1\ell_2),\EE(\ell_1\ell_2))\to (\M,\EE)\end{equation}
For the first, consider the diagram of stacks
\[\begin{tikzcd}
	&	&	{\M_{b\leq ab}}\ar[rd]\ar[lldd, "{\Phi}", bend right=30, swap]	&	\\
	&	{\M(\ell_1\ell_2)}\ar[rr, "{\TT^{\ell_1}_1}", swap]\ar[dd, "{\TT^{\ell_2}_1}"]\ar[ru, "{\simeq}"]\ar[ld, "{\simeq}"]	&&	{\M_b}\ar[dd, "{\TT^{\ell_2}_1}"]	\\
{\M_{a\leq ab}}\ar[rd]	&&&	\\
	&	{\M_a}\ar[rr, "{\TT_1^{\ell_1}}"]	&&	{\M}
\end{tikzcd}\]
where $a=\ell_1^{\nu_{\ell_1}(n)}$ and $b=\ell_2^{\nu_{\ell_2}(n)}$, we have used \Cref{splittingoffibreproduct} to identify various twisted products, and $\Phi$ is the isomorphism of stacks sending a triple $(E,K\leq H)$ to $(E, H[a]\leq H)$---this is an isomorphism as for such a triple $H$ has a canonical decomposition $H\simeq H[a]\oplus H[b]$ as $a$ and $b$ are coprime. It is now clear that the underlying diagram of stacks commutes up to natural isomorphism, as a triple $(E,K\leq H)$ in $\M_{b\leq ab}(S)$ for fixed scheme $S$ is sent to $E$ in either direction. The two morphisms of elliptic curves are equivalent as they are both naturally equivalent to the canonical identification of $\EE_{b\leq ab}$ with the pullback of $\EE$ along the structure map for $\M_{b\leq ab}$.\\

Replacing the two occurrences of $\TT_1^{\ell_2}$ above with $\TT_{\ell_2}^1$, we obtain a similar diagram, which commutes as this time a triple $(E,K\leq H)$ is sent to $E/K$ up to natural isomorphism following the diagram in either direction and the morphisms of elliptic curves are both equivalent to the canonical quotient $\EE_{b\leq ab}\to \EE_{b\leq ab}/\H_b$ by the universal subgroup of order $b$. Finally, replacing both $\TT_1^{\ell_i}$ above with $\TT_{\ell_i}^1$ for $i=1,2$, we obtain another diagram. This commutes as a triple $(E,K\leq H)$ is sent to $E/H$ up to natural isomorphism travelling around the diagram in both directions, and the two morphisms of elliptic curves are naturally equivalent to the quotient map $\EE_{b\leq ab}\to \EE_{b\leq ab}/\H_{ab}$ by the universal subgroup of order $ab$---for this last point, it is crucial that $a$ and $b$ are coprime, so that this quotient map be performed by killing the summands of $\H_{ab}\simeq \H_a\oplus \H_b$ in any order, as in the proofs of \Cref{stablealcompintext,stableheckeandothers}. This process yields the desired natural isomorphisms (\ref{firstcaseofcompositions})--(\ref{thirdcaseofcompositions}) above.\\

This yields our $2$-functor $\TT'_n\colon \Tw(B\ff_n)\to \Isog$ with the desired properties. To check the induced functor of $\infty$-categories is a functor of adequate triples, we should show that $\TT'_n$ preserves ingressive and egressive morphisms. The egressive morphisms in $\Isog$ are all morphisms, and the ingressive morphisms are those in $\fin$, so we are required to check that $\TT'_n$ sends ingressive morphisms in $\Tw(B\ff_n)$ to morphisms in $\fin$. By \cite[Ex.2.9]{spanadjunction}, we see the egressive morphisms of $\Tw(B\ff_n)$ are those (\ref{morphismsidagram}) with $d=1$, and the ingressive are those with $c=1$. In particular, if $d=1$ then $\TT^c_1$ lies in $\fin$. Indeed, if $c=\ell$ is a prime, this is clear by construction and base change, and for general $c$, this holds by base change.\\

In summary, we obtain a functor $\TT_n\colon B\ff_n\to \Sp$ allowing us to define the notation
\[\bigcomp_{\ell| n}\T_{\ell^{\nu_\ell(n)}}\colon \TMF\to \TMF\]
to indicate a collection of naturally homotopic morphisms. The definitions above were made so that there is a natural identification of stacks $\M(\prod_{\ell|n}\ell)\simeq \M_n$ which also preserves the evident $p$- and $q$-maps. In particular, after applying $\OO^\top$, we obtain the promised natural homotopy
\[\T_n\simeq \bigcomp_{\ell| n}\T_{\ell^{\nu_\ell(n)}}\colon \TMF\to \TMF.\qedhere\]
\end{proof}

As a reminder: the functor $\TT_n$ allowing us to define $\bigcomp_{\ell|n}\T_{\ell^{\nu_\ell(n)}}$ does \textbf{not} saying anything about $\T_\ell\circ \T_\ell$---that is the content of \Cref{heckeoperatorscompositionintro}.

\begin{remark}\label{littlediscucsisonondefiniion}
Let us justify our choice not to define $\M(a)$ as either $\M_a$ or perhaps via a non-twisted product over $\M$. Firstly, we cannot define $\M(a)$ as $\M_a$, as if we could prove \Cref{heckeoperatorssplitting} in this manner, we would show that $\T_m\circ \T_n$ is naturally homotopic to $\T_{mn}$ for all pairs $m,n$, which is not true for Hecke operators between classical modular forms, nor is it true for our stable Hecke operators; see \Cref{heckeoperatorscompositionintro}. In essence, \Cref{heckeoperatorssplitting} only refers to the commutativity of coprime stable Hecke operators and says nothing interesting about $\T_{\ell^2}$. To this end, we define $\M(a)$ as formally as possible, which originally led us to take a non-twisted product. In this case though, the author found no reason for $q$-maps to exist, hence the twisted products $\widetilde{\times}$.
\end{remark}

\begin{remark}
The stacks $\M_{a\leq ab}$ and $\M_{b\leq ab}$ appearing towards the end of the above proof and in the proof of \Cref{heckeoperatorssplitting} are abstractly equivalent for all $a,b\geq 1$---not just those coprime pairs. Indeed, the morphism $\Phi_{a,b}$ defined on $S$-points by sending a triple $(E,K\leq H)$ to $(E/H,\ker(\widehat{\pi}_{H/K})\leq \ker(\widehat{\pi}_H))$, where $\pi_{(-)}$ indicates the quotient of $E$ by the finite subgroup $(-)$, is an isomorphism with inverse $\Phi_{b,a}$. This isomorphism is not a morphism over $\M$, and it does not make the two diagrams constructed from (\ref{pullbacksofspansequiations}) by reversing the order of $m$ and $n$, commute inside $\Isog$ in any obvious way. We are currently exploring the possibility of a natural homotopy between $\T_n\circ \T_m$ and $\T_m\circ \T_n$ for general $m,n\geq 1$.
\end{remark}

Combining what we have done so far yields \Cref{naturalsplittings}.

\begin{proof}[Proof of \Cref{naturalsplittings}]
Combine \Cref{adamsoperationspslitting,heckeoperatorssplitting}.
\end{proof}

%%%%%%%%%%%%%%%%%%%%%%%%%%%%%%%%%%%%%%%%%%%%%%%%%%%%%%%%%%

\subsection{Nonexistence of integral stable Hecke operators}\label{nonexistencesection}

It is well-known that the unstable Adams operations $\psi^k$ on $\Omega^\infty \KU$ do not lift to endomorphisms of $\KU$ without inverting $n$; see \cite[\textsection II.13]{bluebook}. Here, we prove \Cref{nonexistenceintro} by providing some computational evidence for the lack of integral stable operations on $\TMF$, as some justification for the number of fractions seen in this article so far. Let us start with the nonexistence of certain stable Adams operations.

\begin{prop}\label{nonexistenceofadams}
Let $p$ be a prime. There exists a map of spectra $\psi^k\colon \TMF_p\to\TMF_p$ which agrees with the operation $\psi^k\colon \TMF[\frac{1}{k}]\to \TMF[\frac{1}{k}]$ of \Cref{adamsdefinition} homotopy groups over $\Q_p$ if and only if $p\nmid k$. 
\end{prop}

The proof below is often used to show that $\psi^k$ only act on $\KU$ with $k$ inverted, with the periodicity class $\Delta^{24}$ switched to the Bott periodicity class in $\pi_2 \KU$.

\begin{proof}
If $\psi^k\colon \TMF_{(p)}\to \TMF_{(p)}$ agrees with our operations of \Cref{adamsdefinition} on homotopy groups, then in particular we have the calculation $\psi^k(\Delta^{-24})=k^{-288}\Delta^{-24}$ inside $\pi_{-576}\TMF_{(p)}$ from \Cref{comparisonintro} and the fact that the localisation map
\[\pi_{-576}\TMF_{(p)}\to \pi_{-576}\TMF_\Q\]
is injective. This calculations makes sense if and only if $k$ is invertible in $\pi_0\TMF_{(p)}\simeq \Z_{(p)}$, in other words, if and only if $p\nmid k$.
\end{proof}

We have a similar statement regarding stable Hecke operators, which (as expected) is more subtle as these operators are not multiplicative nor is the periodicity class $\Delta^{24}$ in $\pi_{576}\TMF$ a $\T_n$-eigenform.

\begin{prop}\label{nonexistenceofhecke}
Let $p=2$ or $3$, choose a positive integer $n$ such that $p|n$, and write $r=\nu_p(n)$. Suppose that for $p=2$ then $n/p^r$ is a square, and for $p=3$ then $n/p^r$ has prime factors $p$ with exponent $e$ such that either $\ell\equiv_3 1$ and $e\equiv_6 0,1,3,4$, or $\ell\equiv_3 2$ and $e$ is even. Then there does not exist a map of spectra $\T_n\colon \TMF_{p}\to\TMF_{p}$ which agrees with the operation $\T_n\colon \TMF[\frac{1}{n}]\to \TMF[\frac{1}{n}]$ of \Cref{stableheckedefinition} on homotopy groups over $\Q_p$.
\end{prop}

In particular, there do not exist stable operations $\T_{2^e}$ (resp. $\T_{3^e}$) on $\TMF_{2}$ (resp. $\TMF_{3}$). It is possible to extend the above proposition to more values of $n$ using \Cref{congruencetheoremintext}, but to simplify our exposition, we will be content with the statement above. \\

For the following proof, we will frequently use the image of the unit map $\pi_\ast\Sph\to \pi_\ast\TMF$, often called the \emph{Hurewicz image}, as calculated in \cite{tmfbook} and highlighted in colour; proofs can be found at the prime $2$ in \cite{hurewicztmf}, the prime $3$ in \cite{tmfthree}, and also in \cite[\textsection 11.11 \& \textsection 13.7]{brunerrognes}.

\begin{proof}
Let us write $n=p^rm$ where $p\nmid m$ and $r\geq 1$. We will now work prime by prime. 

\paragraph{($p=2$ case)}	Consider $\kappa^2\in \pi_{28}\TMF_{2}\simeq (\Z/2\Z)^2$, which lies in Hurewicz image. It follows that for any potential endomorphisms of spectra $\T_n$ on $\TMF_{2}$ we have
\[\T_n(\kappa^2)=\kappa^2 \T_n(1)=\kappa^2\sigma(n)\]
as $\T_n$ is $\Sph$-linear. The multiplicativity of the divisor function we obtain
\[\sigma(n)=\sigma(2^r)\sigma(m)=(2^r+2^{r-1}+\cdots +2+1)\sigma(m)\equiv_2 \sigma(m).\]
Our assumption that $m$ is a square shows $\sigma(m)$ is odd (\Cref{nonvansihingofsigmaattwo}), hence $\T_n(\kappa^2)=\kappa^2$. A representative for $\kappa^2$ on the descent spectral sequence for $\TMF_{2}$ is $\eta[2\Delta h_2]$, where $h_2$ detects $\nu$. As there are no classes in higher filtration, we can compute $\T_n(\kappa^2)$ by computing the map induced by $\T_n$ on the $E_2$-page of this spectral sequence. As this induced map $\T_n$ is linear as a module of the Adams--Novikov spectral sequence for the sphere, and we assume that $\T_n$ agrees with $n\T_n^\cl$ on the $E_2$-page, we see that
\[\T_n(2\Delta h_2)=2\T_n(\Delta) h_2=n\tau(n)2\Delta h_2.\]
Multiplicitivity shows $\tau(n)=\tau(2^r)\tau(m)$, and it is classical that $\tau(2^r)$ is even. Indeed, the behaviour of $\tau(n)$ with respect to products, dictated by the fact $\Delta$ is a Hecke eigenform, shows that if $\tau(p)\equiv 0$ modulo $p$, then $\tau(pn)\equiv 0$ modulo $p$. One then checks that $\tau(2)=-24$ and $\tau(3)=252$ to see that $\tau(2^r)$ is even and $3$ divides $\tau(3^r)$. In particular, this argument also holds if we had assumed that $\T_n$ induces $\T_n^\cl$ on rational homotopy groups, rather than $n\T_n^\cl$. In the end, we obtain a section calculation of $\T_n(\kappa^2)$
\[\T_n(\kappa^2)=\T_n(\eta[2\Delta h_2])=\eta\T_n[2\Delta h_2]=n\tau(n) \eta[2\Delta h_2]=n\tau(n)\kappa^2=0\]
contradicting our previous computation.

\paragraph{($p=3$ case)}	Repeat the above argument, but now for the element $\be^3=\al[\al\Delta]$ in $\pi_{30}\TMF_{3}\simeq \Z/3\Z$, refering to \Cref{nonvansihingofsigmaatthree} for the nonvanishing of $\sigma(n)$ modulo $n$.
\end{proof}

Let us use the \emph{height 2 Adams summand} and \emph{height 2 image of $J$ spectrum} to extend a version of \Cref{nonexistenceofhecke} to all primes; see \cite[\textsection3]{adamsontmf}. Recall that for any prime $p$, we can give $\TMF_p$ the $\Z_p^\times$-action of the stable Adams operations from \Cref{adamsdefinition}. Concretely, this given by recognising $\TMF_p$ as $\O^\top_\BTtwo(\M\times \Spf \Z_p,\EE[p^\infty])$ using \Cref{diagramfromluriestheorem} and noting that $\Z_p^\times$ acts on $\EE[p^\infty]$ through multiplication by $p$-adic units $\lambda\in \Z_p^\times$. 

\begin{mydef}
Define the $\E_\infty$-ring $\U=\TMF_p^{h\F_p^\times}$ using the canonical inclusion $\F_p^\times\to \Z_p^\times$ called the \emph{periodic height $2$ Adams summand}. Given a topological generator $g$ of $\Z_p^\times$, define the $\E_\infty$-ring $\J^2$, called the \emph{periodic height $2$ image of $J$} as the equaliser of $\psi^g\colon \U\to \U$ and the identity.
\end{mydef}

These definitions are periodic versions of those from \cite[\textsection3]{adamsontmf}, meaning that, writing $\u$ and $j^2$ for the connective analogues of these $\E_\infty$-rings from \emph{idem}, there is a commutative diagram of $\E_\infty$-rings
\[\begin{tikzcd}
{j^2}\ar[r]\ar[d]	&	{\u}\ar[d]	\\
{J^2}\ar[r]		&	{\U}
\end{tikzcd}\]
where the vertical maps are induced by the localisation $\tmf\to \TMF$. Now we can precisely state the last piece of \Cref{nonexistenceintro}.

\begin{prop}\label{heckenonexstiencecomplicated}
Let $p$ be a prime and $e\geq 1$ a positive integer. Then there does not exist a map of spectra $\T_{p^e}\colon \TMF_p\to \TMF_p$ which is $\Z_p^\times$-equivariant and agrees with the operation $\T_{p^e}\colon \TMF[\frac{1}{p}]\to \TMF[\frac{1}{p}]$ of \Cref{stableheckedefinition} on homotopy groups over $\Q_p$.
\end{prop}

Note that by \Cref{interrationsintro} our operations $\T_n$ on $\TMF_p$ for $p\nmid n$ are $\Z_p^\times$-equivariant.

\begin{proof}
Suppose that $p\geq 5$ so that $\pi_\ast\TMF_p$ has no torsion. If such an operation $\T_{p^e}$ did exist on $\TMF_p$, then their $\Z_p^\times$-equivariance would mean that they induce morphisms $\T_{p^e}\colon \J^2\to \J^2$. By \cite[Th.F]{adamsontmf}, we know that $\j^2$ detects the whole $\al_{i/j}$ family; see \cite[\textsection5]{greenbook}. In particular, the class $\al=\al_{1/1}\in \pi_{2p-3}\Sph_p$ is $p$-torsion, and is detected by the class $\bar{E}$ in $\pi_{2p-3}\j^2$, itself the image of the Eisenstein series $E_{p-1}\in \pi_{2p-2}\u$ under the boundary map
\[\pi_{2p-2}\u\to \pi_{2p-3}\j^2\]
defining $\j^2$; see \cite[\textsection3.3]{adamsontmf}. The map $\j^2\to \J^2$ is an injection on homotopy groups as this is true for both $\tmf\to \TMF$ and $\u\to \U$, so the class $\bar{E}$ in $\pi_\ast\J^2$ also detects $\al$. This leads us to the calculation
\[\al=(1+\cdots+p^e)\al=\sigma(p^e)\al=\T_{p^e}(\al)=\T_{p^e}(\bar{E})=p^e\sigma_{p-2}(p^e)\al=0\]
inside $\pi_{2p-3}J^2\simeq \Z/p\Z$ and our desired contradiction.
\end{proof}

Together, the above results imply \Cref{nonexistenceintro}.

\begin{proof}[Proof of \Cref{nonexistenceintro}]
This is a combination of \Cref{nonexistenceofadams,nonexistenceofhecke,heckenonexstiencecomplicated}.
\end{proof}

%%%%%%%%%%%%%%%%%%%%%%%%%%%%%%%%%%%%%%%%%%%%%%%%%%%%%%%%%%
%%%%%%%%%%%%%%%%%%%%%%%%%%%%%%%%%%%%%%%%%%%%%%%%%%%%%%%%%%
%%%%%%%%%%%%%%%%%%%%%%%%%%%%%%%%%%%%%%%%%%%%%%%%%%%%%%%%%%

\section{Applications to number theory}\label{applicationssection}

To end this article, we apply the existence of the stable Hecke operators $\T_n$ on $\TMF$ with a little homotopy theory to prove \Cref{congruencetheoremintro,primethreemaedaintro} as well as their generalisations, and a general family of congruences of modular forms.\\

As a sample of what is to come, fix a positive integer $n$ not divisible by $3$, and consider the operation $\T_n$ on $\pi_{30} \TMF[\frac{1}{n}]$. From \cite{bauer}, one can calculate $\pi_{30} \TMF[\frac{1}{n}]\simeq \Z/3\Z$ generated by the image of the class $\be^3$ from $\pi_{30} \Sph$ under the unit map $\Sph\to \TMF$. As $\T_n$ is $\Sph$-linear, we immediately obtain the calculation
\[\T_n(\be^3)=\be^3\T_n(1)=\be^3 \sigma(n)\]
where the fact that $\T_n(1)=\sigma(n)$ is a consequence of \Cref{comparisonintro} and the classical calculation $\T_n^\cl(1)=\frac{\sigma(n)}{n}$. Alternatively, this is a consequence of \Cref{compositionsareconstants} and the fact that the map of stacks $\M_n\to \M_{\Z[\frac{1}{n}]}$ is finite \'{e}tale of degree $\sigma(n)$---the number of subgroups of order $n$ inside $(\Q/\Z)^2$. There is, however, another way to calculate $\T_n(\be^3)$. Appealing to \cite{bauer} again, we see that $\be^3=\al[\al\Delta]$ where $\al$ is the image of the class $\al\in \pi_3\Sph$---the $3$-local part of the stable Hopf map $\nu\colon S^7\to S^4$. By \Cref{comparisonintro} (or better yet, \Cref{etwoaggreement}), we see that $\T_n$ can be calculated on the $E_2$-page of the descent spectral sequence for $\TMF[\frac{1}{n}]$ by $n\T_n^\cl(\al[\al\Delta])$, where we now consider $\Delta\in E_2^{0,24}$ and $\al\in E_2^{1,3}$ as elements on the $E_2$-page. Again, as $\T_n$ is $\Sph$-linear, we see $\T_n(\al)=\al\T_n(1)$, hence $n\T_n^\cl$ must also be $\al$-linear. This leaves us with the expression
\[n\T_n^\cl(\al^2\Delta)=n\al^2\T_n^\cl(\Delta)=n\tau(n)\al^2\Delta\]
representing $\T_n(\be^3)$ on the $E_2$-page, as $\tau(n)$ is the $\T_n^\cl$-eigenvalue of $\Delta$. As there are no nonzero elements on the $E_\infty$-page of the descent spectral sequence for $\TMF[\frac{1}{n}]$ with filtration higher than $\al^2\Delta$, we obtain the equality $\T_n(\be^3)=n\tau(n)\be$ inside $\pi_{30}\TMF[\frac{1}{n}]$. We therefore have two expressions for $\T_n(\be^3)$ inside $\pi_{30}\TMF[\frac{1}{n}]\simeq\Z/3\Z$, in other words, we know that $\sigma(n)$ and $n\tau(n)$ agree modulo $3$---this recovers a famous congruence of Ramanujan.\\

In \Cref{congrunecessection}, we will generalise this argument to obtain a series of congruences between coefficients in the $q$-expansions of various modular forms capitalising on the torsion classes and periodicity of $\TMF$, a special case of which yields \Cref{congruencetheoremintro}. In \Cref{maedasconjecturesubsection}, we will apply these congruences to obtain some previously unknown cases of Maeda's conjecture, proving \Cref{primethreemaedaintro}.

%%%%%%%%%%%%%%%%%%%%%%%%%%%%%%%%%%%%%%%%%%%%%%%%%%%%%%%%%%

\subsection{Congruences between coefficients in $q$-expansions}\label{congrunecessection}

Let us implicitly localise at a prime $p$. For a positive integer $d$, let us write $\b_{12d}$ for the following basis of the $\Z_{(p)}$-module $\mf_{12d}$ of weight $12d$ holomorphic modular forms over $\Z_{(p)}$:
\[\b_{12d}=\{\Delta^d, c_4^3\Delta^{d-1},c_4^6\Delta^{d-2},\ldots, c_4^{3d-6}\Delta^2,c_4^{3d-3}\Delta, c_4^{3d}\}\]
Using the fact that $\TMF_{(2)}$ is $\Delta^8$-periodic and $\TMF_{(3)}$ is $\Delta^3$-periodic, we can extrapolate the example opening this section to general multiples and powers of $\Delta$.

\begin{theorem}\label{congruencetheoremintext}
Let $d$ be a nonnegative integer and $n$ a positive integer, and write $b_n^e$ for the $\Delta^e$-coefficient of $\T_n^\cl(\Delta^e)$ with respect to the basis $\b_{12e}$.
\begin{enumerate}
\item If $n$ is odd, then $\sigma(n)$ is congruent to:
\begin{enumerate}
\item $nb^{8d}_n$ and $nb_n^{8d+1}$ modulo $8$;
\item $nb^{8d+2}_n$, $nb^{8d+4}_n$, $nb^{8d+5}_n$, and $nb^{8d+6}_n$ modulo $4$; and
\item $nb^{8d+3}_n$ and $n b_n^{8d+7}$ modulo $2$.
\end{enumerate}
\item If $n$ is not divisible by $3$, then $\sigma(n)$ is congruent to $nb_n^d$ modulo $3$.
\end{enumerate}
\end{theorem}

The following proof amounts to well-known manipulations in stable homotopy theory and some familiarity with the homotopy groups of $\TMF$; see will use either \cite{bauer} or \cite{brunerrognes} as references for these calculations.

\begin{proof}
Let us start with the $p=2$ case. We will use two lines of argumentation repeatedly, so let us call them the \textbf{Hurewicz argument} and the \textbf{Toda argument}. The first two cases of (a) demonstrate these two arguments:

\paragraph{(Hurewicz argument)}	Consider the following equalities in $\pi_{192d+20}\TMF_{(2)}$:
\begin{equation}\label{hurewiczexample}\sigma(n)\bar{\kappa}\Delta^{8d}=\bar{\kappa}\Delta^{8d}\T_n(1)=\T_n(\bar{\kappa}\Delta^{8d})=\bar{\kappa}\T_n(\Delta^{8d})=\bar{\kappa} [n\T_n^\cl(\Delta^{8d})]\end{equation}
The element $\bar{\kappa}\in\pi_{20}\TMF_{(2)}$ is the image of the well-known element from $\pi_{20}\Sph_{(2)}$ and $\Delta^8\in \pi_{192}\TMF_{(2)}$ is the periodicity class. The first equality comes from the classical calculation $\T_n^\cl(1)=\frac{\sigma(n)}{n}$ and \Cref{comparisonintro}, the second and third from the facts that $\bar{\kappa}\Delta^{8d}\in \pi_{192d+20}\TMF_{(2)}$ lies in the image of the unit map $\Sph\to \TMF$, see \cite{hurewicztmf} or \cite[\textsection11.11]{brunerrognes}, also called the Hurewicz image, and $\T_n$ is $\Sph$-linear, and the fourth from \Cref{etwoaggreement} combined with the facts that $[g]$ represents $\bar{\kappa}$ in the DSS for $\TMF_{(2)}$ and that no nonzero classes live in higher filtration in the $192d$th column on the $E_\infty$-page than $g \Delta^{8d}$. As $\b_{96d}$ is also a basis for $\pi_{192d}\TMF_{(2)}$ and the only element in this basis supporting multiplication by $\bar{\kappa}$ is $\Delta^{8d}$, then (\ref{hurewiczexample}) shows that $nb^{8d}_n$ is equal to $\sigma(n)$ inside $\Z/8\Z$. This is the first case of (a).

\paragraph{(Toda argument)}	For the second case of (a), consider the Toda bracket $\langle 8, \nu, \bar{\kappa}\Delta^{8d}\rangle$ as a subset of $\pi_{192d+24}\TMF_{(2)}$. The element $\nu\in \pi_3\TMF_{(2)}$ is the image of the quaternionic Hopf map $\nu\in \pi_3\Sph_{(2)}$. We will now discuss this bracket in some detail. One can calculate the Massey product $\langle 8,h_2,g\Delta^{8d}\rangle$ on the $E_5$-page of the DSS, where $h_2$ (resp.\ $g$) represent $\nu$ (resp.\ $\bar{\kappa}$). Indeed, on the $E_5$-page we have $8h_2=0$ which yields the equality of sets
\[\langle 8,h_2,g\Delta^{8d}\rangle=\{8x\in E_5^{0,24+192d} | d_5(x)=h_2 g\Delta^{8d}\}\subseteq E_5^{0,24+192d}.\]
Bauer's calculation of the homotopy groups of $\tmf$ \cite{bauer}, in particular, the calculation of all $d_5$-differentials shows that the above set is equal to those elements
\[a\Delta^{8d+1}+a_1\Delta^{8d}c_4^3+a_2 \Delta^{8d-1}c_4^6+\cdots +a_{8d+1}c_4^{24d+3}\in E_5^{0,24+192d}\]
where $a\equiv 1$ modulo $8$. In particular, modulo $8$ and all the elements in the basis $\b_{96d+12}$ \textbf{not} of the form $\Delta^{8d+1}$, the Massey product $\langle 8,h_2,g\rangle$ is the singleton set of $8\Delta$. The Moss convergence theorem (see \cite{moss} for the original statement for the Adams spectral sequence and \cite{evahanamoss} for the adaptation to other multiplicative spectral sequences) shows our desired Toda bracket $\langle 8, \nu,\bar{\kappa}\Delta^{8d}\rangle$ is congruent to $[8\Delta^{8d+1}]$, modulo $8$ and the elements of $\b_{96d+12}-\{\Delta^{8d+1}\}$. This means the indeterminacy of this Toda bracket is generated by $8$ and the elements of $\b_{96d+12}-\{\Delta^{8d+1}\}$. We will now need the following well-known lemma; see \cite[Pr.IV.2.6(iii)]{symspectravthree}.

\begin{lemma}\label{todabracketclaim}
Let $f\colon A\to B$ be a morphism of $\E_\infty$-rings, $M$ an $A$-module, $N$ a $B$-module, and $\phi\colon M\to f_\ast N$ a morphism of $A$-modules. Given elements $x\in \pi_\ast M$ and $y,z\in \pi_\ast A$ such that $xy=0=yz$ and $\langle x,y,z\rangle$ is a Toda bracket in $\pi_\ast M$, then we have the containment
\[\phi\left(\langle x,y,z\rangle\right)\subseteq \langle \phi(x),f(y),f(z)\rangle\]
of subsets inside $\pi_\ast N$. 
\end{lemma}

Using our calculation of $\langle 8,\nu,\bar{\kappa}\Delta^{8d}\rangle$ and this lemma with $A=\Sph$, $B=M=N=\TMF_{(2)}$, $f$ is the unit $\Sph\to \TMF_{(2)}$, and $\phi=\T_n$, we obtain the containments of sets
\[\T_n([8\Delta]\Delta^{8d})\equiv \T_n(\langle 8,\nu,\bar{\kappa}\Delta^{8d}\rangle)\subseteq \langle \T_n(8),\nu,\bar{\kappa}\Delta^{8d}\rangle \subseteq \sigma(n)\langle 8,\nu,\bar{\kappa}\Delta^{8d}\rangle\equiv \sigma(n)[8\Delta]\Delta^{8d}.\]
The first and last expressions are singleton sets, so we obtain the equality $\T_n([8\Delta]\Delta^{8d})\equiv \sigma(n)[8\Delta]\Delta^{8d}$, modulo $8$ and the elements of $\b_{96d+12}-\{\Delta^{8d+1}\}$. We can also use \Cref{etwoaggreement} and the fact there are no elements of higher filtration to see $\T_n([8\Delta]\Delta^{8d})$ is represented by the class $[n\T_n^\cl(8\Delta^{8d+1})]$. In total, this shows that the $\Delta^{8d+1}$-coefficient of $\T_n^\cl(8\Delta^{8k+1})$ multiplied by $n$ is congruent to $\sigma(n)$---our desired result.

\paragraph{(Remaining $2$-local cases)}	The justifications for parts (b) and (c) follow from slight variations on the Hurewicz and Toda arguments given above, so let us only outline the differences in all the cases. First, by $\Delta^8$-periodicity, let us only consider the $d=0$ cases.\\

For (b), we use the Hurewicz argument with respect to the element $\bar{\kappa}[2\Delta^4]$ for the $b_n^4$-case, which is denoted as $\bar{\kappa}D_4$ in \cite[Df.9.22]{brunerrognes}. For $b_n^2$ and $b_n^6$ we make similar arguments with $\nu[h_2\Delta^2]$ and $\nu[h_2\Delta^6]$, respectively. These elements are denoted as $\nu\cdot\nu_2$ and $\nu\cdot\nu_6$ in \cite[Df.9.22]{brunerrognes}. Both of these elements lie in the Hurewicz image, as does $\nu$. We can only make an $E_2$-calculation for $\T_n$ on $[h_2\Delta^2]$ and $[h_2\Delta^6]$ modulo $4$ as these classes are $4$-torsion on the $E_2$-page and these classes are only well-defined up to higher filtration anyhow. %\footnote{One can try to combine the higher filtration elements into our calculations using $\MU$-based synthetic spectra (not the even variant) and analysing $\nu(\TMF)\otimes C(\tau^4)$; see \cite{syntheticspectra} for these terms and notations. Indeed, these small exotic extensions and $\Delta$ are both detected by the above synthetic spectrum. For simplicity, let us forgo this discussion and extension. Thank you to John Rognes, for introducing this argument to us.}
The last case for (b) is to apply the Toda argument to $\langle 4,e_{116}, \nu\rangle$ inside $\pi_{120}\TMF_{(20)}$, where $e_{116}$ is represented by $[2g\Delta^4]$ on the $E_2$-page of the DSS. One then calculates this Toda bracket is equal to $[8\Delta^5]$ modulo $4$ and elements in $\b_{60}$ not of the form $\Delta^5$. This yields the $b_n^5$-case.\\

For (c), we first consider the class $x=\eta\epsilon \Delta^8=\nu^3\Delta^8$. We would like to apply the Hurewicz argument here, but we will use a slight modification. Indeed, this element lies in the Hurewicz image which leads us to the equality
\[\T_n(\nu^3\Delta^8)=\sigma(n)\nu^3\Delta^8.\]
This class has representative $h_1c\Delta^8$ on the $E_2$-page of the DSS, where $h_1 c\Delta$ represents the class $q$ in homotopy, which also lies in the Hurewicz image, and $c$ represents the element $\epsilon\in \pi_8\TMF_{(2)}$ which is the image of the element somtimes written as $c_0\in \pi_8\Sph_{(2)}$. Using this, and the fact that there are no nonzero classes in filtration higher than $x$ in this column of the DSS, we obtain
\[\T_n(\nu^3\Delta^8)=[n\T_n^\cl(h_1c\Delta^8)]=[n h_1c\Delta \T_n^\cl(\Delta^7)]\]
using \Cref{etwoaggreement} and the fact that the DSS for $\TMF_{(2)}$ is linear over the ANSS for $\Sph$. On the $E_2$-page of the DSS, we see the only classes in degree $168$ and filtration zero (where $\Delta^7$ lives) that support multiplication by $c h_1\Delta$ are multiples of $\Delta^7$, leading us the equality
\[\T_n(\nu^3\Delta^8)=[nh_1c\Delta \T_n^\cl(\Delta^7)]=nb_n^7[ h_1 c\Delta \Delta^7]=nb_n^7\nu^3 \Delta^8\]
and our desired congruence $n b_n^{7}\equiv \sigma(n)$ modulo $2$. For the $b_n^3$-case, apply the Toda argument to $\langle 2,e_{70},\eta\rangle$, where $e_{70}$ is the unique nonzero class in $\pi_{70}\TMF_{(2)}$. This Toda bracket is then calculated to be equal to $[8\Delta^3]$ modulo $2$ and the elements of $\b_{36}$ not of the form $\Delta^3$.

\paragraph{(The $3$-local cases)}	In the $3$-local world, we use the Hurewicz argument with respect to the elements $\be\Delta^{3d}$ and $\al[\al\Delta]\Delta^{3d}$ to obtain the $b_n^{3d}$- and $b_n^{3d+1}$-cases, respectively. The element $\al\in \pi_3\TMF_{(3)}\simeq \Z/3\Z$ is the image the element in $\pi_3\Sph_{(3)}$ detected by the first Steenrod power in the Adams spectral sequence. The element $\be\in \pi_{10}\TMF_{(3)}$ is the image of the first element $\be\in \pi_{10}\Sph_{(3)}$ in the divided $\beta$-family. For the $b_n^{3d+2}$-case, use the slightly altered Hurewicz argument (as done in the $b_n^{7}$-case at the prime $2$) applied to the element $\al\be\Delta^{3d+3}$, as this class lies in the Hurewicz image and has $E_2$-representation $[\al\be\Delta^{3d+3}]$ and $\al\be\Delta$ represents a class in the Hurewicz image.
\end{proof}

One can push \Cref{congruencetheoremintext} a little harder and look at various cases. A few more direct examples of this can be found in \cite[Cor.10.1.6]{mythesis} which compute some explicit values of $b^e_n$. For example, 
\[b_n^2=n\left(\tau(n)^2+2\sum_{i=1}^{n-1}\tau(i)\tau(2n-i)\right)\equiv_4 \sigma(n)\]
for odd $n$.\\

We can now immediately prove \Cref{congruencetheoremintro}.

\begin{proof}[Proof of \Cref{congruencetheoremintro}]
As $\Delta$ is a normalised $\T_n^\cl$-eigenform we have $b_n^1=\tau(n)$. \Cref{congruencetheoremintext} implies that $n\tau(n)\equiv_8 \sigma(n)$ when $n$ is even, and $n\tau(n)\equiv_3 \sigma(n)$ when $n$ is not divisible by $3$.
\end{proof}

%%%%%%%%%%%%%%%%%%%%%%%%%%%%%%%%%%%%%%%%%%%%%%%%%%%%%%%%%%

\subsection{Further evidence for Maeda's conjecture}\label{maedasconjecturesubsection}

The congruences of \Cref{congruencetheoremintext} will allow us to prove \Cref{primethreemaedaintro}, which we repeat now for the reader's convenience.

\begin{theorem}\label{primetthreemaedaintext}
Let $d,n\geq 2$ be two coprime integers with $n$ not divisible by $3$ satisfying the following two conditions:
\begin{enumerate}
\item $d\leq 1,000$ and for all $1\leq i\leq d-1$, the coefficient of $q^d$ in the $q$-expansion of $\Delta^i$ is divisible by $3$.
\item For each prime factor $\ell$ of $n$ with exponent $e$, if $\ell\equiv_3 1$ then $e\equiv_6 0,1,3,4$, and if $\ell\equiv_3 2$ then $e$ is even. 
\end{enumerate}
Then $\T_{dn, 12d}(X)$ satisfies Maeda's conjecture, meaning that the characteristic polynomial of $\T_{dn}^\cl$, acting on the rational vector space of weight $12d$ cusp forms $S_{12d}$, is irreducible over $\Q$ and has full Galois group $\Sigma_{d}$.
\end{theorem}

For example, as there are infinitely many primes $\ell$ congruent to $1$ modulo $3$, we see that $\T_{d\ell,12d}(X)$ satisfies Maeda's conjecture for all such primes and all $2\leq d\leq 1000$ satisfying condition $1$ above. Using the \texttt{SAGE} code
\begin{center}\texttt{[d for d in range(1,500) if Set([i for i in range(1,d)\\
if (D\^{ }i).padded\_list (510) [d] \% 3 == 0]).cardinality() == d-1]}
\end{center}
where $\texttt{D}$ is a power series expression for the discriminant modular form $\Delta$, one can calculate the set of all $2\leq d\leq 500$ satisfying condition $1$ above as\footnote{According to \url{https://oeis.org/} this set is possibly defined by the sequence sending $2n$ to $3^n$ and $2n+1$ to $2\cdot 3^n$. Proving this sequence indeed defines the set below seems to involve some combinatorics we have not yet overcome.}
\[\{2,3,6,9,18,27,54,81,162,243,486\}.\]
According to the list of confirmed cases of Maeda's conjecture from \cite[Th.1.5]{expevimaedaconjecture}, this is the first infinite family of Hecke operators $\T_n^\cl$ where $n$ is not a prime which satisfies Maeda's conjecture for a fixed weight. This infinite family is also just a subset of all Hecke operators covered by the above theorem. We will prove the above statement using the $3$-torsion in $\pi_\ast\TMF$. Using the $2$-power torsion, we obtain the following alternative statement.

\begin{theorem}\label{primetwomaedaintext}
Let $d,n\geq 2$ be coprime integers where $n$ is odd satisfying the following two conditions:
\begin{enumerate}
\item $d\leq 1,000$ and for all $1\leq i\leq d-1$, then writing $e_i$ for the $2$-adic valuation of the coefficient of $q^d$ in the $q$-expansion of $\Delta^i$, we require that
\[e_i\geq \begin{cases}
1	&	d\equiv_8 3,7	\\
2	&	d\equiv_8 2,4,5,6	\\
3	&	d\equiv_8 0,1
\end{cases}\]
\item Setting $e=\min(e_i)$, then we require that if:
\begin{enumerate}
\item $e=1$ then $n$ is a square.
\item $e=2$ then $n$ has at most one prime factor with an odd exponent, and in this case, the prime and the exponent are congruent to $1$ modulo $4$.
\item $e\geq 3$ then $n$ satisfies the equivalent conditions in part 3 of \Cref{nonvansihingofsigmaattwo}.
\end{enumerate}
\end{enumerate}
Then $\T_{dn, 12d}(X)$ satisfies Maeda's conjecture.
\end{theorem}

For example, as there are infinitely many primes $\ell\equiv_4 1$, we have an infinite family of $\T_{d\ell,12d}(X)$ satisfying Maeda's conjecture where $2\leq d\leq 1000$ is such that conditions 1 and 2 above hold. Using a \texttt{SAGE} code similar to that following \Cref{primetthreemaedaintext}, we see the collection of such $2\leq d\leq 500$ satisfying these two conditions are precisely those in the following set:
\[\{2,4,6,8,12,16,24,32,48,64,96, 128, 192, 256, 384\}\]
To prove the above theorems we rely on two (collections of) results from the literature. First, a statement by Ahlgren.

\begin{theorem}[{\cite[Th.1.4]{maedaahlgren}}]\label{ahlgrenresult}
Write $S_d'$ for the rational subspace of $S_d$ spanned by modular forms with vanishing constant and linear terms in their $q$-expansion. Let $d$ be a positive integer such that $\dim_\Q S_d'\geq 1$ and suppose there exists an $m$ such that Maeda's conjecture holds for $\T_m^\cl$ acting on $S_d$. then for any positive integer $n\geq 2$, the following are equivalent:
\begin{enumerate}
\item Maeda's conjecture holds for $\T_n^\cl$ acting on $S_d$.
\item There exists a modular form $f\in S_d'$ with $a_n(f)\neq 0$, ie, whose $q$-expansion has nonzero $q^n$-coefficient.
\end{enumerate}
\end{theorem}

The second result is by Ghitza--McAndrew, of which we will only state a fraction of what they prove.

\begin{theorem}[{\cite[Th.1.5]{expevimaedaconjecture}}]\label{lessthantwelvethousandresult}
For a positive integer $d\leq 12,000$, Maeda's conjecture holds for $\T_2^\cl$ acting on $S_d$.
\end{theorem}

The proof outline for \Cref{primetthreemaedaintext,primetwomaedaintext} is to show that various coefficients in the $q$-expansions of $\Delta^d$ do not vanish using \Cref{congruencetheoremintext}.

\begin{proof}[Proof of \Cref{primetthreemaedaintext}]
Recall the basis $\b_{12d}-\{c_3^{3d}\}$ for $S_{12d}$; see \Cref{congrunecessection}. Notice the $q$-expansion of $c_4$ takes the form $1+240(O(q))$, meaning that modulo $3$, the coefficient for $q^m$ in the $q$-expansion of $c_4^{3d-3i}\Delta^i$ vanishes for $1\leq m\leq i-1$ and otherwise the coefficient of $q^i$ is $1$. In particular, using the above basis and condition 2 in our hypotheses, we see the $\Delta^d$-coefficient of a modular form $f\in S_{12d}$ is congruent to the $q^d$-coefficient of $f$ modulo $3$. Using the definition of the classical Hecke operators on $q$-expansions of modular forms, we have the expression for the $q^d$-coefficient in the $q$-expansion of $\T_n^\cl(\Delta^d)$
\[\sum_{e|d,n}e^{12d-1}a_{\frac{dn}{e^2}}(\Delta^d)=a_{dn}(\Delta^d)\]
as $\gcd(d,n)=1$. Using the discussion surrounding the basis above, we see that $a_{dn}(\Delta^d)\equiv b_n^d$ modulo $3$ using the notation of \Cref{congruencetheoremintext}. By this same theorem, we further see this is equivalent to $\frac{\sigma(n)}{n}$ modulo $3$, where we note that $3\nmid n$ by assumption. By \Cref{nonvansihingofsigmaatthree}, we see that $\sigma(n)$ does not vanish modulo $3$ due to our hypotheses on $n$. This implies that the coefficient of $q^{dn}$ in the $q$-expansion of $\Delta^d\in S_{12d}'$ does not vanish modulo $3$, hence it does not vanish integrally, and \Cref{ahlgrenresult} combined with \Cref{lessthantwelvethousandresult} gives the desired result.
\end{proof}

\begin{proof}[Proof of \Cref{primetwomaedaintext}]
The proof is essentially the same as the proof of \Cref{primetthreemaedaintext}: one studies the numbers $b_n^d$ modulo either $2$, $4$, or $8$, depending on the value of $e$ in the hypotheses, and refers to \Cref{nonvansihingofsigmaattwo} to determine the nonvanishing of $\sigma(n)$.
\end{proof}

\begin{remark}\label{higherprimesmaedaremark}
To proofs of both \Cref{primetthreemaedaintext,primetwomaedaintext} have clear limitations to the primes $3$ and $2$, respectively, however, we are optimistic that more can be said. At a prime $p$, for example, one could hope to study spectra such as Behrens' $Q(N)$ spectra (\Cref{behrensqnconstruction}), for some $N$ not divisible by $p$, and the image of the unit map $\Sph\to Q(N)$ on homotopy groups. Using \Cref{interrationsintro}, these spectra $Q(N)$ have stable Hecke operators, so one can hope to play the same game we played in the proof of \Cref{congruencetheoremintext} to obtain extra congruences of coefficients of various modular forms at other primes $p$, and hence extra cases of Maeda's conjecture---a version for arbitrary primes could prove fruitful to study this conjecture.
\end{remark}

%%%%%%%%%%%%%%%%%%%%%%%%%%%%%%%%%%%%%%%%%%%%%%%%%%%%%%%%%%
%%%%%%%%%%%%%%%%%%%%%%%%%%%%%%%%%%%%%%%%%%%%%%%%%%%%%%%%%%
%%%%%%%%%%%%%%%%%%%%%%%%%%%%%%%%%%%%%%%%%%%%%%%%%%%%%%%%%%

\appendix
\section{Combinatorics}\label{combinatorialappendix}

In the main body of this article, there are a few instances where we use simple combinatorial results related to the stacks $\M_A$ and $\M_{A\leq B}$ of \Cref{stacksforheckedecomposition}, and the morphisms between them. None of the results below are new or deep, but we also could not find other clean references.\\

First, we verify the divisibility of instances of the divisor function $\sigma(n)$ at the primes $2$ and $3$. Let us write this down at the primes $2$ and $3$ separately, although the statements and proofs are similar.

\begin{prop}\label{nonvansihingofsigmaatthree}
Let $n$ be a positive integer. If $3\nmid n$, then $\sigma(n)\not\equiv_3 0$ if and only if for every prime factor $\ell$ of $n$ with exponent $e$, if $\ell\equiv_3 1$ then $e\equiv_6 0,1,3,4$, and if $\ell\equiv_3 2$ then $e$ is even.
\end{prop}

\begin{proof}
Recall the divisor function $\sigma(n)$ is multiplicative, meaning that writing $n=\prod \ell^e$ in its prime factorisation, we have the following equality:
\[\sigma(n)=\prod \sigma(\ell^e)\]
From this, we see $\sigma(n)$ does not vanish modulo $3$ if and only if every $\sigma(\ell^e)$ does not vanish modulo $3$. Writing out $\sigma(\ell^e)$ and applying Fermat's little theorem yields the following congruences:
\[\sigma(\ell^e)=\ell^e+\ell^{e-1}+\cdots+\ell+1\equiv_3\begin{cases}
(\ell+1)\frac{e+1}{2}	&	\text{for odd }e\\
(\ell+1)\frac{e}{2}+1	&	\text{for even }e
\end{cases}\]
If $e$ is odd then we need both $e+1\not\equiv_6 0$ and $\ell+1\not\equiv_3 0$, ie, if $e\equiv_6 1,3$ then $\ell\equiv_3 1$. If $e$ is even we have two cases: for $\ell\equiv_3 2$ there are no further restrictions on $e$, but if $\ell\equiv_3 1$ we further require $2\frac{e}{2}+1\not\equiv_3 0$ so $e\not\equiv_3 1$. Collecting these observations, we see that for the nonvanishing of $\sigma(n)$, if $\ell\equiv_3 1$ then demand that $e\equiv_6 0,1,3$, or $4$, and if $\ell\equiv_3 2$ then $e$ can be any even exponent.
\end{proof}

\begin{prop}\label{nonvansihingofsigmaattwo}
Let $n$ be a positive odd integer.
\begin{enumerate}
\item Then $\sigma(n)\not\equiv_2 0$ if and only if $n$ is a square.
\item Then $\sigma(n)\not\equiv_4 0$ if and only if there is at most one prime factor $\ell$ of $n$ whose exponent $e$ is odd, and for such a singular $\ell$ we demand that both $\ell$ and $e$ are congruent to $1$ modulo $4$.
\item Then $\sigma(n)\not\equiv_8 0$ if and only if exactly one of the following two statements is true:
\begin{enumerate}
\item There are at most two prime factors of $n$ whose exponents are odd, and in this case, we demand such primes and their exponents are congruent to $1$ modulo $4$.
\item There is at most one prime factor $\ell$ of $n$ with odd exponent $e$ such that either $e\equiv_8 1$ and $\ell\equiv_8 3$, or $e\equiv_8 3$ and $\ell\equiv_4 1$, or $e\equiv_8 5$ and $\ell\equiv_8 7$.
\end{enumerate}
\end{enumerate}
\end{prop}

\begin{proof}
For part 1, consider the following congruences arising from Fermat's little theorem:
\[\sigma(\ell^e)=\ell^e+\ell^{e-1}+\cdots+\ell+1\equiv_2 (e+1)\]
It follows that for $\sigma(n)\not\equiv_2 0$, we need all prime factors of $n$ to have an even exponent, ie, we need $n$ to be a square. For part 2, we need all the prime factors of $n$ to have even exponent except precisely one prime $\ell$ with exponent $e$, where we allow $\sigma(\ell^e)\equiv_4 2$. To check this condition, suppose $e$ is odd, in which case Euler's theorem provides us with the following congruences:
\[\sigma(\ell^e)=1+\ell+\ell^2+\cdots + \ell^e\equiv_4 (1+\ell)\frac{e+1}{2}\]
From this, we must have $\ell\equiv_4 1$, else the right-hand side vanishes, and in this case, we need $\frac{e+1}{2}$ to be odd, which yields the restriction $e\equiv_4 1$ as well. This yields part 2. For part 3, we may allow $n$ to have either at most two prime factors $\ell$ with exponent $e$ satisfying $\sigma(\ell^e)\equiv_4 2$, or a single prime $\ell$ satisfying $\sigma(\ell^e)\equiv_8 4$. We already know how to characterise the first case, so let us move on to the second. Suppose then that $e$ is odd, then Euler's theorem again yields the congruences
\[\sigma(\ell^e)=\ell^e+\ell^{e-1}+\cdots+\ell+1\equiv_8 
\begin{cases}
\frac{e-1}{4}(\ell^3+\ell^2+\ell+1)+1+\ell		&	e\equiv_4 1	\\
\frac{e+1}{4}(\ell^3+\ell^2+\ell+1)			&	e\equiv_4 3	
\end{cases}\]
which we now analyse. Note the sum $\ell^3+\ell^2+\ell+1$ is always divisible by $4$. If $e\equiv_4 1$, then we consider two further cases: if $e\equiv_8 1$, then $\frac{e-1}{4}$ is even and we demand $\ell+1\equiv_8 4$ so $\ell\equiv_8 3$; if $e\equiv_8 5$, then for similar reasons we demand $\ell\equiv_8 7$. If $e\equiv_4 3$, then we need $\frac{e+1}{4}$ to be odd so we demand $e\equiv_8 3$, and we also need $\ell^3+\ell^2+\ell+1\equiv_8 4$ so we demand $\ell\equiv_4 1$; one can easily check this last condition by hand for $\ell\equiv_8 1,3,5,7$. This finishes part 3.
\end{proof}

Finally, a discussion of the numbers of subgroups of $C_e\times C_{\frac{mn}{e}}$ of a fixed type. Recall that the \emph{Dedekind $\psi$ function} is the multiplicative function $\psi\colon \N\to\N$ defined by taking values $\psi(p^e)=p^e+p^{e-1}$ for primes $p$ and $e\geq 1$, with $\psi(1)=1$. Equivalently, $\psi(n)$ is the degree of the finite \'{e}tale map of stacks $\M_0(n)\to \M_{\Z[\frac{1}{n}]}$ and the number of cyclic subgroups of $C_n\times C_n$ of order $n$.

\begin{prop}\label{beautyofsubgroupnumbers}
Let $d,e,m,n$ be positive integers such that $d^2|m$ and $e^2|mn$, and write $c_{m,n}(d,e)$ for the number of subgroups of $C_e\times C_{\frac{mn}{e}}$ of type $C_d\times C_{\frac{m}{d}}$. Then:
\begin{enumerate}
\item If either $d\nmid e$ or the two conditions $m\nmid de$ and $e\nmid dn$ both hold, then $c_{m,n}(d,e)=0$.
\item If $d|e|dn$ and $m\nmid de$, then $c_{m,n}(d,e)=\frac{e}{d}$.
\item If $d|e$ and $m|de$, then $c_{m,n}(d,e)=\psi(\frac{m}{d^2})$ where $\psi$ denotes the Dedekind $\psi$ function.
\item Regarding $d,e$ as variables, there is the following equality of integral polynomials in $X$:
\[\sum_{\substack{d^2|m \\ e^2|mn}}c_{m,n}(d,e)X^e=\sum_{\substack{b|m,n \\ a^2|\frac{mn}{b^2}}} bX^{ab}\]
\end{enumerate}
\end{prop}

For the following proof, recall the \emph{Euler totient function} is the multiplicative function $\phi\colon \N\to\N$ defined by taking values $\psi(p^e)=p^e-p^{e-1}$ for primes $p$ and $e\geq 1$, with $\psi(1)=1$. Equivalently, $\psi(n)$ is the degree of the finite \'{e}tale map of stacks $\M_1(n)\to \M_0(n)$, the order of the group $(\Z/n\Z)^\times$, or the count of numbers between $1$ and $n$ which are not divisors of $n$.

\begin{proof}
For part 1, notice that if there exists an injection of groups
\[i\colon C_d\times C_{\frac{m}{d}}\to C_e\times C_{\frac{mn}{e}}\]
then either $d|e$ or $d|\frac{mn}{e}$, and either $\frac{m}{d}|e$ or $\frac{m}{d}|\frac{mn}{e}$, as the two factors on the left have to include somewhere on the right. This last pair of conditions immediately reveals that either $m|de$ or $e|dn$, so we will focus our attention on the potential for $d$ to divide $e$ or not. Suppose $d\nmid e$, then the number $y_1$ defined by $i(1,0)=(x_1,y_1)$ must generate the standard cyclic subgroup $C_{\frac{mn}{e}}$ of order $d$, else $i$ would fail to be injective. For similar reason, as $d\nmid e$ then $\frac{m}{d}\nmid e$ and $y_2$ defined by $i(0,1)=(x_2,y_2)$ must generate the standard cyclic subgroup of order $\frac{m}{d}$ in the right-hand factor in the codomain of $i$. As $d|\frac{m}{d}$ by assumption, there is an $a$ such that $ay_1=y_2$. Our assumption that $d\nmid e$ implies that $(e,-ae)$ is nonzero in the domain of $i$, however, from the above considerations we arrive at the equations
\[i(e,-ae)=(e(x_1-ax_2),e(y_1-ay_2))=(0,0)\]
inside the codomain of $i$, in contradiction to the injectivity of $i$, hence $d|e$.\\

For parts 2 and 3, we appeal to the formula for $c_{m,n}(d,e)$ found in \cite[Th.4.5]{subgroupsofproduct}
\[c_{m,n}(d,e)=\sum_{\substack{i|e,j|\frac{mn}{e} \\ m|ij \\ \lcm(i,j)+\frac{m}{d}}}\phi\left(\frac{ij}{m}\right)\]
ranging over all permitted positive integers $i$ and $j$, where $\phi$ is Euler's totient function. This totient function as well as the functions occurring in parts 2 and 3 are multiplicative with respect to prime decompositions, hence we fix a prime $\ell$ and replace $d,e,m,n$ with $\ell^d, \ell^e, \ell^m$, and $\ell^n$, respectively. This changes our fixed assumptions on $d,e,m$, and $n$ to $2d\leq m$, $2e\leq m+n$, and $d\leq e$. With this change of variables, the formula above reads
\begin{equation}\label{littlecombinatorialexcusion} c_{\ell^m,\ell^n}(\ell^d,\ell^e)=\sum_{\substack{0\leq a\leq e \\ 0\leq b\leq m+n-e \\ m\leq a+b \\ \max(a,b)=m-d}}\phi\left(\ell^{a+b-m}\right).\end{equation}
We will readily split the above sum into two parts depending on if either $a$ or $b$ achieves the maximum $m-d$. Now we restrict to our two cases. For part 2, we further assume that $e<m-d$ which when combined with $a\leq e$ means $a$ cannot achieve the maximum, so $b=m-d$. We also know $m\leq a+b=a+m-d$ which means $d\leq a\leq e$. The formula (\ref{littlecombinatorialexcusion}) then becomes
\[\sum_{d\leq a\leq e} \phi\left(\ell^{a-d}\right)=1+\sum_{d+1\leq a\leq e} \left(\ell^{a-d}-\ell^{a-d-1}\right)=\ell^{e-d}\]
as desired. For part 3 we have the added assumption that $m-d\leq e$. Consider the half of (\ref{littlecombinatorialexcusion}) where the maximum is obtained by $a$ and potentially $a=b$, so we have $a=m-d$ and $b\leq m-d$, and the variable $b$ ranges over $0\leq b\leq m+n-e$. Using the facts that $m\leq a+b=m-d+b$ we see that $d\leq b\leq m-d$. The added assumption for part 3 that $m\leq d+e$ and the continued assumption that $2e\leq m+n$ combine to show that $2e\leq d+e+n$ hence $e\leq d+n$. Adding $m$ to both sides gives $m+e\leq m+n+d$ and hence $m-d\leq m+n-e$. This shows that the range of the variable $b$ is then precisely $d\leq b\leq m-d$. This first half of (\ref{littlecombinatorialexcusion}) for part 3 where $a=m-d$ then simplifies to
\[\sum_{d\leq b\leq m-d}\phi\left(\ell^{b-d}\right)=\ell^{m-2d}\]
exactly as in part 2. Consider now the second half of (\ref{littlecombinatorialexcusion}) for part 3 where $b=m-d$ and $a$ is strictly less than $m-d$, so $a\leq m-d-1$. Copying the same arguments as above, we see that $a$ ranges over the values $d\leq a\leq m-d-1$---our assumption for part 3 that $m-d\leq e$ shows this range is strictly tighter than the condition $a\leq e$ from (\ref{littlecombinatorialexcusion}). The second half of (\ref{littlecombinatorialexcusion}) for part 3 is then given by $\ell^{m-2d-1}$, which produces the desired result for part 3:
\[c_{\ell^m,\ell^n}(\ell^d,\ell^e)=\ell^{m-2d}+\ell^{m-2d-1}=\psi\left(\ell^{m-2d}\right)\]
For part 4, we note that as the coefficients on either side of the desired equality are multiplicative it suffices to prove the equation
\begin{equation}\label{ladicpolynomailequiation}
\sum_{\substack{0\leq 2d\leq m \\ 0\leq 2e \leq m+n}} c^\ell(d,e)X^{\ell^ e} = \sum_{\substack{0\leq t\leq \min(m,n) \\ 0\leq 2u\leq m+n-2t}}\ell^ t X^{\ell^{t+u}}
\end{equation}
holds, where $c^\ell(d,e)=c_{\ell^m,\ell^n}(\ell^d,\ell^e)$ for some prime $\ell$. We claim this equality is clear if one studies the following $\mathbf{N}^2$-table of values of $c^\ell(-,-)$; see \Cref{firstexamplenumbetheoryttable} for two examples tables, and \Cref{firstnumbetheoryttable,secondnumbetheoryttable} at the end of this section for the general cases.\\

In more detail, the restrictions that $0\leq 2d\leq m$ and $0\leq 2e\leq m+n$ show these tables are all supported in a bounded domain of $\mathbf{N}^2$---let us visualise $d$ on the horizontal axis and $e$ on the vertical axis. By part 1, we also know that $c^\ell(d,e)=0$ for $e<d$, so a triangle in the lower-right corner of our domain vanishes. Part 1 also states that $c^\ell(d,e)$ vanishes if both the inequalities $m> d+e$ and $e-d> n$ hold. If $m\leq n$, we note that $2e\leq m+n\leq 2n$ so $e\leq n$, so we see that our domain is already contained in the region where $e-d\leq n$. If $n\leq m$, then there is another triangle, now in the upper-left corner of our domain, which vanishes.\\

One can now arrive at the desired formula (\ref{ladicpolynomailequiation}) as follows:
\begin{itemize}
\item the left-hand side can be obtained by summing each row separately, so for each fixed value of $e$ we sum all $c^\ell(d,e)$ for all possible $d$ as the coefficient of $X^{\ell^e}$. This takes the form
\[1\cdot X+(1+\ell)X^{\ell}+(1+\ell+\ell^2)X^{\ell^2}+\cdots+ (1+\cdots +\ell^8)X^{\ell^8}+(1+\cdots +\ell^8)X^{\ell^9}+(1+\cdots +\ell^8)X^{\ell^{10}}\]
in the $(m,n)=(8,12)$ example in \Cref{firstexamplenumbetheoryttable}.
\item the right-hand side can be obtained by fixing a particular fixed value of $t$, so a particular form of coefficient $\ell^t$ of some $X^{\ell^{t+u}}$, and summing up all of the possible $X^{\ell^{t+u}}$ which can occur. Again, this takes the form
\[1\cdot (X+X^\ell+\cdots+X^{\ell^{10}})+\ell(X^{\ell}+\cdots+X^{\ell^{10}})+\cdots+\ell^8(X^{\ell^{8}}+X^{\ell^{9}}+X^{\ell^{10}})\]
in the $(m,n)=(8,12)$ example in \Cref{firstexamplenumbetheoryttable}. The range of all such possible $t$ is precisely $0\leq t\leq \min(m,n)$ and for a fixed $t$ the range of all possible $u$ is precisely $0\leq 2u\leq m+n-2t$.
\end{itemize}
This yields the desired formula (\ref{ladicpolynomailequiation}) and ends our proof.
\end{proof}

\begin{table}[h!]\begin{center}
\begin{tabular}{|c|c|c|c|c|}
\hline
8+7 & 6+5 & 4+3 & 2+1 & 0 \\ \hline
8+7 & 6+5 & 4+3 & 2+1 & 0 \\ \hline
8+7 & 6+5 & 4+3 & 2+1 & 0 \\ \hline
7   & 6+5 & 4+3 & 2+1 & 0 \\ \hline
6   & 5   & 4+3 & 2+1 & 0 \\ \hline
5   & 4   & 3   & 2+1 & 0 \\ \hline
4   & 3   & 2   & 1   & 0 \\ \hline
3   & 2   & 1   & 0   &   \\ \hline
2   & 1   & 0   &     &   \\ \hline
1   & 0   &     &     &   \\ \hline
0   &     &     &     &   \\ \hline
\end{tabular}\qquad
\begin{tabular}{|c|c|c|c|c|l|l|}
\hline
  &   &   &   & 4+3 & 2+1 & 0 \\ \hline
  &   &   & 4 & 3   & 2+1 & 0 \\ \hline
  &   & 4 & 3 & 2   & 1   & 0 \\ \hline
  & 4 & 3 & 2 & 1   & 0   &   \\ \hline
4 & 3 & 2 & 1 & 0   &     &   \\ \hline
3 & 2 & 1 & 0 &     &     &   \\ \hline
2 & 1 & 0 &   &     &     &   \\ \hline
1 & 0 &   &   &     &     &   \\ \hline
0 &   &   &   &     &     &   \\ \hline
\end{tabular}
\caption{\label{firstexamplenumbetheoryttable} Values of $c^\ell(d,e)$ for $(m,n)=(8,12)$ and $(12,4)$, respectively. An $x$ (resp.\ $x+y$) above refers to $c^\ell(d,e)=\ell^x$ (resp.\ $\ell^x + \ell^y$) where $d$ runs along the horizontal axis and $e$ along the vertical.}
\end{center}\end{table}

\begin{table}[h!]\begin{center}
\begin{tabular}{|c|c|c|c|c|c|c|c|}
\hline
$\scriptscriptstyle m+(m-1)$ & $\scriptscriptstyle(m-2)+(m-3)$ & $\scriptscriptstyle(m-4)-(m-5)$ & $\cdots$  & $\scriptscriptstyle 6+5$ & $\scriptscriptstyle 4+3$ & $\scriptscriptstyle 2+1$ & $\scriptscriptstyle0$ \\ \hline
    $\scriptscriptstyle\vdots$    &      $\scriptscriptstyle\vdots$       &    $\scriptscriptstyle\vdots$         &  &   $\scriptscriptstyle\vdots$  & $\scriptscriptstyle\vdots$    &  $\scriptscriptstyle\vdots$   & $\scriptscriptstyle\vdots$  \\ \hline
$\scriptscriptstyle m+(m-1)$ & $\scriptscriptstyle(m-2)+(m-3)$ & $\scriptscriptstyle(m-4)-(m-5)$ &  &     &     &     &   \\ \hline
$\scriptscriptstyle m+(m-1)$ & $\scriptscriptstyle(m-2)+(m-3)$ & $\scriptscriptstyle(m-4)-(m-5)$ &  &     &     &     &   \\ \hline
$\scriptscriptstyle m-1$     & $\scriptscriptstyle(m-2)+(m-3)$ & $\scriptscriptstyle(m-4)-(m-5)$ &  &     &     &     &   \\ \hline
$\scriptscriptstyle m-2$     & $\scriptscriptstyle m-3$         & $\scriptscriptstyle(m-4)-(m-5)$ &  &     &     &     &   \\ \hline
$\scriptscriptstyle m-3$     & $\scriptscriptstyle m-4$         & $\scriptscriptstyle m-5$         &  &     &     &     &   \\ \hline
  $\scriptscriptstyle\vdots$      &     $\scriptscriptstyle\vdots$        &      $\scriptscriptstyle\vdots$       & $\scriptscriptstyle\ddots$ &   $\scriptscriptstyle\vdots$  &  $\scriptscriptstyle\vdots$   &  $\scriptscriptstyle\vdots$  & $\scriptscriptstyle\vdots$  \\ \hline
        &             &             &  & $\scriptscriptstyle 6+5$ & $\scriptscriptstyle 4+3$ & $\scriptscriptstyle 2+1$ & $\scriptscriptstyle0$ \\ \hline
        &             &             &  & $\scriptscriptstyle5$   & $\scriptscriptstyle4+3$ & $\scriptscriptstyle2+1$ & $\scriptscriptstyle0$ \\ \hline
        &             &             &  & $\scriptscriptstyle4$   & $\scriptscriptstyle3$   & $\scriptscriptstyle2+1$ & $\scriptscriptstyle0$ \\ \hline
        &             &             &  & $\scriptscriptstyle3$   & $\scriptscriptstyle2$   & $\scriptscriptstyle1$   & $\scriptscriptstyle0$ \\ \hline
        &             &             &  & $\scriptscriptstyle2$   & $\scriptscriptstyle1$   & $\scriptscriptstyle0$   &   \\ \hline
        &             &             &  & $\scriptscriptstyle1$   & $\scriptscriptstyle0$   &     &   \\ \hline
        &             &             &  & $\scriptscriptstyle0$   &     &     &   \\ \hline
  $\scriptscriptstyle\vdots$      &     $\scriptscriptstyle\vdots$        &  $\scriptscriptstyle\vdots$           & \reflectbox{$\scriptscriptstyle\ddots$}  &     &     &     &   \\ \hline
$\scriptscriptstyle2$       & $\scriptscriptstyle1$           & $\scriptscriptstyle0$           &  &     &     &     &   \\ \hline
$\scriptscriptstyle1$       & $\scriptscriptstyle0$           &             &  &     &     &     &   \\ \hline
$\scriptscriptstyle0$       &             &             &  &     &     &     &   \\ \hline
\end{tabular}
\caption{\label{firstnumbetheoryttable} Values of $c^\ell(d,e)$ where $m\leq n$, concentrated in $0\leq 2d\leq m$ and $0\leq 2e\leq m+n$. We have also assumed $m$ is even for the above picture, however, the $m$ is odd case simply has $1+0$ in the final column instead of simply $0$'s, and the other columns are shifted appropriately. Each $x$ above corresponds to $c^\ell(d,e)=\ell^x$ and $x+y$ to $c^\ell(d,e)=\ell^x+\ell^y$, where the horizontal axis is the $d$-axis and the vertical axis is the $e$-axis.}
\end{center}\end{table}

\begin{table}[h!]\begin{center}
\begin{tabular}{|c|c|c|c|c|c|c|c|c|c|c|}
\hline
        &         &  &         &         & $\scriptscriptstyle n+(n-1)$ & $\scriptscriptstyle (n-2)+(n-3)$ & $\cdots$  & $\scriptscriptstyle 4+3$ & $\scriptscriptstyle 2+1$ & $\scriptscriptstyle 0$ \\ \hline
        &         &  &         & $\scriptscriptstyle n$     & $\scriptscriptstyle (n-1)$   & $\scriptscriptstyle (n-2)+(n-3)$ &  &   $\scriptscriptstyle\vdots$    &   $\scriptscriptstyle\vdots$    &   $\scriptscriptstyle\vdots$  \\ \hline
        &         &  & $\scriptscriptstyle n$     & $\scriptscriptstyle (n-1)$ & $\scriptscriptstyle (n-2)$   & $\scriptscriptstyle (n-3)$       & $\cdots$ & $\scriptscriptstyle 4+3$ & $\scriptscriptstyle 2+1$ & $\scriptscriptstyle 0$ \\ \hline
        &         & \reflectbox{$\scriptscriptstyle\ddots$} & $\scriptscriptstyle (n-1)$ & $\scriptscriptstyle (n-2)$ & $\scriptscriptstyle (n-3)$   & $\scriptscriptstyle (n-4)$       & $\cdots$ & $\scriptscriptstyle 3$   & $\scriptscriptstyle 2+1$ & $\scriptscriptstyle 0$ \\ \hline
        & $\scriptscriptstyle n$     & \reflectbox{$\scriptscriptstyle\ddots$} & $\scriptscriptstyle (n-2)$ & $\scriptscriptstyle (n-3)$ & $\scriptscriptstyle (n-4)$   &           $\scriptscriptstyle\vdots$    & \reflectbox{$\scriptscriptstyle\ddots$} & $\scriptscriptstyle 2$   & $\scriptscriptstyle 1$   & $\scriptscriptstyle 0$ \\ \hline
$\scriptscriptstyle n$     & $\scriptscriptstyle (n-1)$ &\reflectbox{$\scriptscriptstyle\ddots$}  & $\scriptscriptstyle (n-3)$ & $\scriptscriptstyle (n-4)$ &     $\scriptscriptstyle\vdots$      & $\scriptscriptstyle 2$           & \reflectbox{$\scriptscriptstyle\ddots$} & $\scriptscriptstyle 1$   & $\scriptscriptstyle 0$   &     \\ \hline
$\scriptscriptstyle (n-1)$ & $\scriptscriptstyle (n-2)$ & \reflectbox{$\scriptscriptstyle\ddots$} & $\scriptscriptstyle (n-4)$ &    $\scriptscriptstyle\vdots$     & $\scriptscriptstyle 2$       & $\scriptscriptstyle 1$           &  \reflectbox{$\scriptscriptstyle\ddots$}& $\scriptscriptstyle 0$   &       &     \\ \hline
$\scriptscriptstyle (n-2)$ & $\scriptscriptstyle (n-3)$ & \reflectbox{$\scriptscriptstyle\ddots$} &    $\scriptscriptstyle\vdots$     & $\scriptscriptstyle 2$     & $\scriptscriptstyle 1$       & $\scriptscriptstyle 0$           &  &       &       &     \\ \hline
$\scriptscriptstyle (n-3)$ & $\scriptscriptstyle (n-4)$ &  & $\scriptscriptstyle 2$     & $\scriptscriptstyle 1$     & $\scriptscriptstyle 0$       &               &  &       &       &     \\ \hline
$\scriptscriptstyle (n-4)$ &     $\scriptscriptstyle\vdots$    & \reflectbox{$\scriptscriptstyle\ddots$} & $\scriptscriptstyle 1$     & $\scriptscriptstyle 0$     &           &               &  &       &       &     \\ \hline
  $\scriptscriptstyle\vdots$      & $\scriptscriptstyle 2$     & \reflectbox{$\scriptscriptstyle\ddots$} & $\scriptscriptstyle 0$     &         &           &               &  &       &       &     \\ \hline
$\scriptscriptstyle 2$     & $\scriptscriptstyle 1$     & \reflectbox{$\scriptscriptstyle\ddots$} &         &         &           &               &  &       &       &     \\ \hline
$\scriptscriptstyle 1$     & $\scriptscriptstyle 0$     &  &         &         &           &               &  &       &       &     \\ \hline
$\scriptscriptstyle 0$     &         &  &         &         &           &               &  &       &       &     \\ \hline
\end{tabular}
\caption{\label{secondnumbetheoryttable} Values of $c^\ell(d,e)$ where $n\leq m$, which is concentrated in $0\leq 2d\leq m$ and $0\leq 2e\leq m+n$. As in \Cref{firstnumbetheoryttable}, we have assumed $m$ and $n$ are even---the other cases are similar.}
\end{center}\end{table}

\scriptsize
%old laptop
\bibliography{C:/Users/jackd/Dropbox/Work/references} 

%new mac
%\bibliography{/Users/jackdavies/Dropbox/Work/references} 

\bibliographystyle{alpha}

\end{document}